\documentclass[10 pt,a4paper,leqno]{article}
\usepackage{latexsym,amsmath,amscd,amsfonts}
\usepackage[francais,english]{babel}
\usepackage[all]{xy}

\title{Galois theory of fuchsian $q$-difference equations}

\author{Jacques Sauloy\\
{\small Laboratoire Emile Picard - UMR 5580}\\
{\small Universit\'e Paul Sabatier,118,route de Narbonne,31062,Toulouse Cedex,
France.}\\
{\small E-mail : sauloy@picard.ups-tlse.fr}\\
}

\date{}

\newcommand{\Diag}{\text{diag}}

\newcommand{\Ker}{\text{Ker }}

\newcommand{\Otimes}{\underline{\otimes}}

\def \gt {\symbol{62}}

\def \lt {\symbol{60}}

\begin{document}

\maketitle

\hrule \bigskip

\selectlanguage{francais}

\begin{abstract} 
\emph{
{\small Nous proposons une approche analytique de la th\'eorie 
de Galois des syst\`emes aux $q$-diff\'erences lin\'eaires
singuliers r\'eguliers. Nous combinons la dualit\'e de Tannaka 
avec la m\'ethode de classification de Birkhoff \`a l'aide de
la matrice de connexion pour d\'efinir et d\'ecrire leurs groupes
de Galois. Puis nous d\'ecrivons des \emph{sous-groupes
fondamentaux} qui donnent lieu \`a une correspondance de 
Riemann-Hilbert et \`a un th\'eor\`eme de densit\'e de type
Schlesinger.}
}
\end{abstract}

\selectlanguage{english}

\bigskip \hrule \bigskip

\begin{abstract}
{\small We propose an analytical approach to the Galois theory
of singular regular linear $q$-difference systems. We use Tannaka
duality along with Birkhoff's classification scheme with the
connection matrix to define and describe their Galois groups. 
Then we describe \emph{fundamental subgroups} that give rise 
to a Riemann-Hilbert correspondence and to a density theorem
of Schlesinger's type.}
\end{abstract}

\bigskip \hrule \bigskip



\tableofcontents

\hspace{-0.7cm}
{\footnotesize \textbf{\emph{This preprint contains three appendices
(pp 50-60) that won't appear in the published version.}}}

\bigskip \hrule \bigskip

\begin{quote}

``Je suis convaincu que, tout comme pour les fonctions sp\'eciales solutions
d'\'equations diff\'erentielles, les formules int\'eressantes d\'erivent de
consid\'erations ``g\'eom\'etriques'' simples'' 
(Jean-Pierre Ramis,\cite{RamisJPRTraum}).
\end{quote}


\setcounter{section}{-1}

\setcounter{equation}{-1}

\section{Introduction}


\subsection{Rational linear $q$-difference systems and rational equivalence}

Let $q$ be a fixed complex number such that $|q| \gt 1$.
Let $\sigma_{q}$ denote the dilatation $z \mapsto q z$
of a complex coordinate $z$ on the Riemann sphere
$\mathbf{S} = \mathbf{P}^{1}\mathbf{C}$ , viewed 
as an operator on functions of $z$. 
In \cite{JSAIF}, we studied the classification 
of \emph{rational linear $q$-difference systems}
\begin{equation}
\sigma_{q} X = A X.
\end{equation}
over $\mathbf{S}$ under \emph{rational equivalence}. 
Here, $A \in GL_{n}(\mathbf{C}(z))$ is a given matrix
and the unknown $X$ is a matrix with $n$ rows and 
(possibly multivalued) holomorphic entries on some
$q$-invariant open subset of $\mathbf{S}$. 
The gauge group $GL_{n}(\mathbf{C}(z))$ operates on
the left on solutions of such systems, hence on the
systems themselves. This gives rise to the rational
equivalence relation:
\begin{equation*}
A \sim \left(\sigma_{q} F\right)^{-1} A F \;,\; 
F \in GL_{n}(\mathbf{C}(z)).
\end{equation*}
In this paper, we study the Galois theory of such systems.


\subsection{The classification theorem of Birkhoff}

In analogy with Riemann's classification scheme for
complex fuchsian differential equations, by local data
at the singularities and monodromy transformations
arising from the analytic continuation of local
solutions (the so-called ``connection formulae''),
Birkhoff defined in \cite{Birkhoff1} the ``generalized
Riemann problem'' for $q$-difference equations. He
solved it for \emph{fuchsian systems} under innocuous
(generically true) assumptions. We now describe his
classification scheme, as slightly revised by us in
\cite{JSAIF}. \\

We assume the system (0) to be \emph{fuchsian over 
$\mathbf{S}$}. This condition is precisely defined in 1.2.1;
it is essentially equivalent to the existence 
of fundamental solutions with moderate growth at $0$ 
and $\infty$, as opposed to theta-like growth, 
like $z^{\log z}$ (\emph{see} \cite{RamisGrowth}). 
The idea of Birkhoff is to use local solutions at the
two only $q$-invariant points of $\mathbf{S}$,
$0$ and $\infty$, and to classify fuchsian systems 
using linear data at $0$ and $\infty$ and \emph{one}
matrix connecting a local solution at $0$ to one
at $\infty$. Here are the main steps.
Note the close similitude of steps 1 to 3 with the
classical Frobenius-Fuchs method for ordinary
differential equations.

\begin{enumerate}
\item{Any fuchsian system is meromorphically equivalent 
(near $0$) to one with constant coefficients. Here,
meromorphic equivalence is defined by letting the gauge
transformation $F \in GL_{n}(\mathbf{C}(\{z\}))$
($F$ will then automatically belong to
$GL_{n}(\mathcal{M}(\mathbf{C}))$).}

\item{Any constant coefficient system can be reduced by
linear algebra to $1$-dimensional systems: 
$\sigma_{q} f = c f$ (where $c$ is an exponent) and, 
if its matrix is not semi-simple, $2$-dimensional unipotent
systems: $\sigma_{q} g = g + 1$. One can build solutions 
to such elementary systems that are meromorphic 
on $\mathbf{C}^{*}$ and have moderate growth at $0$
and $\infty$, relying on Jacobi's theta function
$\Theta_{q}$. We thus obtain the \emph{$q$-characters}
$e_{q,c}$ and the \emph{$q$-logarithm} $l_{q}$. From this,
we get a canonical fundamental solution $e_{q,A}$ for each
constant coefficient system $A \in GL_{n}(\mathbf{C})$.}

\item{Each fuchsian system is therefore endowed with 
a local solution at $0$: $X^{(0)} = M^{(0)} e_{q,A^{(0)}}$,
where $M^{(0)} \in GL_{n}(\mathcal{M}(\mathbf{C}))$ and
$A^{(0)} \in GL_{n}(\mathbf{C})$, and, symmetrically, 
a local solution at 
$\infty$: $X^{(\infty)} = M^{(\infty)} e_{q,A^{(\infty)}}$,
The Jordan structures of 
$A^{(0)}$ and $A^{(\infty)}$ are the the required ``linear data''.}

\item{The \emph{connection matrix of Birkhoff} is then defined to be:
$P = \left(X^{(\infty)}\right)^{-1} X^{(0)}$.
It has coefficients in the $q$-invariant subfield
$\mathcal{M}(\mathbf{C}^{*})^{\sigma_{q}}$ 
of $\mathcal{M}(\mathbf{C}^{*})$, that is, in the field 
of meromorphic functions on the elliptic curve 
$\mathbf{E}_{q} =  \mathbf{C}^{*}/q^{\mathbf{Z}}$.}

\end{enumerate}

Now, Birkhoff's theorem says that, given linear data 
at $0$ and at $\infty$ and an invertible elliptic matrix,
one can recover a system (0) well defined up to rational
equivalence. A precise formulation can be found in
\cite {JSAIF}. \\

Birkhoff's method seems to lend itself easily to a
reformulation in the style of Riemann-Hilbert's
classification scheme, via local systems and 
representations of the fundamental group. The
local linear data should allow one to compute
local monodromy transformations around $0$ and
$\infty$, so that a \emph{groupoid with two base points}
naturally appears. Then these base points should be
connected by paths from $0$ to $\infty$: here, the
representation has a natural counterpart, the connection
matrix. But the latter has elliptic coefficients,
so that one does not end up with algebraic matrices
over the complex numbers, as in the classical case:
the field of ``constants'' of the theory, 
$\mathcal{M}(\mathbf{C}^{*})^{\sigma_{q}} =
\mathcal{M}(\mathbf{E}_{q})$, is too big.


\subsection{Adapting Picard-Vessiot theory}

The first break through the difficulties caused by the 
big constant field $\mathcal{M}(\mathbf{E}_{q})$ 
was made by Etingof in \cite{Etingof}.
He adapted Picard-Vessiot theory to the case
of \emph{regular} $q$-difference systems, 
those such that $A(0) = A(\infty) = I_{n}$. 
Etingof defines and builds Picard-Vessiot
extensions and shows the related $q$-difference Galois group
to be generated by the \emph{values} (at authorized points)
$P(a)^{-1} P(b) \in GL_{n}(\mathbf{C})$, where $P$ is the
connection matrix. However, in the case of non regular systems,
elliptic functions spontaneously arise in yet
another way: if $c,d \in \mathbf{C}^{*}$, then, along
with the ``legal'' solution $e_{q,cd}$ to the equation
$\sigma_{q} f = cd f$, there is also $e_{q,c} e_{q,d}$,
so that $\phi(c,d) = \frac{e_{q,c} e_{q,d}}{e_{q,cd}}$
is elliptic. Any field containing these basic solutions
will contain \emph{all} elliptic functions, and it is easy to prove
that there is no way to trivialize the cocycle $\phi(c,d)$
while using ``true functions'', so that there seems to be
no hope for a Picard-Vessiot theory with a group defined
over $\mathbf{C}$ even in the fuchsian case. 
The second breakthrough was accomplished by van der Put and 
Singer in \cite{SVdP}, using \emph{symbolic solutions}.
They build a Picard-Vessiot theory with constant field $\mathbf{C}$. 
They are then able to solve the problem in
total generality, including the case of \emph{irregular}
(non fuchsian) systems.  \\

However, our work is part of a program 
that requires a function theoretic attack at these problems.
First, there is the link of $q$-difference
equations to ordinary differential equations through
$q$-analogies. When $q \to 1$, $q$-analogs of special
functions (like Heine's \emph{basic hypergeometric series},
\emph{see} the ``bible'' \cite{GR}) ``tend to'' their classical 
counterpart. One may wish to follow Galois groups along such a 
\emph{confluence}. Some results in this direction were expounded
in \cite{JSAIF}. They receive here substantial extensions.
This is of obvious interest in mathematical physics,
with the present ubiquity of $q$-deformations. 
Second, the appearance (and central importance) of Jacobi's theta
functions, elliptic functions and complex elliptic curves
in the landscape unveil rich geometric structures.
We build explicit galoisian automorphisms and give them
a geometric interpretation, which allows us to exhibit a reasonable 
candidate for the role of \emph{fundamental group}, that is, a finitely
generated and finitely presented Zariski-dense subgroup of
the Galois group. This is in analogy with Schlesinger's
theorem (\emph{see} \cite{Bertrand}, \cite{CanoRamis}) 
and with the topological flavour of the 
classical Riemann-Hilbert correspondence, where fuchsian differential
equations are classified by monodromy representations
(\emph{see} \cite{Deligne})
\footnote{To these arguments, one should add that
the oldest historical motivation for the $q$-world lies in
magic identities by Gauss, Euler, Jacobi, Ramanujan ...
(\emph{see} \cite{Ramanujan}). These involve classical 
analytical functions and one may hope for a \emph{geometric}
understanding of them. One must also mention that 
$q$-difference equations are a possible intermediate
step to understand the mysterious analogy between 
irregular linear ordinary differential equations and wildly 
ramified phenomena in positive characteristic 
(\emph{see} \cite{RamisGalDiff}).}. \\

As noted before, the use of ``true functions'' as $q$-characters 
forces on us a big constant field $\mathcal{M}(\mathbf{E}_{q})$. 
Yet, in our version of the 
classification theorem of Birkhoff, automorphisms of fuchsian
objects are classified by complex matrices. 
The root of this fact is that we authorize as unique ``legal
model'' for the equation $\sigma_{q} f = c f$ \emph{the one
function $e_{q,c}$ and nobody else}, thereby rigidifying 
a lot the situation. The corresponding drawback is that we
cannot multiply solutions: $e_{q,c}e_{q,d}$ is not legal,
only $e_{q,cd}$ is ; we do not even have an algebra of 
solutions.  
Therefore, to produce a Galois group, we turn to Tannaka duality
(\emph{see} \cite{DM}, \cite{DF}).
This has already been used in this context twice:
by van der Put and Singer in \cite{SVdP} for one, then by
Yves Andr\'e in his work \cite{YA2}, where deformation results
are proved for difference and differential Galois groups. 
\emph{The goal of this paper is therefore to give a tannakian 
formulation of the classification theorem of Birkhoff, while 
using as basic objects uniform analytic functions.} 
We now list our main results (they are detailed in section 0.4).

The \emph{local} category $\mathcal{E}_{f}^{(0)}$ of 
$q$-difference systems is naturally equivalent to                    
the $\mathbf{C}$-linear neutral tannakian category 
$Fib_{p}(\mathbf{E}_{q})$ of flat vector 
bundles over the elliptic curve $\mathbf{E}_{q}$. There is a naturally 
defined local Galois groupoid $G^{(0)}$ of $\mathcal{E}_{f}^{(0)}$
with base set $\mathbf{C}^{*}$ and we compute it explicitly, as well 
as the local Galois group, also called $G^{(0)}$.
We build explicit elements of the group $G^{(0)}$, and we want to see them 
as \emph{loops}
\footnote{By nicknaming ``loop'' a galoisian automorphism
(i.e. a tensor automorphism of a fibre functor), we just go one 
little step beyond the terminology introduced by Katz in
\cite{Katz}, 1.1.2.1.}
;then, we single out two commuting loops with a nice topological 
interpretation as ``fundamental loops of an infinitesimal elliptic 
curve''. The group they generate is Zariski-dense in $G^{(0)}$.
We see it as the \emph{local fundamental group}. \\

The \emph{global} category $\mathcal{E}_{f}$ of $q$-difference 
systems is  equivalent to the $\mathbf{C}$-linear neutral tannakian 
category of triples $(A^{(0)},M,A^{(\infty)})$, made up of 
two flat vector bundles and a meromorphic isomorphism between. 
Evaluating such isomorphisms at non singular points provides us
with galoisian isomorphisms (``paths'') in the global Galois 
groupoid $G$. The local groupoids $G^{(0)}$ and $G^{(\infty)}$ together 
with these paths generate a Zariski-dense subgroupoid of $G$.
For regular abelian objects with prescribed singular locus,
the Galois group is 
Zariski-generated by the values of the connection matrix. Using methods 
from geometric class field theory, one can classify all regular abelian 
representations with prescribed singular locus $S$ of the global Galois 
group as representations of an explicit affine group. \\

Last, we describe explicitly the confluence of the generators
of the global Galois group to elements of a differential Galois
group when $q$ tends to $1$.


\subsection{Contents of this paper}

Let us now describe more precisely the organisation of
this paper. In chapter 1, we review some basic properties
of linear $q$-difference systems with rational coefficients. 
In section 1.1, we briefly recall general
algebraic properties of $q$-difference systems and $q$-difference
modules, mostly adapted from \cite{SVdP}.
In section 1.2, we define the category $\mathcal{E}_{f}$ 
of \emph{fuchsian} $q$-difference systems, 
a neutral tannakian category over $\mathbf{C}$, and we summarize the 
first part our previous work \cite{JSAIF}, about local solutions and 
classification of such systems. \\

Chapters 2 and 3 contain the core of this paper, the construction
and tentative description of the local and global Galois groups
and groupoids of the category of fuchsian $q$-difference systems.
Chapter 2 deals with the local setting and chapter 3 with the glueing 
of the local descriptions at $0$ and $\infty$. 
In 2.1, we study the \emph{local} category $\mathcal{E}_{f}^{(0)}$; 
we consider here as localisation at $0$ 
the action to allow for morphisms with coefficients defined 
\emph{locally for the transcendant topology}. We find a particularly 
simple equivalent category $\mathcal{P}^{(0)}$ of local models related 
with the category $\mathcal{R}$ of complex representations of 
$\mathbf{Z}$ . In section 2.2, we exploit this link to exhibit 
a $\mathbf{C}^{*}$-indexed family of fibre functors 
$\omega^{(0)}_{z_{0}}$ extending the canonical fibre functor $\omega$ 
on $\mathcal{R}$, allowing us to compute the \emph{local Galois 
groupoid $G^{(0)}$ of $\mathcal{E}_{f}^{(0)}$}, with base set 
$\mathbf{C}^{*}$ from the knowledge of the proalgebraic hull 
$\mathbf{Z}^{alg} = 
Hom_{grp}(\mathbf{C}^{*},\mathbf{C}^{*}) \times \mathbf{C}$ 
of $\mathbf{Z}$: $G^{(0)}$ is thereby identified with a subgroupoid 
of $\mathbf{Z}^{alg} = Aut^{\otimes}(\omega)$. We prove: \\

\textbf{2.2.2.1 Theorem: the local Galois groupoid. -}
\emph{With the previous identification of $Aut^{\otimes}(\omega)$
with $\mathbf{Z}^{alg}$,
$$
Iso^{\otimes}(\omega^{(0)}_{z_{0}},\omega^{(0)}_{z_{1}}) = 
\{(\gamma,\lambda) \in \mathbf{Z}^{alg} \;/\; \gamma(q)z_{0} = z_{1} \}.
$$}
The Galois group (also called $G^{(0)}$) is then immediately deduced
in \textbf{corollary 2.2.2.2}. Two commuting algebraically independent
elements (``loops'') $\gamma_{1}$ and $\gamma_{2}$ in the semi-simple 
component of $G^{(0)}$ are built in 2.2.3. The following density
theorem thus provides an analog to the local fundamental group: \\

\textbf{2.2.3.5 Theorem. -}
\emph{The subgroup of
$Hom_{grp}(\mathbf{C^{*}},\mathbf{C^{*}}) \times \mathbf{C}$
whose unipotent component is $\mathbf{Z} \subset \mathbf{C}$
and whose semi-simple component is generated by $\gamma_{1}$
and $\gamma_{2}$ is Zariski-dense in the local Galois group.} \\

According to Weil's correspondence between the degree $0$ vector bundles
on a compact Riemann surface and the representations of its fundamental
group, our category $\mathcal{P}^{(0)}$ of local models is shown to be
equivalent to the category of flat vector bundles over the elliptic curve 
$\mathbf{E}_{q}$ in 2.3. Our solutions can be interpreted as sections 
of these bundles and the singling out of our fundamental solutions is 
equivalent to a choice of frames. \\

In chapter 3, we start global Galois theory. In 3.1, we define a 
category $\mathcal{C}$ of \emph{connection triples} 
$(A^{(0)},M,A^{(\infty)})$, made up of two local (flat) systems and 
an isomorphism $M$ between them that is meromorphic over $\mathbf{C^{*}}$.
In \textbf{proposition 3.1.1.3},
we prove that $\mathcal{E}_{f}$ and $\mathcal{C}$ are equivalent
tensor categories. The natural projections to the local categories
at $0$ and $\infty$ equip them with two $\mathbf{C}^{*}$-indexed 
families of fibre functors $\omega^{(0)}_{z_{0}}$ and
$\omega^{(\infty)}_{z_{0}}$. This defines a Galois groupoid $G$ of 
$\mathcal{C}$ with base set $\mathbf{C}^{*} \amalg \mathbf{C}^{*}$.
Evaluating $M$ at a point $z_{0}$ defines a galoisian isomorphism 
$\Gamma_{z_{0}}$ between the restrictions of $\omega^{(0)}_{z_{0}}$ 
and $\omega^{(\infty)}_{z_{0}}$ to the tannakian 
subcategory $\mathcal{C}_{\Sigma}$ of systems with singular locus 
carried by $\Sigma$; we want to see such an element as a ``path''.
We then get another density result: \\

\textbf{3.1.2.3 Theorem. -}
\emph{The local groupoids at $0$ and at $\infty$ 
(defined and computed in chapter 2) together with the
paths $\Gamma_{z_{0}} \;,\; z_{0} \not\in \Sigma$
generate a Zariski-dense subgroupoid of the Galois groupoid
of $\mathcal{C}_{\Sigma}$.} \\

In 3.2, we follow more literally Birkhoff and get stuck in many 
complications due to the bad multiplicative properties of solutions, 
precisely, the fact that $e_{q,cd} \not= e_{q,c} e_{q,d}$, leading us 
to a \emph{twisted tensor structure} and a \emph{twisted connection 
matrix}. However, the grubby computations of 3.2 give a more concrete 
approach and a simple structural description of the global Galois group. 
Moreover, it is better fitted for the important confluence results of 
chapter 4. The relation with the point of view of 3.1 is explained in 
3.2.3. In both approaches, we have exhibited a lot of ``connecting'' 
galoisian isomorphisms (from $0$ to $\infty$), built from the values of 
the connection matrix, and we have proven a density lemma; but we want
to reduce the \emph{uncountable} family of generators thus obtained
and to make explicit the relations between them. We solve this difficult 
problem in 3.3 for \emph{regular abelian} objects.
The Galois group is then reduced to its connection component, which, 
after 3.1.2, is Zariski-generated by the values of the connection 
matrix. Relative to a prescribed singular locus $S$ we explicitly
define and compute in 3.3.2 (equation (3), 3.3.2.1 and 3.3.2.2 )
affine algebraic groups $L_{S,s}$, $L_{S,s}'$ and $L_{S,u}$ and prove: \\

\textbf{3.3.2.3 Theorem. -}
\emph{The abelian regular objects with singularities in $S$ are 
classified by the representations of the following algebraic group:
$$
\pi_{ab,S,reg}^{1} = \frac{L_{S,s}}{L_{S,s}'} \times L_{S,u}.
$$
}

In chapter 4, we study, 
along the lines of our previous work \cite{JSAIF}, the confluence 
of $q$-difference galoisian automorphisms to differential galoisian 
automorphisms when $q \to 1$. This can be seen as an ``internal'', 
maybe more explicit, illustration of results by 
Yves Andr\'e in \cite{YA2}, relating a family of $q$-difference Galois 
groups to a differential Galois group.


\subsection*{General facts and conventions}

We fix for the whole paper a complex number $q \in \mathbf{C}$
such that $|q| > 1$ and a number $\tau \in \mathcal{H}$ 
(Poincar\'e's half plane) such that $q = e^{- 2 \imath \pi \tau}$.
The only exception is section 4.1, where $q$ and $\tau$ 
will be allowed to vary. For any $c \in \mathbf{C}^{*}$, there is 
a unique pair $(m,d) \in \mathbf{Z} \times \mathbf{C}$ with
$c = q^{m} d$, where $d$ belongs to the \emph{fundamental annulus}: 
$1 \leq |d| < |q|$; we then put $\epsilon(c) = m$ and 
$\overline{c} = d$, so that $\epsilon(c)$ is the integral part of
$\frac{\log |c|}{\log |q|}$ and $c = q^{\epsilon(c)} \overline{c}$. \\

For any complex regular matrix $A \in GL_{n}(\mathbf{C})$, we write 
$A = A_{s} A_{u}$ its \emph{(multiplicative) Dunford decomposition}:
$A_{s}$ is semi-simple, $A_{u}$ is unipotent and they commute.
Such a decomposition is unique and both factors are polynomials
in $A$. Let $f$ be any map: $\mathbf{C}^{*} \rightarrow \mathbf{C}^{*}$.
Write $A_{s} = Q \; \Diag(c_{1},\ldots,c_{n}) \; Q^{-1}$. Then,
the matrix
$Q \; \Diag(f(c_{1}),\ldots,f(c_{n})) \; Q^{-1}$ depends on $A_{s}$
only and we write it $f(A_{s})$. Except otherwise explicitly stated,
we shall then write $f(A) = f(A_{s})$. One exception to this last 
convention is that $\overline{A} = \overline{A_{s}} A_{u}$, so that 
$A = q^{\epsilon(A)} \overline{A}$. Another exception appears when
we define $e_{q,A}$ in 1.2.2. 
Also, note the following general fact: if $S A = B S$, then,
for any map $f:\mathbf{C}^{*} \rightarrow \mathbf{C}^{*}$
and any $\lambda \in \mathbf{C}$, one has $S f(A_{s}) = f(B_{s}) S$ 
and $S A_{u}^{\lambda} = B_{u}^{\lambda} S$; here, of course,
$A_{u}^{\lambda} = 
\underset{k \geq 0}{\sum} {\begin{pmatrix} \lambda \\ k \end{pmatrix}} 
\left(A_{u} - I_{n}\right)^{k}$ (actually, a finite sum). 

{\small
\subsection*{Acknowledgements}

Part of the present work constitutes the second part of my thesis, 
that was prepared under the direction of Jean-Pierre Ramis. 
But the whole investigation has been deeply influenced by his vision.
I thank him heartily for his precious support. I have also been greatly 
helped by suggestions of the referees of that thesis, Yves Andr\'e 
and Marius van der Put. I am grateful for many fruitful discussions 
with Bertrand Toen, Joseph Tapia, Lucia Di Vizio and Changgui Zhang. 
Last, but not least, I owe thanks to the referee of this paper, whose 
outstanding patience allowed to increase dramatically the precision 
and understandability of the text. \\

A first version of this work was circulated in autumn 2000, with
substantially the same contents. 
Then, in june 2001, Bernard Malgrange gave me a copy of the paper
\cite{BG} of 1996, by Baranovsky and Ginzburg, emphasizing the use of
vector bundles over an elliptic curve to classify $q$-difference 
systems. The point of view is strikingly similar, there
are however important differences (see section 2.3.4).
}



\section{Preliminary results}


\subsection{Difference systems and difference modules}

Most of the general formalism here is expounded in \cite{SVdP}.
\footnote{A different formalism is presented in \cite{YA2}, 
which introduces a notion of \emph{non commutative connection}.}.
Let $(K,\sigma)$ be a difference field: $K$ is a field 
and $\sigma$ is an automorphism of $K$. We shall also, 
without further notice, denote by $\sigma$ the canonical 
extensions to the vector spaces of matrices, or of row or 
column vectors. The \emph{$q$-difference equation of order $n$}:
\begin{equation}
\sigma^{n} f + a_{1} \sigma^{n-1} f + \cdots + a_{n} f = 0, \quad
a_{1},\ldots,a_{n} \in K \;,\; a_{n} \not=0
\end{equation}
can be put into system form as a \emph{$q$-difference system of 
rank $n$}:
\begin{equation}
\sigma X = A X, 
\quad A \in GL_{n}(K).
\end{equation}
Conversely, any such system is equivalent to such an equation
via the gauge equivalence defined by:
$A \sim  \left(\sigma F \right)^{-1} A F, \quad F \in GL_{n}(K)$.
This is a consequence of \emph{Birkhoff's cyclic vector lemma}
(\emph{see} \cite{JSAIF}, appendix B or \cite{LDVpreprint}). 
As a consequence, from now on, we won't distinguish between
equations and systems.


\subsubsection*{1.1.1~~~The category of difference modules}

The system (2) can in turn be modelled more intrinsically as
a \emph{difference module} $(K^{n},\Phi)$ by putting
$\Phi: K^{n} \rightarrow K^{n} \;\; X \mapsto A^{-1} \sigma X$,
where a difference module over the difference field $K$ (more properly,
over $(K,\sigma)$) is a finite dimensional $K$-vector space $M$
equipped with a $\sigma$-linear automorphism $\Phi_{M}$ (that is,
a group automorphism such that $\Phi_{M}(x m) = \sigma(x) \Phi_{M}(m)$).
Then $\Phi_{M}$ is actually linear over the \emph{constant subfield}:
$$
C_{K} = K^{\sigma} = \{x \in K \;/\; \sigma(x) = x \}.
$$
A morphism $f: (M, \Phi_{M}) \rightarrow (N, \Phi_{N})$ is a 
$K$-linear map such that $\Phi_{N} \circ f = f \circ \Phi_{M}$.
We shall usually write $M$, $f: M \rightarrow N$, etc ...,
the difference module structure being implicit. 
Also, we write $r(M)$ for the \emph{rank} of the difference
module $M$, that is, its dimension as a $K$-vector space.
We thus obtain the category $DiffMod(K,\sigma)$ of difference modules 
over the difference field $(K,\sigma)$. According to \cite{SVdP}, 
this is a $C_{K}$-linear rigid abelian tensor category 
(\emph{see} \cite{DM} and \cite{DF}). Clearly, forgetting the
difference structure (i.e. the automorphism $\Phi_{M}$)
provides us with a fibre functor from $DiffMod(K,\sigma)$
to the category $Vect_{K}^{f}$ of finite dimensional
$K$-vector spaces, thus making $DiffMod(K,\sigma)$ a
$C_{K}$-linear tannakian category neutralized by $K$ 
(\emph{see} \cite{DM}).


\subsubsection*{1.1.2~~~The category of difference systems}

Choosing an ordered basis for each finite dimensional vector
space over $K$ allows one to replace the category $Vect_{K}^{f}$
by its essential full subcategory with objects the 
$K^{n}$ ($n \in \mathbf{N}$); then the morphisms 
$K^{n} \rightarrow K^{p}$ can be identified with the matrices
in $M_{p,n}(K)$. This subcategory is equivalent, as an abelian
category, to $Vect_{K}^{f}$. To have an equivalence \emph{qua} tensor
categories, it is enough to consistently choose an order on
the product of any two ordered bases. For instance, choosing
the lexicographic order gives bijections:
$$
\begin{cases}
\{1,\ldots,n_{1}\} \times \{1,\ldots,n_{2}\}
\rightarrow \{1,\ldots,n_{1}n_{2}\} \\
(i_{1},i_{2}) \mapsto i_{1} + n_{1} (i_{2} - 1)
\end{cases}
.$$
We thus obtain well defined isomorphisms:
$K^{n_{1}} \underset{K}{\otimes} K^{n_{2}} \rightarrow K^{n_{1}n_{2}}$
and
$M_{p_{1},n_{1}}(K) \underset{K}{\otimes} M_{p_{2},n_{2}}(K)
\rightarrow M_{p_{1}p_{2},n_{1}n_{2}}(K)$.
The resulting tensor category has trivial (i.e. identity)
associativity and unity constraints (but it is not so for
the commutativity constraint). \\

In the same spirit, define the \emph{category $DiffEq(K,\sigma)$ of 
difference equations over the difference field $(K,\sigma)$}:
it has as objects the pairs $(K^{n},A)$ where $n \in \mathbf{N}$ 
and $A \in Gl_{n}(K)$; and, as morphisms from $(K^{n},A)$ to 
$(K^{p},B)$, the matrices $F \in M_{p,n}(K)$ 
such that $(\sigma F) A = B F$ (the composition is the natural one).
We shall often simply denote by $A$ the object $(K^{n},A)$ and
identify it with the difference equation $\sigma X = A X$;
the main reason to make the base space $K^{n}$ explicit is to
give a more natural notation to the forgetful functor 
$(K^{n},A) \leadsto K^{n}$.
To obtain $DiffEq(K,\sigma)$ as a \emph{tensor}
model of $DiffMod(K,\sigma)$, we define the tensor
product of two objects by:
$(K^{n_{1}},A_{1}) \otimes (K^{n_{2}},A_{2}) =
(K^{n_{1} n_{2}},A_{1} \otimes A_{2})$,
with the previous identification of $A_{1} \otimes A_{2}$ to
a matrix in $M_{n_{1} n_{2}}(K)$; and the tensor product of 
two morphisms 
$F_{i}: (K^{n_{i}},A_{i}) \rightarrow (K^{p_{i}},B_{i})$
($i = 1,2$) as $F_{1} \otimes F_{2}$, similarly identified with 
a matrix in $M_{p_{1} p_{2},n_{1} n_{2}}(K)$. 
>From 1.1.1, we draw that
\emph{the above constructions make $DiffEq(K,\sigma)$ into 
a rigid $C_{K}$-linear abelian tensor category equivalent to
$DiffMod(K,\sigma)$. It is tannakian and neutralized by $K$.}
The basic relevant linear and tensor constructions are detailed 
in \cite{JSthese} and in \cite{JSAIF2}. In particular, the unit
$\underline{1}$ is $(K,1)$.


\subsubsection*{1.1.3~~~Functors of solutions}

The \emph{functor of global sections} on $DiffMod(K,\sigma)$
is the functor $\Gamma = Hom(\underline{1},-)$.
The elements of $\Gamma(M)$ are precisely the fixed vectors 
of $\Phi_{M}$ in $M$. We clearly get a left exact functor to 
the category of $C_{K}$-vector spaces. From the inequality 
$\dim_{C_{K}} \Gamma(M) \leq r(M)$ (which follows from the 
``$q$-analogue of the Wronskian lemma'', \emph{see} 
\cite{LDVpreprint}, I.1.2), follows that this functor actually goes
to the category $Vect_{C_{K}}^{f}$.  
Now let $(K',\sigma')$ be an extension of $(K,\sigma)$,
that is, $K'$ is an extension of $K$ and 
$\sigma'_{|K} = \sigma$. The naturally defined base change 
functor from $DiffMod(K,\sigma)$ to $DiffMod(K',\sigma')$
is exact and $\otimes$-preserving. Combining these constructions 
yields a functor $M \leadsto (M \otimes K')^{\sigma'}$. \\

In matrix terms,
we associate to a system (2) the $C_{K'}$-space $S_{K'}(A)$
of solutions in ${K'}^{n}$.
We thereby obtain a functor from $DiffEq(K,\sigma)$ to 
$Vect_{C_{K'}}^{f}$
defined by:
$$
\begin{cases}
A \leadsto S_{K'}(A) \\
(F: A \rightarrow B) \leadsto 
\left(U \mapsto F U: S_{K'}(A) \rightarrow S_{K'}(B)\right)
\end{cases}
$$
Call \emph{fundamental (matrix) solution} of the system with
matrix $A \in GL_{n}(K)$ over the extension $K'$ 
a matrix solution $X \in Gl_{n}(K')$.
The rank of the $C_{K'}$ vector space $S_{K'}(A)$
is exactly $n$ if and only if there is a fundamental solution.
If all systems have a fundamental solution in $K'$,
then the functor of solutions is a fibre functor. 
However, in general, functors
of solutions are neither right exact, nor faithful, nor
$\otimes$-compatible.


\subsection{Fuchsian equations}


\subsubsection*{1.2.1~~~The category $\mathcal{E}_{f}$ of fuchsian equations}

We shall define here the category $\mathcal{E}$ of linear
$q$-difference equations with rational coefficients and
its subcategory $\mathcal{E}_{f}$ of fuchsian equations. 
We shall use the following fields of functions: $\mathbf{C}(z)$, 
the field of rational functions; $\mathcal{M}(\mathbf{C})$, 
the field of meromorphic functions over $\mathbf{C}$;
$\mathcal{M}(\mathbf{C}_{\infty})$, the field of meromorphic 
functions over $\mathbf{C}_{\infty} = \mathbf{S} - \{0\}$; 
and $\mathcal{M}(\mathbf{C}^{*})$, the field of meromorphic 
functions over $\mathbf{C}^{*}$. Each of these function fields, 
endowed with the automorphism $\sigma_{q}: f(z) \mapsto f(qz)$, 
is a difference field. To any of them, we can specialize the 
preceding constructions. \\

We are particularly interested in the category of 
(linear) rational $q$-difference equations, obtained by taking 
$K = \mathbf{C}(z)$ and $\sigma = \sigma_{q}$.
We shall call it $\mathcal{E} = DiffEq(\mathbf{C}(z),\sigma_{q})$.
Since the constant field is, in this case,
$C_{K} = \mathbf{C}(z)^{\sigma_{q}} = \mathbf{C}$,
$\mathcal{E}$ is a $\mathbf{C}$-linear tannakian category
neutralised by $\mathbf{C}(z)$. \\

We shall say that a system with matrix $A \in Gl_{n}(\mathbf{C}(z))$ 
is \emph{strictly fuchsian at $0$} if $A(0) \in Gl_{n}(\mathbf{C})$. 
We shall then call \emph{fuchsian at $0$} a system that is 
meromorphically (that is, through a gauge transformation 
with coefficients in $\mathcal{M}(\mathbf{C})$) equivalent to 
a strictly fuchsian one. It was proved in \cite{JSAIF}, Annexe B,
that this definition is equivalent to the the classical one
(using the Newton polygon). 
Considering $A(\infty)$ and gauge transformations with 
coefficients in $\mathcal{M}(\mathbf{C}_{\infty})$, we similarly define 
systems fuchsian (resp. strictly fuchsian) at $\infty$.
It was also proved in \emph{loc. cit.} that an equation 
fuchsian at $0$ and at the same time fuchsian at $\infty$
is \emph{rationally} (that is, through a gauge transformation 
with coefficients in $\mathbf{C}(z)$) equivalent to one that is
strictly fuchsian over $\mathbf{S}$ (i.e. at $0$ \emph{and} $\infty$). 
Such equations we call \emph{fuchsian over $\mathbf{S}$},
or merely fuchsian. They form a strictly full subcategory 
$\mathcal{E}_{f}$ of $\mathcal{E}$. \\

\textbf{1.2.1.1 Theorem. -} 
\emph{The category $\mathcal{E}_{f}$ is a tannakian subcategory 
of $\mathcal{E}$ over $\mathbf{C}$.} \\

To see that $\mathcal{E}_{f}$ is closed under tensor operations
(including unit, dual and internal $Hom$) it is plainly enough
to consider the case of \emph{strictly fuchsian} objects, 
and then it is obvious. 
Now, from the lemma herebelow, follows that the
kernel in $\mathcal{E}$ of any morphism between fuchsian
objects is itself fuchsian. Therefore, it is a kernel in
$\mathcal{E}_{f}$. Since duality in the tannakian category 
$\mathcal{E}$ exchanges kernels with cokernels (this follows 
from \cite{DM}, p.112), we conclude that $\mathcal{E}_{f}$ is 
indeed an abelian subcategory of $\mathcal{E}$. 
\hfill $\Box$ \\

\textbf{1.2.1.2 Lemma. -} 
\emph{In $\mathcal{E}$, any subobject of an object that is 
fuchsian at $0$ is so.} \\

This is, an immediate consequence of the properties
of the Newton polygon studied in \cite{JSCRAS3},\cite{JSAIF2},
\cite{JSGTQDIF}.
For a more explicit analytic proof
\emph{see} \cite{JSthese} (\emph{see} also \cite{Praagman}). 
\hfill $\Box$


\subsubsection*{1.2.2~~~Local reduction and local solutions}

We recall, here and in the following section,
some results from \cite{JSAIF} and \cite{JSthese}. 
First, define Jacobi's theta function:
$$
\Theta_{q}(z) = \sum_{n\in\mathbf{Z}} q^{- \frac{n(n-1)}{2}} z^{n}.
$$
It is holomorphic over $\mathbf{C}^{*}$ with simple zeroes 
on $q^{\mathbf{Z}}$. It satisfies the $q$-difference equation
$\Theta_{q}(qz) = - qz \Theta_{q}(z)$.
It will be our main brick to build everything.
First, one defines the \emph{$q$-logarithm}:
$$
l_{q}(z) = z \frac{\Theta_{q}'(z)}{\Theta_{q}(z)},
$$
which is meromorphic over $\mathbf{C}^{*}$ with simple poles
on $q^{\mathbf{Z}}$ and satisfies the $q$-difference equation
$l_{q}(qz) = l_{q}(z) + 1$.
Then, for each $c \in \mathbf{C}^{*}$, one defines the
\emph{$q$-character with exponent $c$}. First, if $c$
lies in the \emph{fundamental annulus} 
$\{z \in \mathbf{C} \;/\; 1 \leq |z| \lt |q|\}$, one puts:
$$
e_{q,c}(z) = \frac{\Theta_{q}(z)}{\Theta_{q}(c^{-1} z)}.
$$
For $c$ arbitrary, one writes $c = q^{\epsilon(c)} \overline{c}$,
where $\epsilon(c) \in \mathbf{Z}$ and $\overline{c}$ belongs to the
fundamental annulus, and one puts
\footnote{This is an innocuous modification with respect to \cite{JSAIF}.}:
$$
e_{q,c} = z^{\epsilon(c)} e_{q,\overline{c}}.
$$
Then $e_{q,q^{n}} = z^{n}$ (if $n \in \mathbf{Z}$) and each 
non trivial $e_{q,c}$ is meromorphic over $\mathbf{C}^{*}$ 
with simple zeroes on $q^{\mathbf{Z}}$ and simple poles on 
$c q^{\mathbf{Z}}$. It satisfies the $q$-difference equation
$e_{q,c}(qz) = c e_{q,c}(z)$,
as well as various relations as a family: for instance,
$e_{q,qc}(z) = z e_{q,c}(z)$, etc... \\

Now, let $A \in GL_{n}(\mathbf{C})$ with Dunford decomposition
$A = A_{s} A_{u}$.
If $A_{s} = Q \; \Diag(c_{1},\ldots,c_{n}) \; Q^{-1}$, 
it makes sense to define:
$$
e_{q,A_{s}} = Q \; \Diag(e_{q,c_{1}},\ldots,e_{q,c_{n}}) \; Q^{-1}.
$$
Similarly, defining:
$$
e_{q,A_{u}} = A_{u}^{l_{q}} =
\sum_{k \geq 0} {\begin{pmatrix} l_{q} \\ k \end{pmatrix}} 
\left(A_{u} - I_{n}\right)^{k}
$$
makes sense, since $A_{u}$ is unipotent. One then has
$\sigma_{q} \left(e_{q,A_{s}}\right) =A_{s} e_{q,A_{s}}$ and
$\sigma_{q} \left(e_{q,A_{u}}\right) =A_{u} e_{q,A_{u}}$,
and defining
$$
e_{q,A} = e_{q,A_{s}} e_{q,A_{u}},
$$
one gets the \emph{canonical fundamental solution}
of the constant coefficients system $A$. The above
equality, besides, is a Dunford 
decomposition. From the relation $e_{q,qc}(z) = z e_{q,c}(z)$
stems the equality $e_{q,A} = z^{\epsilon(A)} e_{q,\overline{A}}$. \\

We shall build solutions with coefficients in the field:
$$
\mathbf{K_{0}} = 
\mathcal{M}(\mathbf{C})(l_{q},(e_{q,c})_{c \in \mathbf{C}^{*}}).
$$
As noted in the introduction, the constant subfield 
$\mathbf{K_{0}}^{\sigma_{q}}$ is precisely equal to 
$\mathcal{M}(\mathbf{C}^{*})^{\sigma_{q}}$, the field 
$\mathcal{M}(\mathbf{E}_{q})$ of elliptic functions:
indeed, one inclusion is obvious and the other comes
from the classical fact that the family of the cocycle values:
$$
\phi(c,d) = \frac{e_{q,c} e_{q,d}}{e_{q,cd}} 
          = \frac{\Theta_{q}(z) \Theta_{q}(c^{-1}d^{-1}z)}
                 {\Theta_{q}(c^{-1}z) \Theta_{q}(d^{-1}z)}
$$
generates the group $\mathcal{M}(\mathbf{E}_{q})^{*}$.
Extending the cocycle $\phi$ to semi-simple matrices,
one gets a cocycle $\Phi$ of invertible elliptic
matrices such that, for any two semi-simple matrices
$C,C'$, one has
$e_{q,C} \otimes e_{q,C'} = e_{q,C \otimes C'} \Phi(C,C')$.
On the other hand, it is clear that, for unipotent
matrices $U,U'$, one has
$e_{q,U} \otimes e_{q,U'} = e_{q,U \otimes U'}$.
Thus, for any two invertible matrices $A,A'$, 
$\Phi(A_{s},A'_{s})$ is exactly the defect of
$\otimes$-compatibility of the formation of our
canonical solutions $e_{q,-}$:
$$
e_{q,A} \otimes e_{q,A'} = e_{q,A \otimes A'} \Phi(A_{s},A'_{s}).
$$

To build explicit solutions, we follow closely 
the classical way for differential equations: 
\emph{see} \cite{Ince}, \cite{Wasow}. We consider local
reduction at $0$, the case of $\infty$ being similar. \\

First, any fuchsian system reduces by definition to a strictly 
fuchsian one through a rational gauge transformation.
Any strictly fuchsian system reduces similarly
to a \emph{non resonant} one, that is, such that no two distinct 
\emph{exponents} (eigenvalues of $A(0)$) are congruent modulo 
$q^{\mathbf{Z}}$. This process involves some non canonical choices. \\

Second, any non resonant system $A$ is equivalent to the constant 
coefficients system $A(0)$. This is obtained by solving the 
functional equation with initial condition:
$$
\begin{cases}
F(0) = I_{n} \\
\left(\sigma_{q} F\right) A(0) = A F
\end{cases}
$$
with $F$ a convergent power series. The equivalent
equation $\sigma_{q} F = A F \left(A(0)\right)^{-1}$
then entails that $F$ is actually meromorphic over
$\mathbf{C}$ (the polarity will be precised below).
The transformation matrix $F$ is unique. \\

It follows that any system $A$ that is fuchsian 
at $0$ has a fundamental solution $X^{(0)} = M^{(0)} e_{q,A^{(0)}}$, 
where $M^{(0)} \in GL_{n}(\mathcal{M}(\mathbf{C}))$ and
$A^{(0)} \in GL_{n}(\mathbf{C})$ ; the latter can
be assumed to be non resonant and have all its
eigenvalues in the fundamental annulus. One can also
add more rigid normalising conditions, like sorted
eigenvalues (with respect to an arbitrary order on the fundamental
annulus) and sorted sizes of the Jordan blocks (\emph{see}
\cite{JSAIF}). Defining the \emph{singular locus of a matrix $M$}
to be:
$$
\mathcal{S}(M) = 
\{ \text{poles of } M \} \cup \{ \text{poles of } M^{-1} \} =
\{ \text{poles of } M \} \cup \{ \text{zeroes of } \det M \},
$$
one finds the singular locus of the meromorphic part:
$$
\mathcal{S}(M^{(0)}) = q^{- \mathbf{N}^{*}} \mathcal{S}(A).
$$
The singular locus of the ``log-car'' part $e_{q,A^{(0)}}$
is easily deduced from its definition.


\subsubsection*{1.2.3~~~Connection matrix and global classification}

We now assume $A$ to be fuchsian over $\mathbf{S}$
and attach to it local solutions $X^{(0)}$ and $X^{(\infty)}$ 
as defined in 1.2.2. We then define \emph{Birkhoff's connection matrix}
to be:
$$
P = \left(X^{(\infty)}\right)^{-1} X^{(0)}.
$$
It is clearly an elliptic invertible matrix. 
We attach to the system $A$ the (non uniquely defined)
triple $(A^{(0)},P,A^{(\infty)})$ (using the notations 
of 1.2.2). Changing the non unique choices, or changing
$A$ to a rationally equivalent system $B$ may produce
another triple $(B^{(0)},Q,B^{(\infty)})$.
If, in addition, one assumes normalised log-car parts 
(as described above), one must have 
$(A^{(0)},A^{(\infty)}) = (B^{(0)},B^{(\infty)})$
and there exists constant invertible matrices 
$R,S$ such that $R Q = P S$ plus another commutation condition 
more precisely explained in \cite{JSAIF}. 
There is then a bijective correspondence between classes
of fuchsian systems modulo rational equivalence
and classes of triples. The surjectivity follows from 
\emph{Birkhoff's lemma}, originating in \cite{Birkhoff1},
a modern version of which can be found in \cite{Arnold}.
The way we use it is detailed in \cite{JSAIF}, 2.2 and 2.3.
The tannakian counterpart of this bijection will be the content 
of 3.1 and 3.2.


\subsubsection*{1.2.4~~~Functors of solutions}

Specializing the constructions of 1.1.3 to the extension
$K' = \mathcal{M}(\mathbf{C}^{*})$ of $K = \mathbf{C}(z)$ provides us with a functor 
of solutions $\omega^{*} = S_{K'}$ from $\mathcal{E}$ to 
$Vect_{\mathcal{M}(\mathbf{E}_{q})}^{f}$.
In order to make $\omega^{*}$ a fibre functor, it would be enough
to guarantee, for each equation $A$, the existence of 
a fundamental matrix solution $X \in Gl_{n}(\mathcal{M}(\mathbf{C}^{*}))$. 
It follows from 1.2.3 that, for $A$ in $\mathcal{E}_{f}$, 
there is indeed such a fundamental matrix solution
and the rank of $\omega^{*}(A)$ is equal to the order of $A$.
Thus, the restriction of $\omega^{*}$ to $\mathcal{E}_{f}$ is a
fibre functor on $\mathcal{E}_{f}$ over the field 
$\mathcal{M}(\mathbf{E}_{q})$. However, we are rather looking 
for a fibre functor on $\mathcal{E}_{f}$ over 
the field $\mathbf{C}$. \\

We shall \emph{nearly} build such functors from the local solutions 
at $0$ and $\infty$:
$X^{(0)} = M^{(0)} e_{q,A(0)}$ and 
$X^{(\infty)} = M^{(\infty)} e_{q,A(\infty)}$. 
These solutions are not canonical,
but, by choosing such a pair for every equation, one is led
to the definition of two vector spaces of finite rank
\emph{over the complex numbers}. We thus define
the \emph{functors of solutions at $0$ and $\infty$} 
on $\mathcal{E}_{f}$ with values in $Vect_{\mathbf{C}}^{f}$ as:
$$
\omega^{(0)}:
\begin{cases}
A \leadsto X^{(0)} \mathbf{C}^{n} \\
(F: A \rightarrow B) \leadsto (U \mapsto F U)
\end{cases}
\qquad
\omega^{(\infty)}:
\begin{cases}
A \leadsto X^{(\infty)} \mathbf{C}^{n} \\
(F: A \rightarrow B) \leadsto (U \mapsto F U)
\end{cases}
$$
Let $F: A_{1} \rightarrow A_{2}$ be a morphism in $\mathcal{E}_{f}$. 
Then, with obvious notations, $F X_{1}^{(0)}$ is a solution of $A_{2}$,
hence of the form $X_{2}^{(0)} S^{(0)}$ where $S^{(0)}$ is elliptic ;
similarly at $\infty$. That the above definition makes sense comes from
the fact that $S^{(0)}$ and $S^{(\infty)}$ actually have coefficients
in $\mathbf{C}$. This follows from the following lemma, 
which we formulate in a slightly more general guise for further use. \\

\textbf{1.2.4.1 Lemma. -} 
\emph{Let $A,A'$ have constant coefficients and $F,M,M'$ have 
coefficients in $\mathcal{M}(\mathbf{C})$ and suppose that:
$$
F (M e_{q,A}) = (M' e_{q,A'}) S,
$$
where $S$ is elliptic. Then $S$ has constant coefficients.\\}

>From the conjugacy properties of matrices $e_{q,X}$
(viz, $e_{q,QXQ^{-1}} = Qe_{q,X}Q^{-1}$), one sees that 
$A$ and $B$ can be assumed to be upper triangular.
We write $A = A_{s} A_{u}$, $A' = A_{s}' A_{u}'$ their
Dunford decompositions, so that $e_{q,A} = e_{q,A_{s}} e_{q,A_{u}}$
and $e_{q,A'} = e_{q,A_{s}'} e_{q,A_{u}'}$. From the formulae
$e_{q,c} = z^{\epsilon(c)} e_{q,\overline{c}}$, we see that
we can assume all exponents (eigenvalues) of $A$, $A'$
to lie in the fundamental annulus $1 \leq |c| \lt |q|$.
We then write:
$$
(M'e_{q,A_{u}'})^{-1} F (M e_{q,A_{u}}) = e_{q,A_{s}'} S (e_{q,A_{s}})^{-1}.
$$
The $(i,j)$ coefficient of the right hand side is 
$s_{i,j} \frac{e_{q,c'_{i}}}{e_{q,c_{j}}}$, a chararacter
(since $s_{i,j}$ is elliptic) with coefficients in 
$\mathcal{M}(\mathbf{C})[l_{q}]$, because of the left hand
side. It follows from the independance properties of $q$-characters 
(\cite{JSAIF}, appendix) and from the location of $c'_{j},c_{i}$
in the fundamental annulus, that it must be trivial: $s_{i,j}$
must be a constant.
\hfill $\Box$ \\ 

Note that $s_{i,j}$ must be null if $c'_{j} \not= c_{i}$, which
amounts to say that $A_{s}' S = S A_{s}$. This being true after 
reduction to the fundamental annulus, the more general equality is:
$$
\overline{A_{s}'} S = S \overline{A_{s}}.
$$
We then have
$$
({M'}^{-1} F M) e_{q,A_{u}} = e_{q,A_{u}'} S 
$$
expanding over $\mathcal{M}(\mathbf{C})[l_{q}]$ and identifying
corresponding powers of $l_{q}$ (which is transcendental over
$\mathcal{M}(\mathbf{C})$, \emph{see} \emph{loc. cit.}) entails 
$$
({M'}^{-1} F M) A_{u} = A_{u}' S.
$$
These formulae are not affected by the conjugacies we used at the
beginning, but do depend on the reduction to the fundamental
annulus. \\

On the other hand, the functors $\omega^{(0)},\omega^{(\infty)}$
are not $\otimes$-compatible. In chapter 2, we shall eventually build
more easily fibre functors for fuchsian equations while avoiding 
the choice of particular solutions.


\subsubsection*{1.2.5~~~Singularities and exponents}

In order to compute the ``connection component'' 
of the Galois groupoid in chapter 3, we shall need 
to evaluate the connection matrix $P$ or its meromorphic
component $M$ at various points of $\mathbf{C}^{*}$, 
avoiding their singularities. Since $P$ is elliptic
and $M$ satisfies a simple $q$-difference equation,
these singularities are defined
modulo $q^{\mathbf{Z}}$ and may actually be seen 
as points in $\mathbf{E}_{q}$. They either come from 
the equation $A$ at stake, or from its local linear data
at $0$ and $\infty$. 
Let $(A^{(0)},P,A^{(\infty)})$ be a triple
coming from an object $A$ of $\mathcal{E}_{f}$ and write, as usual,
$X^{(0)} = M^{(0)} e_{q,A^{(0)}}$ and
$X^{(\infty)} = M^{(\infty)} e_{q,A^{(\infty)}}$
the local solutions used to compute $P$. Then, the singularities
of $P$ respectively come from:

\begin{enumerate}

\item{$\mathcal{S}(A)$ for the meromorphic part 
$M = (M^{(\infty)})^{-1} M^{(0)}$. Precisely, 
$\mathcal{S}(M) \subset q^{\mathbf{Z}} \mathcal{S}(A)$.}

\item{$Sp(A(0)),Sp(A(\infty))$ for the semi-simple components 
of the $e_{q,A(-)}$ parts. Since reduction to a constant matrix 
is not unique, these exponents of $A$ are actually defined
up to a factor in $q^{\mathbf{Z}}$.}

\item{$l_{q}$, the $q$-logarithm, in case $A(0)$ or $A(\infty)$
is not semi-simple.}

\end{enumerate}

On the other hand, we want to define $\otimes$-stable
categories only, so that we have to allow for multiplication of 
the exponents (and the inverse, to have stability under
dualisation). We thus define two categories 
$\mathcal{E}_{f,\Sigma}$ and $\mathcal{E}_{f,\Sigma,C}$.
Let $\Sigma$ be a finite subset of $\mathbf{C}^{*}$.
Then $\mathcal{E}_{f,\Sigma}$ is the full subcategory 
of $\mathcal{E}_{f}$ whose objects have all singularities 
in $q^{\mathbf{Z}} \Sigma$.
Let moreover $C$ be a finitely generated subgroup 
of the group $\mathbf{C}^{*}$.
Then $\mathcal{E}_{f,\Sigma,C}$ is the full subcategory 
of $\mathcal{E}_{f,\Sigma}$ whose objects have all
exponents in $q^{\mathbf{Z}} C$. 
They are obviously (strictly full) tannakian subcategories 
of $\mathcal{E}_{f}$ which is their inductive limit. 
More precisely, the tannakian subcategory $\left<A\right>$ 
generated by an arbitrary object $A$ (\emph{see} \cite{DM})
is contained in such a $\mathcal{E}_{f,\Sigma,C}$. 
We observe that, for $A$ in $\mathcal{E}_{f,\Sigma}$,
$\mathcal{S}(M) \subset q^{\mathbf{Z}} \Sigma$;
and for $A$ in $\mathcal{E}_{f,\Sigma,C}$,
$\mathcal{S}(P) \subset q^{\mathbf{Z}} \left(\Sigma \cup C \right)$.
Note that we impose (yet) no control on the $q$-logarithm
and unipotent parts. This will be done further below. 
Also note that we do not control separately the exponents at $0$ and $\infty$,
since we mainly seek to avoid the singularities of $P$. \\

\textbf{1.2.5.1 Proposition. -}
\emph{$A$ is in $\mathcal{E}_{f,\Sigma,C}$ if and only if its exponents
are in $C$ and $\mathcal{S}(M) \subset q^{\mathbf{Z}} \Sigma$.} \\

This follows from a precise use of Birkhoff's lemma
(\emph{see} \cite{JSAIF}, 2.3.1).
Here, of course, $M = e_{q,A^{(\infty)}} P (e_{q,A^{(0)}})^{-1}$. 
\hfill $\Box$ \\

We now consider regular equations at $0$, that is, such that 
$A(0) = I_{n}$.
For such an equation, the product formula:
$$
\underset{i \geq 1}{\prod} A(q^{-i}z) \underset{def}{=}
A(q^{-1}z) A(q^{-2}z) \cdots
$$ 
defines a solution holomorphic at $0$ and meromorphic on $\mathbf{C}$. 
Thus, if $A$ is regular at $0$ 
and at $\infty$, we obtain an explicit formula for 
the connection matrix $P$ (rather similar to Jacobi's triple product
formula for the function $\theta$), showing it to have as a singular
locus exactly $q^{\mathbf{Z}} \mathcal{S}(A)$. In this case, $M = P$.
Now, call more generally regular (at $0$) an equation that is equivalent to 
an $A$ such that $A(0) = I_{n}$. We define the category 
$\mathcal{E}_{f,reg}$ as the strictly full subcategory of 
$\mathcal{E}_{f}$ whose objects are equations regular at $0$ and 
$\infty$. It is clearly a tannakian subcategory of $\mathcal{E}_{f}$,
the one studied by Etingof in \cite{Etingof}. 



\section{Local theory}

The usual method in differential Galois theory for building
fibre functors is to use spaces of solutions. 
But our solutions have bad multiplicative properties: 
\emph{any} choice of solutions uniform over $\mathbf{C}^{*}$ 
will lead to some complicated twisting of the tensor product and 
of the connection matrix. We thus avoid to single out special 
fundamental solutions; in essence, this amounts to use sections 
of vector bundles without expressing them in frames.


\subsection{Localisation at $0$ and at $\infty$}


\subsubsection*{2.1.1~~~The categories $\mathcal{E}_{f}^{(0)}$
and $\mathcal{E}_{f}^{(\infty)}$}

The category $\mathcal{E}_{f}^{(0)}$ has the same objects 
as $\mathcal{E}_{f}$, 
but morphisms from $A$ (of rank $n$) to $B$ (of rank $p$)
are all $F \in M_{p,n}(\mathbf{C}(\{z\}))$ such that
$(\sigma_{q} F) A = B F$. The functional equation
$\sigma_{q} F = B F A^{-1}$ entails that the coefficients
of $F$ actually belong to $\mathcal{M}(\mathbf{C})$. \\

We define similarly $\mathcal{E}_{f}^{(\infty)}$
(morphisms with coefficients in $\mathcal{M}(\mathbf{C}_{\infty})$)
and $\mathcal{E}_{f}^{(*)}$ (morphisms with coefficients in 
$\mathcal{M}(\mathbf{C}^{*})$). These are plainly abelian rigid tensor 
categories, but $\mathcal{E}_{f}^{(*)}$ is 
$\mathcal{M}(\mathbf{E}_{q})$-linear while
$\mathcal{E}_{f}^{(0)}$ and $\mathcal{E}_{f}^{(\infty)}$
are only $\mathbf{C}$-linear (since a solution of $\sigma_{q} f = f$
meromorphic at $0$ or at $\infty$ has to be a constant).
It is clear that the natural embeddings 
$\mathcal{E}_{f} \rightarrow \mathcal{E}_{f}^{(0)}$,
$\mathcal{E}_{f} \rightarrow \mathcal{E}_{f}^{(\infty)}$,
$\mathcal{E}_{f}^{(0)} \rightarrow \mathcal{E}_{f}^{(*)}$ 
and $\mathcal{E}_{f}^{(\infty)} \rightarrow \mathcal{E}_{f}^{(*)}$
are $\mathbf{C}$-linear exact faithful $\otimes$-functors.
We intend to use this fact in the following way: any fibre
functor $\omega$ on $\mathcal{E}_{f}^{(0)}$ will restrict to
a fibre functor $\omega'$ on $\mathcal{E}_{f}$. This realizes
$Gal(\mathcal{E}_{f}^{(0)}) = Aut^{\otimes}(\omega')$ as this
subgroup of $Gal(\mathcal{E}_{f}) = Aut^{\otimes}(\omega)$
made up of elements natural with respect to \emph{all} morphisms 
in $\mathcal{E}_{f}^{(0)}$ (and not only in $\mathcal{E}_{f}$).


\subsubsection*{2.1.2~~~The categories $\mathcal{P}^{(0)}$ and
$\mathcal{P}^{(\infty)}$}

We call \emph{flat} an object of $\mathcal{E}_{f}^{(0)}$
(resp. $\mathcal{E}_{f}^{(\infty)}$, $\mathcal{E}_{f}^{(*)}$) 
if it is a matrix with coefficients in $\mathbf{C}$. These will be 
identified with flat vector bundles on $\mathbf{E}_{q}$ in 2.3.
The category $\mathcal{P}^{(0)}$ (resp. $\mathcal{P}^{(\infty)}$,
$\mathcal{P}^{(*)}$) is the full subcategory of $\mathcal{E}_{f}^{(0)}$
(resp. $\mathcal{E}_{f}^{(\infty)}$, $\mathcal{E}_{f}^{(*)}$)
whose objects are flat objects. These subcategories are obviously
$\mathbf{C}$-linear and stable under tensor operations. \\

\textbf{2.1.2.1 Lemma. -}
\emph{They are essential subcategories, that is, the inclusion
functors $\mathcal{P}^{(0)} \leadsto \mathcal{E}_{f}^{(0)}$,
$\mathcal{P}^{(\infty)} \leadsto \mathcal{E}_{f}^{(\infty)}$ and
$\mathcal{P}^{(*)} \leadsto \mathcal{E}_{f}^{(*)}$
are essentially surjective.} \\

Indeed, this is the content of the reduction to a constant coefficient
system by a meromorphic gauge transformation (see 1.2.2). 
\hfill $\Box$ \\

The following is an immediate consequence: \\

\textbf{2.1.2.2 Proposition. -}
\emph{$\mathcal{P}^{(0)}$ is a neutral tannakian category 
over $\mathbf{C}$, equivalent
to $\mathcal{E}_{f}^{(0)}$. Similar statements hold for 
$\mathcal{P}^{(\infty)}$, $\mathcal{P}^{(*)}$.
As a corollary, $\mathcal{P}^{(0)}$ and $\mathcal{P}^{(\infty)}$ 
have the same Galois group as $\mathcal{E}_{f}^{(0)}$ 
and $\mathcal{E}_{f}^{(\infty)}$ respectively.} 
\hfill $\Box$ 


\subsubsection*{2.1.3~~~Flat objects and equivariant morphisms of 
representations of $\mathbf{Z}$}

We do everything at $0$, 
the case of $\infty$ being similar.
We call $\mathcal{R}$ the category $Rep_{\mathbf{C}}(\mathbf{Z})$
of finite dimensional complex representations of of $\mathbf{Z}$.
These can be considered as pairs $(\mathbf{C}^{n},A)$ ($A$ a regular 
complex matrix of rank $n$), and the morphisms
$F: (\mathbf{C}^{n},A) \rightarrow (\mathbf{C}^{p},B)$ as
matrices $F \in M_{p,n}(\mathbf{C})$ such that $F A = B F$.
The category $\mathcal{R}$ is a $\mathbf{C}$-linear neutral tannakian 
category with the forgetful functor $\omega$ as fibre functor
and its Galois group: 
$$
\mathbf{Z}^{alg} = Aut^{\otimes}(\omega)
$$
is the \emph{proalgebraic hull} of $\mathbf{Z}$. Its structure and
the way it operates are recalled in 2.2.1. \\

\textbf{2.1.3.1 Proposition. -}
\emph{The natural embedding $\mathcal{R} \rightarrow \mathcal{P}^{(0)}$
is a faithful (but not full), essentially surjective exact
$\mathbf{C}$-linear $\otimes$-functor.} \\

The categories at stake have the same objects but 
$\mathcal{R}$ has only \emph{constant} morphisms,
while $\mathcal{P}^{(0)}$ has all \emph{equivariant}
morphisms (the word will be justified in 2.3).
The proof is therefore immediate.
\hfill $\Box$ \\

\textbf{2.1.3.2 Lemma. -}
\emph{Any morphism $F: A \rightarrow B$ in $\mathcal{P}^{(0)}$
is a Laurent polynomial:
$$
F = \sum F_{k} z^{k},
$$
a finite sum where each $F_{k}$ has coefficients in $\mathbf{C}$.} \\

Indeed, one has $\sigma_{q} F = B F A^{-1}$, so that $F$ has only
a pole at $0$: otherwise, these would be propagated along 
a $q$-spiral and would accumulate at $0$, contradicting the
meromorphy of $F$. Now, $F$ has a Laurent \emph{series} $F = \sum F_{k} z^{k}$.
Identifying the $z^{k}$ components of the two sides of the functional
equation, we get: $q^{k} F_{k} A = B F_{k}$. This is possible
with $F_{k} \not= 0$ only if $Sp(q^{k} A) \cap Sp(B) \not= \emptyset$
(see, for instance, the lemma in 1.1.3 of \cite{JSAIF}, also proven
in \cite{Wasow} or \cite{Weil}). Hence the finiteness of the sum.
\hfill $\Box$ \\

\textbf{2.1.3.3 Proposition. -}
\emph{$(\mathbf{C}^{n},A) \leadsto \mathbf{C}^{n}$, 
$F \leadsto F(z_{0})$ gives a fibre functor
$\omega^{(0)}_{z_{0}}$ on $\mathcal{P}^{(0)}$
for any $z_{0} \in \mathbf{C}^{*}$.} \\

The only non trivial point is the faithfulness. Suppose $F(z_{0}) = 0$.
Then the functional equation implies 
$\forall k \in \mathbf{Z} \;,\; F(q^{k} z_{0}) = 0$, thus $F = 0$
($F$ being a Laurent polynomial).
\hfill $\Box$ \\

\textbf{2.1.3.4 Corollary. -}
\emph{The local Galois group $G^{(0)}$ is a closed subgroup
of $\mathbf{Z}^{alg}$.} \\

Since the previous fibre functor restricts to the forgetful functor 
on $\mathcal{R}$, it follows from \cite{DM}, prop. 2.21 that we have 
realized $G^{(0)}$ as a subgroup of the Galois group $\mathbf{Z}^{alg}$ 
of $\mathcal{R}$. Precisely, 
restriction of the elements of $Aut^{\otimes}(\omega^{(0)}_{z_{0}})$
to $\mathcal{R}$ gives a morphism of groups
$G^{(0)} \rightarrow \mathbf{Z}^{alg}$ and,
according to \emph{loc. cit},
this is a
closed immersion of proalgebraic groups.
\hfill $\Box$


\subsection{The local Galois groupoid and the local monodromy}


\subsubsection*{2.2.1~~~The proalgebraic hull of $\mathbf{Z}$}

We summarize here the main results we shall need. More details 
and complete proofs can be found in \cite{JSthese} and in \cite{JSGAL}. \\

Our basic tool for building a fundamental group with
some topological meaning is the ``hole'':
$\mathbf{Z} = \pi_{1}(\mathbf{C}^{*})$.
Its proalgebraic hull $\mathbf{Z}^{alg}$ (see the beginning of 2.1.3)
is commutative. It is the product 
$\mathbf{Z}^{alg}_{s} \times \mathbf{Z}^{alg}_{u}$ 
of its semi-simple part:
$$
\mathbf{Z}^{alg}_{s} = Hom_{grp}(\mathbf{C}^{*},\mathbf{C}^{*}),
$$
the group of characters of the abstract group
$\mathbf{C}^{*}$, and of its unipotent part:
$$
\mathbf{Z}^{alg}_{u} = \mathbf{C}.
$$
The \emph{fundamental loop in $\pi_{1}(\mathbf{C}^{*})$} 
is represented by $1 \in \mathbf{Z}$; it is
a Zariski-generator of $\mathbf{Z}^{alg}$, where
it splits into a semi-simple component: 
$Id_{\mathbf{C}^{*}} \in \mathbf{Z}^{alg}_{s}$ and
a unipotent component $1 \in \mathbf{Z}^{alg}_{u}$.
We talk of \emph{fundamental semi-simple loop} and 
\emph{fundamental unipotent loop}. \\

A representation of $\mathbf{Z}$ is the same thing as
a pair $(V,f)$ of a finite dimensional complex vector
space $V$ and an automorphism $f \in GL(V)$. So let
$(\mathbf{C}^{n},A)$, $A \in GL_{n}(C)$, be a (standard
model of a) generic object of $\mathcal{R}$.
Let $A = A_{s} A_{u}$ be the Dunford
decomposition (see the end of the introduction). Let 
$(\gamma,\lambda) \in 
Hom_{grp}(\mathbf{C}^{*},\mathbf{C}^{*}) \times \mathbf{C}$
be a Galoisian automorphism. Then we put:
$$
A^{(\gamma,\lambda)} = \gamma(A_{s}) \; A_{u}^{\lambda} = 
A_{u}^{\lambda} \; \gamma(A_{s}).
$$
Here, $\gamma$ acts on $A_{s}$ through its eigenvalues:
if $A_{s} = Q \; \Diag(c_{1},\ldots,c_{n}) \; Q^{-1}$,
one has 
$\gamma(A_{s}) = Q \; \Diag(\gamma(c_{1}),\ldots,\gamma(c_{n})) \; Q^{-1}$.
Since $A^{(\gamma,\lambda)} \in GL_{n}(C)$, it defines
an automorphism of $\mathbf{C}^{n} = \omega(\mathbf{C}^{n},A)$.
This is precisely how $(\gamma,\lambda)$ is incarnated as
an element of $Aut^{\otimes}(\omega)$.


\subsubsection*{2.2.2~~~The local Galois groupoid}

Any fibre functor $\omega^{(0)}_{z_{0}}$ (as defined in 2.1.3.3) 
obviously restricts to the
forgetful functor $\omega$ on $\mathcal{R}$. Therefore,
$Aut^{\otimes}(\omega^{(0)}_{z_{0}})$ is a subgroup (and
$Iso^{\otimes}(\omega^{(0)}_{z_{0}},\omega^{(0)}_{z_{1}})$
a subset) of $Aut^{\otimes}(\omega)$. \\

\textbf{2.2.2.1 Theorem: the local Galois groupoid. -}
\emph{With the previous identification of $Aut^{\otimes}(\omega)$
with $\mathbf{Z}^{alg}$,
$$
Iso^{\otimes}(\omega^{(0)}_{z_{0}},\omega^{(0)}_{z_{1}}) = 
\{(\gamma,\lambda) \in \mathbf{Z}^{alg} \;/\; \gamma(q)z_{0} = z_{1} \}.
$$}

We thus obtain a local Galois groupoid at $0$, which
we call $G^{(0)}$. Its base set is $\mathbf{C}^{*}$, 
the arrow sets 
$G^{(0)}(a,b) = Iso^{\otimes}(\omega^{(0)}_{a},\omega^{(0)}_{b})$
being described in the theorem. It is clearly connected
on the base, so that we can take as local Galois group 
any group $G^{(0)}(a,a)$. The proalgebraic structure
on each of the $G^{(0)}(a,b)$ is then induced by that
of $\mathbf{Z}^{alg}$, of which they are Zariski
closed subsets. \\

\textbf{2.2.2.2 Corollary: the local Galois group. -}
\emph{With the same identifications,
$$
Aut^{\otimes}(\omega^{(0)}_{z_{0}}) = 
\{(\gamma,\lambda) \in \mathbf{Z}^{alg} \;/\; \gamma(q) = 1 \}.
$$}
\hfill $\Box$ \\

Let us now proceed to prove the theorem.
Let $(\gamma,\lambda)$ be an element of 
$Iso^{\otimes}(\omega^{(0)}_{z_{0}},\omega^{(0)}_{z_{1}})$. 
Applied to the morphism
$z: (\mathbf{C},1) \rightarrow (\mathbf{C},q)$
in $\mathcal{P}^{(0)}$,
it gives rise to the commutative diagram:

\begin{equation*}
\begin{CD}
\mathbf{C} @>z_{0}>> \mathbf{C} \\
@V{\gamma(1)}VV                          @VV{\gamma(q)}V      \\
\mathbf{C} @>z_{1}>> \mathbf{C}
\end{CD}
\end{equation*}

hence $\gamma(q) z_{0} = z_{1}$. \\

Conversely, suppose $\Phi = (\gamma,\lambda)$ with $\gamma(q) = 1$.
We must check its compatibility with all morphisms in 
$\mathcal{P}^{(0)}$ (and not only in $\mathcal{R}$).
That is, $F: (\mathbf{C}^{n},A) \rightarrow (\mathbf{C}^{p},B)$ 
being such a morphism, we want to show that the diagram:

\begin{equation*}
\begin{CD}
\mathbf{C}^{n} @>F(z_{0})>> \mathbf{C}^{p}  \\
@V{\Phi(A)}VV                          @VV{\Phi(B)}V      \\
\mathbf{C}^{n} @>F(z_{1})>> \mathbf{C}^{p}
\end{CD}
\end{equation*}

is commutative, where $z_{0},z_{1} \in \mathbf{C}^{*}$
and $a z_{0} = z_{1} \;,\; a =\gamma(q)$. But we have 
$\left(\sigma_{q} F \right) A = B F$ thus, as noticed before,
$\forall k \in \mathbf{Z} \;,\; F_{k} (q^{k} A) = B F_{k}$
so that $F_{k}: q^{k} A \rightarrow B$ is a morphism in
$\mathcal{R}$. For this morphism we check the naturality of $\Phi$.
Then:
\begin{eqnarray*}
F_{k} \Phi(q^{k}A)            & = & \Phi(B) F_{k} \\ 
\Rightarrow F_{k} a^{k} \Phi(A)           & = & \Phi(B) F_{k} \\ 
\Rightarrow F_{k} a^{k} z_{0}^{k} \Phi(A) & = & \Phi(B) F_{k} z_{0}^{k}  
\end{eqnarray*}
Then, summing on $k$ gives:
$$
F(a z_{0}) \Phi(A) = \Phi(B) F(z_{0}).
$$
\hfill $\Box$ \\

\textbf{2.2.2.3 Incarnation of the local Galois groupoid. -}
To incarnate $G^{(0)}$, one must show how its elements 
operate on an object $A$ of $\mathcal{E}_{f}^{(0)}$. 
So we take two base points $a,b \in \mathbf{C}^{*}$
and an element $(\gamma,\lambda) \in G^{(0)}(a,b)$.
We must obtain an isomorphism from $\omega^{(0)}_{a}(A)$
to $\omega^{(0)}_{b}(A)$. Both vector spaces are equal to 
the space $\mathbf{C}^{n}$ underlying an object $A^{(0)}$
of $\mathcal{P}^{(0)}$ equivalent to $A$ (2.1.2.1). From 2.2.1, 
we see that the matrix describing the isomorphism we are looking for 
is $\gamma(A^{(0)}_{s}) \left(A^{(0)}_{u}\right)^{\lambda}$.
Note that the semi-simple factor actually depends only on 
$\overline{A^{(0)}_{s}}$, not on $A^{(0)}_{s}$.


\subsubsection*{2.2.3~~~The local fundamental group}

We stick to our overloaded notation $G^{(0)}$ for
the local Galois groupoid at $0$ and at the same time
for any of the local Galois groups, realized as:
$$
G^{(0)}(a,a) \simeq
Hom_{grp}(\mathbf{C^{*}}/q^{\mathbf{Z}},\mathbf{C^{*}}) \times \mathbf{C}.
$$
We shall now exhibit explicit elements in the local Galois group, 
aiming at a Zariski-dense finitely presented discrete group 
with a topological flavour. In the same way as the fundamental loop 
$1 \in \mathbf{Z} \subset \mathbf{Z}^{alg} =
Hom_{grp}(\mathbf{C^{*}},\mathbf{C^{*}}) \times \mathbf{C}$
splits into semi-simple and a unipotent parts, we shall
consider separately the semi-simple and unipotent components
of these elements.
As for the latter, we shall take the obvious candidate:
$1 \in \mathbf{C}$, which corresponds to the automorphism
$\mathcal{X} \leadsto A_{u}^{(0)}$ of $\omega^{(0)}$.
It generates a Zariski-dense subgroup of the unipotent component 
$G_{u}^{(0)}$. We concentrate hereafter on the semi-simple component 
$G^{(0)}_{s}$. \\

\textbf{2.2.3.1 Explicit elements of $G_{s}^{(0)}$. -}
The topological group $\mathbf{C^{*}}/q^{\mathbf{Z}}$ being compact, 
all the elements of $G_{s}^{(0)}$ built from \emph{continuous}
morphisms have their image in the biggest compact subgroup of 
$\mathbf{C}^{*}$, the unit circle $\mathbf{U}$. 
We now proceed to describe them. \\

We write $q = e^{- 2 \imath \pi \tau}$, $Im(\tau) > 0$,
whence the splitting
$\mathbf{C}^{*} = \mathbf{U} \times q^{\mathbf{R}}$, where we put, 
for $y \in \mathbf{R}$, $q^{y} = e^{- 2 \imath \pi \tau y}$
and $q^{\mathbf{R}} = \{q^{y} \;/\; y \in \mathbf{R}\}$. 
This comes (through the lift $\mathbf{C} \rightarrow \mathbf{C}^{*}$,
$x \mapsto e^{2 \imath \pi x}$) from the splitting
$\mathbf{C} = \mathbf{R} \oplus \mathbf{R} \tau$. 
We then define: 
$$
\gamma_{1}:
\begin{cases}
\mathbf{C}^{*} \rightarrow \mathbf{C}^{*} \\
u q^{y} \mapsto u
\end{cases}
\qquad
\gamma_{2}:
\begin{cases}
\mathbf{C}^{*} \rightarrow \mathbf{C}^{*} \\
u q^{y} \mapsto e^{2 \imath \pi y}
\end{cases}
$$
More generally, one can define, for $b \in \mathbf{C}$, a group
morphism $\gamma_{2}^{b} : u q^{y} \mapsto e^{2 \imath \pi b y}$. 
This will be used in 4.1.1. \\

\textbf{2.2.3.2 Lemma. -}
\emph{The subgroup of continuous elements of $G^{(0)}_{s}(a,a)$
is generated by $\gamma_{1}$ and $\gamma_{2}$.} \\

All continuous group morphisms from $\mathbf{C}$ to $\mathbf{C}$
have the form $x + y \tau \mapsto \alpha x + \beta y$ with 
$\alpha,\beta \in \mathbf{C}$. Continuous group morphisms 
from $\mathbf{C}^{*}$ 
to $\mathbf{C}^{*}$ are obtained (through the same lifting as above)
from those that send $\mathbf{Z}$ to itself, that is, those such that 
$\alpha \in \mathbf{Z}$. Such a morphism factors through 
$\mathbf{C^{*}}/q^{\mathbf{Z}}$ (i.e., it maps $q$ to $1$)
if and only if $\beta \in \mathbf{Z}$. 
\hfill $\Box$ \\

\textbf{2.2.3.3 Lemma: a density criterion. -}
\emph{Let $\Gamma$ be a subset of 
$G_{s}^{(0)} = Hom_{grp}(\mathbf{C^{*}/q^{\mathbf{Z}}},\mathbf{C^{*}})$.
Then $\Gamma$ generates a Zariski-dense subgroup of
$G_{s}^{(0)}$ if and only if
$\underset{\gamma \in \Gamma}{\bigcap} \Ker \gamma = \{1\}$.} \\

By Chevalley's criterion , used in a similar way
in \cite{DM} (p. 129, proof of prop. 2.8; the precise formulation 
is given in \emph{loc. cit.}, first lecture, p. 40, prop. 3.1.b 
and p. 41, remark 3.2.a),
the density condition amounts to the following: 
for any object $A$ of $\mathcal{E}_{f}^{(0)}$, 
if a line of $\omega^{(0)}(A)$ 
is stabilised by all $\gamma(A)$, $\gamma \in \Gamma$, 
it is stabilised by all $\gamma(A)$, $\gamma \in G_{s}^{(0)}$.
This means that any common eigenvector of all the
$\gamma(A^{(0)}_{s})$, $\gamma \in \Gamma$ 
is actually a common eigenvector of all the
$\gamma(A^{(0)}_{s})$, $\gamma \in G_{s}^{(0)}$. \\

We now fix such an object $A$ and a non-zero
$x \in \omega^{(0)}(A) = \mathbf{C}^{n}$.
Up to conjugacy, one may assume that 
$A^{(0)}_{s} = \Diag(c_{1},\ldots,c_{n})$. 
Then, $x$ is an eigenvector of $\gamma(A^{(0)}_{s})$
if and only if, for any indices $i \not= j$
such that $x_{i} x_{j} \not= 0$ (let us say
that such indices are \emph{linked}), one has
$\gamma(\overline{c_{i}}) = \gamma(\overline{c_{j}})$
(recall from the introduction that $\overline{c}$ denotes
the image of $c$ in $\mathbf{C^{*}/q^{\mathbf{Z}}}$).
Thus, the elements of $G_{s}^{(0)}$ whose action 
on $\omega^{(0)}(A)$ fixes the line $\mathbf{C} x$
are precisely those such that $\overline{c_{i}/c_{j}} \in \Ker \gamma$ 
for any pair of linked indices $(i,j)$. 
The lemma follows immediately.
\hfill $\Box$ \\

\textbf{2.2.3.4 Two generating loops.} 
According to the previous lemma, we are to choose elements
of $Hom_{grp}(\mathbf{C^{*}},\mathbf{C^{*}})$
sending $q$ to $1$ and such that the intersection of their
kernels is $q^{\mathbf{Z}}$. This can't be done
with one element, since $\mathbf{C^{*}}/q^{\mathbf{Z}}$ does
not embed in $\mathbf{C}^{*}$ (for instance, because the former 
has $4$ elements of order $2$, while the latter has $2$).
To have anything explicit, we have to select among the
morphisms found in 2.3.1. Since it is clear that 
$\Ker\gamma_{1} \cap \Ker \gamma_{2} = q^{\mathbf{Z}}$,
these two elements fit our needs. Note that $\gamma_{1}$
and $\gamma_{2}$ generate a free abelian group. \\

\textbf{2.2.3.5 Theorem. -}
\emph{The subgroup of
$Hom_{grp}(\mathbf{C^{*}},\mathbf{C^{*}}) \times \mathbf{C}$
whose unipotent component is $\mathbf{Z} \subset \mathbf{C}$
and whose semi-simple component is generated by $\gamma_{1}$
and $\gamma_{2}$ is Zariski-dense in the local Galois group.}
\hfill $\Box$ \\

This subgroup can therefore be considered as the \emph{local 
fundamental group}. As a consequence, all the  semi-simple parts 
of our ``monodromy matrices'' will be unitary, in contrast to 
the classical (differential) case. \\

Since continuous elements of the Galois groupoid
form a torsor under the action of any one of the
groups just determined, it is enough, to get them
all, to find \emph{one} such element $g_{b,c}$
in each $G^{(0)}_{s}(b,c)$. To make it a groupoid, 
it is therefore enough to put $\Pi(a,b) = g_{a,b} \Pi(a,a)$,
where $\Pi(a,a)$ is the group we just defined. \\

\textbf{2.2.3.6 Corollary: the local fundamental groupoid. -}
\emph{Choose, for each $a \in \mathbf{C}^{*}$,
a logarithm $\alpha \in \frac{1}{2 \imath \pi} \log a$.
Put $g_{a} = \gamma_{2}^{\alpha}$. One can then take
$g_{b,c} = g_{\frac{c}{b}}$.}
\hfill $\Box$ \\

One can't expect $g_{a}$ to be a continuous function of $a$,
at the best, the choice of a logarithm involves a cut: 
the discontinuity at the cut will be interpreted in 2.2.4.
There is also another interpretation of $g_{a}$
as an automorphism of the ``field of solutions'', corresponding
to the translation $z \mapsto a z$ (in multiplicative notation)
of the elliptic curve $\mathbf{E}_{q}$, \emph{see} \cite{JSGAL}.
We postpone the geometrical interpretation of the local
fundamental group and groupoid to 2.4.


\subsection{Flat vector bundles over the elliptic curve $\mathbf{E}_{q}$}

We give here a geometrical interpretation of $\mathcal{P}^{(0)}$,
close to Weil's correspondence between flat vector bundles 
on a curve and  representations of its fundamental
group (\emph{see} \cite{Weil}). \\

Write $Fib(X)$ the category of holomorphic vector bundles over 
a compact Riemann surface $X$ and $Fib_{p}(X)$ the full subcategory
of those which are \emph{flat}, i.e. those whose transition matrices
can be taken to be constant (vs holomorphic) for some adequate
covering. Note that the \emph{morphisms} between such flat bundles
are \emph{not} required to be constant. 


\subsubsection*{2.3.1~~~Holomorphic sections of a flat bundle over $\mathbf{E}_{q}$}

Let $A \in Gl_{n}(\mathbf{C})$, thus an object of $\mathcal{P}^{(0)}$
as well as a linear $q$-difference system with constant coefficients.
One introduces the equivalence relation $\sim_{A}$ on 
$\mathbf{C}^{*} \times \mathbf{C}^{n}$ generated by the relations:
$$
\forall (z,X) \in \mathbf{C}^{*} \times \mathbf{C}^{n} \;,\;
(z,X) \sim_{A} (qz,AX).
$$
The first projection 
$\mathbf{C}^{*} \times \mathbf{C}^{n} \rightarrow \mathbf{C}^{*}$
is compatible with the action of $q^{\mathbf{Z}}$ on $\mathbf{C}^{*}$
and, factoring out, we define a holomorphic vector bundle of rank $n$:
$$
F_{A} \underset{def}{=} \frac{\mathbf{C}^{*} \times \mathbf{C}^{n}}{\sim_{A}}
\rightarrow \mathbf{E}_{q} = \frac{\mathbf{C}^{*}}{q^{\mathbf{Z}}}.
$$
This is a particular case of Weil's correspondence alluded to above:
the fundamental group of $\mathbf{E}_{q}$ is the lattice:
$$
\pi_{1}(\mathbf{E}_{q}) = \Lambda_{\tau} = \mathbf{Z} + \mathbf{Z} \tau,
$$
where $q = e^{- 2 \imath \pi \tau}$ and one takes the representation
$\Lambda_{\tau} \rightarrow Gl_{n}(\mathbf{C})$ sending $1$ to $I_{n}$
and $\tau$ to $A$. \\

\textbf{2.3.1.1 Lemma. -}
\emph{There is a natural bijection between the solutions of the system
with matrix $A$ meromorphic on $\mathbf{C}$ and the holomorphic 
sections of $F_{A}$.} \\

There is an obvious bijection between the holomorphic sections
of $F_{A}$ and the solutions of the system with matrix $A$ 
holomorphic on $\mathbf{C}^{*}$.
We have to show that the latter uniquely extend to meromorphic 
solutions on $\mathbf{C}$ (this, without having to prescribe growth
conditions). We have already seen (2.1.3.2)
that a solution meromorphic 
on $\mathbf{C}$ has to be holomorphic on $\mathbf{C}^{*}$ to prevent 
the accumulation of poles at $0$. So, let $F$ be holomorphic 
on $\mathbf{C}^{*}$ and such that $F(qz) = A F(z)$. 
Call $M$ the maximum of $||F||$ on the compact
fundamental annulus $1 \leq |z| \leq |q|$ 
(this, for an arbitrary norm $||-||$)
and, for $z \in \mathbf{C}^{*}$, put 
$k = \left\lfloor\frac{\ln |z|}{\ln |q|}\right\rfloor$
(where $\lfloor x \rfloor)$ denotes the biggest integer 
less than or equal to $x$).
One has:
\begin{eqnarray*}
||F(q^{-k}z)|| \leq M & \Rightarrow & ||F(z)|| \leq M |||A|||^{k} \\
& \Rightarrow & ||F(z)|| \leq |z|^{\frac{\ln |||A|||}{\ln |q|}}
\end{eqnarray*}
($|||-|||$ the subordinate norm). This entails polynomial growth
at $0$, hence a pole.
\hfill $\Box$ \\

\textbf{2.3.1.2 Remark. -}
For a \emph{unitary} bundle, one can prove that all sections 
are actually constant (\emph{see} \cite{Weil}, \cite{Seshadri}).


\subsubsection*{2.3.2~~~Comparison of the categories 
$\mathcal{P}^{(0)}$ and $Fib_{p}(\mathbf{E}_{q})$}

Since every morphism in $\mathcal{P}^{(0)}$ is holomorphic 
on $\mathbf{C}^{*}$, it defines a holomorphic map between the 
corresponding vector bundles and we clearly have a $\otimes$-functor 
from $\mathcal{P}^{(0)}$ to $Fib_{p}(\mathbf{E}_{q})$. \\

\textbf{2.3.2.1 Theorem. -}
\emph{This is a $\otimes$-equivalence.} \\

The full faithfulness comes from the existence of internal
$Hom$ in both categories, implying that the morphisms
$A \rightarrow B$ (resp. $F_{A} \rightarrow F_{B}$) are in
natural bijection with solutions of $\underline{Hom}(A,B)$
(resp. $\underline{Hom}(F_{A},F_{B})$) and from lemma 2.3.1.1.
Now, as regards essential surjectivity, let a vector bundle 
over $\mathbf{E}_{q}$ correspond to a representation of 
$\Lambda_{\tau}$ that maps $1$ to $A$
and $\tau$ to $B$, these being commuting elements of 
$Gl_{n}(\mathbf{C})$. One writes $A = e^{2 \imath \pi U}$
and uses the gauge transformation $e^{2 \imath \pi x U}$
to reduce this representation to one that sends $1$ to $I_{n}$
(and $\tau$ to $A^{-\tau}B$). From the construction in 2.3.1,
this comes from an object of $\mathcal{P}^{(0)}$.
\hfill $\Box$ 


\subsubsection*{2.3.3~~~Fibre functors}

The category $Fib_{p}(\mathbf{E}_{q})$ is a thickening of 
$Rep_{\mathbf{C}}(\pi_{1}(\mathbf{E}_{q}))$ since it has
the same objects but more morphisms: if 
$\phi: \pi_{1}(\mathbf{E}_{q}) \rightarrow Gl(V)$ and
$\psi: \pi_{1}(\mathbf{E}_{q}) \rightarrow Gl(W)$ are
two such representations, a morphism between the associated
bundles gives rise to an \emph{equivariant} morphism
$\phi \rightarrow \psi$, that is, a holomorphic map
$F: \mathbf{C} \rightarrow \mathcal{L}(V,W)$ such that:
$$
\forall \gamma \in \pi_{1}(\mathbf{E}_{q}) \;,\; 
\forall x \in \mathbf{C} \;,\;
\phi(\gamma) \circ F(x) = F(\gamma x) \circ \psi(\gamma).
$$
This is a morphism in $Rep_{\mathbf{C}}(\pi_{1}(\mathbf{E}_{q}))$
if and only if $F$ is a constant map. Here, $\gamma$ operates
on $x$ via the action of $\pi_{1}(\mathbf{E}_{q})$ on the
universal covering $\mathbf{C}$ of $\mathbf{E}_{q}$ (that is, 
the translation action of $\Lambda_{\tau}$ on $\mathbf{E}_{q}$). \\

\textbf{2.3.3.1 Punctual fibre functors. -}
Therefore, any fibre functor on $\mathcal{P}^{(0)}$ and 
$Fib_{p}(\mathbf{E}_{q})$ naturally induces the usual 
fibre functor on $Rep_{\mathbf{C}}(\pi_{1}(\mathbf{E}_{q}))$ 
and, by the very same trick we already used in 2.1.3.4 and 2.2.2, 
we obtain
$Gal(Fib_{p}(\mathbf{E}_{q}))$ as a proalgebraic subgroup of
$\pi_{1}(\mathbf{E}_{q})^{alg}$. The latter is isomorphic to
$\mathbf{Z}^{alg} \oplus \mathbf{Z}^{alg}$. Actually, only
one factor is really involved here since every bundle has been shown
in 2.3.2.1 to be isomorphic to one on which $1$ acts trivially (that
is, one which already trivializes on the quotient covering
$\mathbf{C}^{*}$ of the universal covering $\mathbf{C}$). \\

\textbf{2.3.3.2 The global fibre functor. -}
On the other hand, the equivalence of $\mathcal{P}^{(0)}$ 
with $Fib_{p}(\mathbf{E}_{q})$ may itself be viewed as
a fibre functor on $\mathcal{E}_{f}^{(0)}$ in the following way. 
Call $S$ the unique curve (scheme) over $\mathbf{C}$ with
underlying analytic space $S^{an} = \mathbf{E}_{q}$. 
Then, holomorphic vector bundles over $\mathbf{E}_{q}$ ``are'' 
locally free sheaves and we get a fibre functor over $S$
in the sense of \cite{DF} . 
The theorem 1.12 of \emph{loc. cit.} then provides us with a groupoid 
structure over $S$: \\ 

\textbf{2.3.3.3 Theorem. -}
\emph{The tensor category $\mathcal{E}_{f}^{(0)}$
is equivalent to the category of representations
of a groupoid $G^{(0)}$ that is faithfully flat 
over $S \times S$.}
\hfill $\Box$ \\

In essence, this says that $G^{(0)}$ acts transitively
on the base $S$ and that composition of paths is a morphism with
respect to the proalgebraic structure on each
$G^{(0)}(a,b)$, but also with respect to the 
algebraic structure on the base $S$.


\subsubsection*{2.3.4~~~The classification theorem of Baranovsky and Ginzburg}

Let $G$ denote a complex algebraic group and write $G((z))$ for 
the group of $\mathbf{C}((z))$ rational points of $G$, where 
$\mathbf{C}((z))$ is the field of formal Laurent series, a kind
of ``loop group''. Write $G[[z]]$ for its subgroup of
$\mathbf{C}[[z]]$ rational points. Then define a ``twisted'' 
conjugation action of $G((z))$ on itself by putting, for 
$a(z),g(z) \in G((z))$:
$$
{}^{g}a(z) = g(q.z).a(z).g(z)^{-1}.
$$
This twisted conjugation action can actually be seen as a plain
conjugation action in some larger group: putting
$a(z) \mapsto a(t.z)$ defines a $\mathbf{C}^{*}$-action
on the loop group $G((z))$ (the ``rotation of the loop''),
hence a semi-direct product, and twisted conjugacy classes are
actually ordinary conjugacy classes in a Kac-Moody group. \\

The following is stated and proved in \cite{BG}: \\

\textbf{Theorem. -}
\emph{If $G$ is connected and semi-simple, there is a natural 
bijection between the set of integral twisted conjugacy classes
in $G((z))$ and the set of isomorphism classes of semi-stable
holomorphic principal $G$-bundles on the elliptic curve
$E = \mathbf{C}^{*}/q^{\mathbf{Z}}$.} \\

Here, integral twisted conjugacy classes are those which contain
an element of $G[[z]]$: they are analogous to our fuchsian
equations. One subtlety of this theorem (and the main difficulty
in its proof) is that it provides an \emph{analytic}
classification of \emph{formal} objects. The authors give a 
tannakian extension of this theorem. They define a tensor
category $\mathcal{M}_{q}$ of formal integral $q$-difference
modules and prove: \\

\textbf{Theorem. -}
\emph{It is equivalent to the tensor category of degree zero
semi-stable holomorphic vector bundles on $E$.} \\

In the electronic (IMRN) version of their paper, they moreover
quote a computation by Kontsevich of the corresponding Galois
group: the result is the same as our corollary 2.2.2.2.


\subsection{Heuristic topological interpretation and perspectives}

\subsubsection*{2.4.1~~~A ``local'' elliptic curve}

Recall that our constructions aim at a \emph{geometric} understanding 
of $q$-difference equations. Extending the class of 
morphisms of the category $\mathcal{E}_{f}$ as we did in 2.1 amounts to 
a localisation for the transcendental topology. Accordingly,
our vector bundles on $\mathbf{E}_{q}$ can be considered as
induced by equivariant vector bundles on the germ of complex space
$(\mathbf{C}^{*},0)$ and the curve $\mathbf{E}_{q}$ itself as
the quotient of the germ $(\mathbf{C}^{*},0)$ by the action
of the dilatation $\sigma_{q}$. We visualize this curve as
``local at $0 \in \mathbf{S}$'', since it carries the 
\emph{local geometry of fuchsian $q$-difference equations}.
The loops found in 2.2.3 can been interpreted as the two fundamental 
loops of $\mathbf{E}_{q}$. They are algebraically free
and generate a free abelian group of rank $2$. However, they
satisfy a ``transcendental relation'', as predicted in
\cite{RamisJPRTraum}:
$$
\gamma_{1} \gamma_{2}^{-\tau} = Id_{\mathbf{C}^{*}},
$$
the semi-simple fundamental loop of $\mathbf{Z}^{alg}$.
Suppose the logarithms involved in the definition
of the $g_{a}$ have been chosen continuously, but 
for some cut: for instance, the main determination
(alternatively, see 3.2.2).
Let $a$ turn counterclockwise once around $0$. Then 
$g_{a}(u q^{y}) = e^{- 2 \imath \pi \tau \alpha y}$
is multiplied by $e^{-2 \imath \pi y}$, that is,
by $\gamma_{2}(u q^{y})$. This means that $\gamma_{2}$
represents the plain classical loop around $0$ 
in $\mathbf{C}^{*}$. This can be seen yet another
way: suppose we change our choice of
a logarithm of $q$, writing $q = e^{- 2 \imath \pi \tau'}$,
where $\tau' = \tau + m$, $m \in \mathbf{Z}$.
Then, the formulas in 2.2.3.1 produce modified loops
$\gamma'_{1}$ and $\gamma'_{2}$ and one checks that
$\gamma'_{1} = \gamma_{1} \gamma_{2}^{m}$ 
and $\gamma'_{2} = \gamma_{2}$. 
To interpret $\gamma_{1}$ as the second generating loop
of $\mathbf{E}_{q}$ is not so easy while
staying in a strict local context. It will be seen in 3.2.2.2 
to be the defect of ellipticity of the twisted connection matrix, 
which is the generic Galois isomorphism linking $0$ to $\infty$.
This suggests the interpretation of $\gamma_{1}$ as the beginning 
of the movement from $0$ to $\infty$ along a $q$-spiral,
or even a precession.

\subsubsection*{2.4.2~~~Irregular equations and
``infinitesimal'' elliptic curves}

The next step in our program is to tackle the irregular local
theory. The classification problem is solved in a common work
with Jean-Pierre Ramis and Changgui Zhang (\emph{see} \cite{RSZ}).
It uses a new discrete summation method for divergent solutions
(\emph{see} \cite{ZhangGroningen}), the existence of a canonical
filtration by the slopes for $q$-difference modules (\emph{see}
\cite{JSCRAS3}, \cite{JSAIF2} and \cite{JSGTQDIF}) and sheaf
theoretic methods due to Malgrange, Martinet and Ramis. 
The latter allow us to extend and enrich the geometric setting 
of the present paper. 
In \cite{RSZ}, we use a dynamical interpretation (due to
J.-L. Martins) of classical asymptotical developments. This 
version can be discretized in the following sense. While the
sheaf of Malgrange (\emph{see}, for instance, \cite{CanoRamis})
is defined on the \emph{horizon} 
$\mathbf{S}^{1} = \mathbf{C}^{*} / \Sigma$
of the action of the semigroup $\Sigma = ]0;1[$, 
taking instead the semigroup $\Sigma_{q} = q^{- \mathbf{N}}$,
we get the horizon $\mathbf{C}^{*} / \Sigma_{q}$, whence 
a sheaf defined on the elliptic curve $\mathbf{E}_{q}$.
Our vector bundles are related to this sheaf in the fuchsian
case. The elliptic curve $\mathbf{E}_{q}$ could also be viewed as the 
\emph{quotient of an infinitesimal neighborhood of $0$},
predicted by Ramis in \cite{RamisJPRTraum}. Actually, as shown in \cite{RSZ},
there is a whole family of such infinitesimal neighborhoods related to 
all the possible slopes and the corresponding sheaves of functions).
The results presented here can easily be extended to the
category of ``tamely irregular modules''; these are direct sums
of pure modules. One thus computes a Galois group $G_{mi}^{(0)}$.
Then, the graded functor $gr$ associated  to the canonical filtration
being faithful, exact and $\otimes$-preserving, one realizes the general local
Galois group $G_{i}^{(0)}$ as a semi-direct product of $G_{mi}^{(0)}$
by a pro-unipotent group, generated by the Stokes operators.
This will be detailed in \cite{JSIRR}.



\section{Global theory}


\subsection{The global Galois groupoid}


\subsubsection*{3.1.1~~~Birkhoff's classification revisited}

We shall give a galoisian meaning to Birkhoff's classification 
theorem. \\

\textbf{3.1.1.1 The category $\mathcal{C}$ of connection data. -}
We introduce a categorical variant $\mathcal{C}$ of the set 
of classifying data introduced in 1.3.3. The objects are triples 
$(A^{(0)},M,A^{(\infty)})$ where, for some $n \in \mathbf{N}^{*}$,
$A^{(0)},A^{(\infty)} \in Gl_{n}(\mathbf{C})$,
$M \in Gl_{n}(\mathcal{M}(\mathbf{C}^{*}))$ and moreover:
$$
(\sigma_{q} M) A^{(0)} = A^{(\infty)} M.
$$
Said otherwise, $M: A^{(0)} \rightarrow A^{(\infty)}$ 
is an isomorphism in $\mathcal{E}_{f}^{(*)}$. Morphisms
$(A^{(0)},M,A^{(\infty)}) \rightarrow (B^{(0)},N,B^{(\infty)})$
are pairs $(S^{(0)},S^{(\infty)})$ such that:
$$
\begin{cases}
S^{(0)}: A^{(0)} \rightarrow B^{(0)}
\text{ is a morphism in } \mathcal{E}_{f}^{(0)} 
\text{ or, what amounts to the same, in } \mathcal{P}^{(0)} \\
S^{(\infty)}: A^{(\infty)} \rightarrow B^{(\infty)}
\text{ is a morphism in } \mathcal{E}_{f}^{(\infty)}
\text{ or, what amounts to the same, in } \mathcal{P}^{(\infty)}
\end{cases}
$$
and, moreover, the following square commutes:
\begin{equation*}
\begin{CD}
A^{(0)}  @>M>> A^{(\infty)}  \\
@V{S^{(0)}}VV                          @VV{S^{(\infty)}}V      \\
B^{(0)}  @>N>> B^{(\infty)}
\end{CD}
\end{equation*}
In the same vein as in 2.3, objects of $\mathcal{C}$ 
can be interpreted as triples $(F^{(0)},f,F^{(\infty)})$ 
where $(F^{(0)},F^{(\infty)})$ are holomorphic vector bundles 
over $\mathbf{E}_{q}$ and $f: F^{(0)} \rightarrow F^{(\infty)}$ 
is a \emph{meromorphic} map between them. \\

Then, we make $\mathcal{C}$ a tensor category by endowing 
it with the natural (componentwise) tensor product;
here, we use the conventions of 1.1.2 for the tensor product
of matrices. The resulting category is plainly an abelian 
$\mathbf{C}$-linear neutral tannakian category. 
Moreover, after \cite{DM}, prop. 2.21, the projections
to $\mathcal{E}_{f}^{(0)}$ and $\mathcal{E}_{f}^{(\infty)}$
induce closed embeddings of $G^{(0)}$ and $G^{(\infty)}$
into the Galois group of $\mathcal{C}$. Our goal is 
to build an equivalence of $\mathcal{C}$ with $\mathcal{E}_{f}$. \\

\textbf{3.1.1.2 The category $\mathcal{S}$ of solutions. -}
In Birkhoff's method, one encodes a fuchsian equation
$\sigma X = A X$, by its local solutions at $0$ and $\infty$:
$X^{(0)} = M^{(0)} e_{q,A^{(0)}}$ and 
$X^{(\infty)} = M^{(\infty)} e_{q,A^{(\infty)}}$.
In more intrinsic terms, we shall use the local
flat forms $A^{(0)}$ and $A^{(\infty)}$ together
with the meromorphic gauge transformations $M^{(0)}$ 
(reducing $A$ to $A^{(0)}$) and $M^{(\infty)}$ 
(reducing $A$ to $A^{(\infty)}$). Due to the non canonicity 
of all these local data (except from the generic case 
of strictly fuchsian non resonant equations), we shall 
eventually map solutions to equations and not the other 
way round. \\

We shall therefore consider local pairs at $0$ and at
$\infty$:
\begin{eqnarray*}
(A^{(0)},M^{(0)}) & \in & 
GL_{n}(\mathcal{M}(\mathbf{C})) \times GL_{n}(\mathbf{C}) \\
(A^{(\infty)},M^{(\infty)}) & \in & 
GL_{n}(\mathcal{M}(\mathbf{C}_{\infty})) \times GL_{n}(\mathbf{C})
\end{eqnarray*}
We shall say that two such pairs are \emph{connected}
if one of the following (obviously) equivalent
conditions is realized:
\begin{enumerate}

\item{The following two expressions are equal:
$$
\left(\sigma M^{(\infty)}\right) A^{(\infty)}
\left(M^{(\infty)}\right)^{-1} = 
\left(\sigma M^{(0)}\right) A^{(0)}
\left(M^{(0)}\right)^{-1}
$$
}

\item{The matrix 
$M = \left(M^{(\infty)}\right)^{-1} M^{(0)}
\in GL_{n}(\mathcal{M}(\mathbf{C}^{*}))$
is such that $\left(\sigma M\right) A^{(0)} = A^{(\infty)} M$.}

\end{enumerate}
In this case, we shall call $A$ the common value of 
the two expressions appearing in the first condition. 
Being meromorphic on both $\mathbf{C}$ and $\mathbf{C}_{\infty}$,
it is meromorphic on $\mathbf{S}$, that is, rational:
$A \in GL_{n}(\mathbf{C}(z))$. Moreover, $M^{(0)}$
(resp. $M^{(\infty)}$) can be viewed as a morphism
from $A^{(0)}$ (resp. from $A^{(\infty)}$) to $A$. \\

We now define our category of solutions.
\begin{itemize}

\item{\textbf{Objects of $\mathcal{S}$. -} 
They are the quadruples:
$$
(A^{(0)},M^{(0)},A^{(\infty)},M^{(\infty)}) \in
Gl_{n}(\mathbf{C}) \times Gl_{n}(\mathcal{M}(\mathbf{C})) \times 
Gl_{n}(\mathbf{C}) \times Gl_{n}(\mathcal{M}(\mathbf{C}_{\infty}))
$$
such that the two component pairs $(A^{(0)},M^{(0)})$
and $(A^{(\infty)},M^{(\infty)})$ are connected.
}

\item{\textbf{Morphisms in $\mathcal{S}$. -} 
The morphisms from
$(A^{(0)},M^{(0)},A^{(\infty)},M^{(\infty)})$
to
$(B^{(0)},N^{(0)},B^{(\infty)},N^{(\infty)})$
in $\mathcal{S}$ are the triples $(F,S^{(0)},S^{(\infty)})$
such that 
$$
\begin{cases}
S^{(0)}: A^{(0)} \rightarrow B^{(0)}
\text{ is a morphism in } \mathcal{E}_{f}^{(0)} \\
S^{(\infty)}: A^{(\infty)} \rightarrow B^{(\infty)}
\text{ is a morphism in } \mathcal{E}_{f}^{(\infty)} \\
\end{cases}
$$
and, moreover, the following squares commute:
\begin{equation*}
\begin{CD}
A^{(0)}      @>M^{(0)}>>   A   @<M^{(\infty)}<<     A^{(\infty)}  \\
@V{S^{(0)}}VV           @V{F}VV                @VV{S^{(\infty)}}V \\
B^{(0)}      @>N^{(0)}>>   B   @<N^{(\infty)}<<     B^{(\infty)}
\end{CD}
\end{equation*}
Here, $A$ and $B$ are defined according to our
previous convention. One then notes, first that
\begin{eqnarray*}
F & = & N^{(\infty)} S^{(\infty)} \left(M^{(\infty)}\right)^{-1} \\
  & = & N^{(0)} S^{(0)} \left(M^{(0)}\right)^{-1}
\end{eqnarray*}
is meromorphic on both $\mathbf{C}$ and $\mathbf{C}_{\infty}$,
therefore rational; second, that, by any of these two
expressions for $F$, $\left(\sigma F \right) A = BF$,
that is, $F$ is a morphism $A \rightarrow B$
in $\mathcal{E}_{f}$.
}

\item{\textbf{Tensor structure on $\mathcal{S}$. -} 
The tensor product of objects (resp. morphisms) is defined
componentwise on the quadruples (resp. triples), using the
usual identifications. }

\end{itemize}
Again, one has obtained an abelian $\mathbf{C}$-linear 
neutral tannakian category such that the projections
to $\mathcal{E}_{f}^{(0)}$ and $\mathcal{E}_{f}^{(\infty)}$
induce closed embeddings of $G^{(0)}$ and $G^{(\infty)}$
into the Galois group of $\mathcal{S}$. \\

\textbf{3.1.1.3 Proposition: equivalence of the tensor categories 
$\mathcal{E}_{f}$, $\mathcal{S}$ and $\mathcal{C}$. -}
\emph{One keeps the previous conventions for $A$ and $M$.
Then, taking:
$$
\begin{cases}
(A^{(0)},M^{(0)},A^{(\infty)},M^{(\infty)}) \leadsto A \\
(F,S^{(0)},S^{(\infty)}) \leadsto F
\end{cases}
$$
and
$$
\begin{cases}
(A^{(0)},M^{(0)},A^{(\infty)},M^{(\infty)}) \leadsto
(A^{(0)},M,A^{(\infty)}) \\
(F,S^{(0)},S^{(\infty)}) \leadsto (S^{(0)},S^{(\infty)})
\end{cases}
$$
provides us with $\mathbf{C}$-linear $\otimes$-equivalences
from $\mathcal{S}$ to $\mathcal{E}_{f}$ and 
from $\mathcal{S}$ to $\mathcal{C}$.} \\

It is clear that one has defined two $\otimes$-functors
and that the first one is fully faithful. That it is 
also essentially surjective stems from the existence of flat
local reductions at $0$ and at $\infty$ for any fuchsian 
equation, as recalled in 1.2. \\

That the second functor is fully faithful comes from
the fact that the equalities:
\begin{eqnarray*}
F & = & N^{(\infty)} S^{(\infty)} \left(M^{(\infty)}\right)^{-1} \\
  & = & N^{(0)} S^{(0)} \left(M^{(0)}\right)^{-1}
\end{eqnarray*}
give a unique antecedent to a pair $(S^{(0)},S^{(\infty)})$.
For essential surjectivity, one starts from an object
$(A^{(0)},M,A^{(\infty)})$ of $\mathcal{C}$. Since
$M \in GL_{n}(\mathcal{M}(\mathbf{C}^{*}))$, Birkhoff's
lemma (\emph{see} \cite{JSAIF}, 2.2 and 2.3) allows us
to write:
$$
M = \left(M^{(\infty)}\right)^{-1} M^{(0)},
\text{ where } (M^{(0)},M^{(\infty)}) \in
Gl_{n}(\mathcal{M}(\mathbf{C})) \times Gl_{n}(\mathcal{M}(\mathbf{C}_{\infty})).
$$
It is then clear that 
$(A^{(0)},M^{(0)},A^{(\infty)},M^{(\infty)})$
is a connected quadruple and an antecedent of
$(A^{(0)},M,A^{(\infty)})$.
\hfill $\Box$. \\

\textbf{3.1.1.4 Singularities and exponents. -}
Let $\Sigma$ be a finite subset of $\mathbf{C}^{*}$.
We shall have to consider the full subcategory
$\mathcal{C}_{\Sigma}$ of $\mathcal{C}$ whose
objects are the triples $(A^{(0)},M,A^{(\infty)})$
such that $\mathcal{S}(M) \subset q^{\mathbf{Z}} \Sigma$.
It is stable by all tensor and abelian constructions,
hence a tannakian subcategory. From 1.2.5 and from the remarks 
in \emph{loc. cit.}, 2.3.1, one draws that the equivalence
shown in 3.1.1.3 induces an equivalence of the tannakian
categories $\mathcal{E}_{f,\Sigma}$ and $\mathcal{C}_{\Sigma}$.


\subsubsection*{3.1.2~~~The global Galois groupoid and 
the global fundamental groupoid}

Composing the above projections with fibre functors
$\omega^{(0)}_{z_{0}},\omega^{(\infty)}_{z_{0}}$ 
provides us with two fibre functors on $\mathcal{C}$. 
We shall call these restrictions by the same names. \\

\textbf{3.1.2.1 Definition. -} 
\emph{The Galois groupoid of $\mathcal{C}$ 
is the groupoid $G$ having as base set:
$$
\{\omega^{(0)}_{a} \;/\; a \in \mathbf{C}^{*}\} \cup
\{\omega^{(\infty)}_{a} \;/\; a \in \mathbf{C}^{*}\},
$$
and such that, for any two $a,b \in \mathbf{C}^{*}$:
\begin{eqnarray*}
G(\omega^{(0)}_{a},\omega^{(0)}_{b}) & = &
Iso^{\otimes}(\omega^{(0)}_{a},\omega^{(0)}_{b}) \\
G(\omega^{(\infty)}_{a},\omega^{(\infty)}_{b}) & = &
Iso^{\otimes}(\omega^{(\infty)}_{a},\omega^{(\infty)}_{b}) \\
G(\omega^{(0)}_{a},\omega^{(\infty)}_{b}) & = &
\begin{cases}
\text{if } a = b \;:\;
Iso^{\otimes}(\omega^{(0)}_{a},\omega^{(\infty)}_{b}) \\
\text{if } a \not= b \;:\; \emptyset
\end{cases}
\end{eqnarray*}
}

The local groupoids computed in chapter 2 embed in the
corresponding subgroupoids of $G$, giving many elements in
the groups $Aut^{\otimes}(\omega^{(0)}_{z_{0}})$ and 
$Aut^{\otimes}(\omega^{(\infty)}_{z_{0}})$ for all
$z_{0} \in \mathbf{C}^{*}$ and in the sets
$Iso^{\otimes}(\omega^{(0)}_{z_{0}},\omega^{(0)}_{z_{1}})$ 
and 
$Iso^{\otimes}(\omega^{(\infty)}_{z_{0}},\omega^{(\infty)}_{z_{1}})$
for all $z_{0},z_{1} \in \mathbf{C}^{*}$. To complete this
and connect $G$, we want to build sufficiently many elements
in the sets
$Iso^{\otimes}(\omega^{(0)}_{z_{0}},\omega^{(\infty)}_{z_{0}})$
for all $z_{0} \in \mathbf{C}^{*}$. For instance, one gets 
such a $\otimes$-isomorphism from $\omega^{(0)}_{z_{0}}$
to $\omega^{(\infty)}_{z_{0}}$ by taking $(A^{(0)},M,A^{(\infty)})$
to $M$. But this is not defined over $\mathbf{C}^{*}$, so we change 
our way. This can be done by evaluating
$M$ at points $z_{0} \not\in \mathcal{S}(M)$. For that,
we fix a finite subset $\Sigma$ of $\mathbf{C}^{*}$
and restrict to the full subcategory $\mathcal{C}_{\Sigma}$
of $\mathcal{C}$. \\

\textbf{3.1.2.2 Proposition. -}
\emph{For any such point $z_{0}$, the natural transformation
$\Gamma_{z_{0}} : (A^{(0)},M,A^{(\infty)}) \leadsto M(z_{0})$ 
is an element of 
$Iso^{\otimes}(\omega^{(0)}_{z_{0}},\omega^{(\infty)}_{z_{0}})$.} \\

Here and in the following, we keep the same names for the restrictions
to $\mathcal{C}_{\Sigma}$ of the fibre functors $\omega^{(0)}_{z_{0}}$,
$\omega^{(\infty)}_{z_{0}}$ ($z_{0} \not\in q^{\mathbf{Z}} \Sigma$). 
The proof of the proposition is then more or less tautological. 
Tensor preservation comes from the definition of the tensor structure 
componentwise, plus the obvious fact that 
$\left(M \otimes M'\right)(z_{0}) = M(z_{0}) \otimes M'(z_{0})$.
Functoriality comes from the commuting square in the
definition of morphisms in 3.1.1.1 plus the obvious 
computation:
\begin{eqnarray*}
N(z_{0}) S^{(0)}(z_{0}) & = & \left(N S^{(0)}\right)(z_{0}) \\
                        & = & \left(S^{(\infty)} M\right)(z_{0}) \\
                        & = & S^{(\infty)}(z_{0}) M(z_{0}) \\
\end{eqnarray*}
\hfill $\Box$ \\

We see $\Gamma_{z_{0}}$ as a path connecting the points
$\omega^{(0)}_{z_{0}},\omega^{(\infty)}_{z_{0}}$
of the groupoid $G$. \\

\textbf{3.1.2.3 Theorem. -}
\emph{The local groupoids at $0$ and at $\infty$ 
(defined and computed in chapter 2) together with the
paths $\Gamma_{z_{0}} \;,\; z_{0} \not\in \Sigma$
generate a Zariski-dense subgroupoid of the Galois groupoid
of $\mathcal{C}_{\Sigma}$.} \\

We appeal again to the criterion of Chevalley (see the proof of 2.2.3.3). 
It can easily be extended to the case 
of a groupoid in the following way. We choose an object 
$\mathcal{X} = (A^{(0)},M,A^{(\infty)})$, and, for each
basepoint $\omega^{(0)}_{a}$ (resp. $\omega^{(\infty)}_{a}$),
$a \not\in q^{\mathbf{Z}} \Sigma$, a line 
$D^{(0)}_{a} \subset \omega^{(0)}_{a}(\mathcal{X})$
(resp.
$D^{(\infty)}_{a} \subset \omega^{(\infty)}_{a}(\mathcal{X})$)
and we assume that this family of lines is globally stable under 
the action of $G^{(0)}$, $G^{(\infty)}$ and our special paths.
It is then sufficient to check that this family of lines is 
actually stable under the action of the whole Galois groupoid. \\

By Tannaka duality for the category $\mathcal{E}_{f}^{(0)}$
and for the groupoid $G^{(0)}$, we see that the family of
lines $D^{(0)}_{a}$ induces a subrepresentation of rank $1$
of the representation defined by the object $A^{(0)}$, hence
comes from a subobject of rank $1$ of $A^{(0)}$. 
This subobject is an injection 
$x^{(0)}: a^{(0)} \rightarrow A^{(0)}$ 
in $\mathcal{E}_{f}^{(0)}$ and we may take it to lie
in $\mathcal{P}^{(0)}$. This means that 
$a^{0} \in \mathbf{C}^{*}$, that $x^{(0)}$ is
a function holomorphic on $\mathbf{C}^{*}$ and that,
for all $a \not\in q^{\mathbf{Z}} \Sigma$, the line
$D^{(0)}_{a}$ is the image of the linear map
$\omega^{(0)}_{a}(x^{(0)})$, that is:
$D^{(0)}_{a} = \mathbf{C} x^{(0)}(a)$. \\

The same argument on the $\infty$ side shows that there
exists a subobject 
$x^{(\infty)}: a^{(\infty)} \rightarrow A^{(\infty)}$ 
such that, for all $a \not\in q^{\mathbf{Z}} \Sigma$,
$D^{(\infty)}_{a} = \mathbf{C} x^{(\infty)}(a)$.
The condition of stability under our special paths
says that 
$\forall a \not\in q^{\mathbf{Z}} \Sigma \;,\;
M(a) D^{(0)}_{a} = D^{(\infty)}_{a}$, so that,
out of $q^{\mathbf{Z}} \Sigma$, there exists 
a holomorphic function $m$ such that
$M(a) x^{(0)}(a) = m(a) x^{(\infty)}(a)$. 
This amounts to say that $\phi = (x^{(0)},x^{(\infty)})$
is a morphism from $\mathcal{X'} = (a^{0},m,a^{(\infty)})$ 
to $\mathcal{X} = (A^{(0)},M,A^{(\infty)})$ in $\mathcal{C}$. \\

We now take an arbitrary galoisian isomorphism, that is, an element 
$h \in Iso^{\otimes}(\omega_{a}^{(0)},\omega_{a}^{(\infty)})$.
The functoriality condition gives rise to a commutative diagram:
\begin{equation*}
\begin{CD}
\omega_{a}^{(0)}(\mathcal{X}') @>h(\mathcal{X}')>> \omega_{a}^{(\infty)}(\mathcal{X}') \\
@V{\omega_{a}^{(0)}(\phi)}VV                     @VV{\omega_{a}^{(\infty)}(\phi)}V      \\
\omega_{a}^{(0)}(\mathcal{X})  @>h(\mathcal{X})>>  \omega_{a}^{(\infty)}(\mathcal{X}) 
\end{CD}
\end{equation*}

Then $h(\mathcal{X})x^{(0)}(a)  = h(\mathcal{X}') x^{(\infty)}(a)$ ;
since $h(\mathcal{X}') \in \mathbf{C}^{*}$, 
this shows the stability of our family of lines
under the action of $Iso^{\otimes}(\omega^{(0)},\omega^{(\infty)})$,
hence also under the action of the whole Galois groupoid.
\hfill $\Box$ 


\subsubsection*{3.1.3~~~The case of a regular equation}

\textbf{3.1.3.1 Regular triples. -}
Let us now consider the case of a \emph{regular} equation.
Recall from chapter 1 that $A$ is said to be regular at $0$
(resp. at $\infty$) if it is equivalent to the identity matrix
$I_{n}$ at $0$ (resp. at $\infty$). Then $A$ has a local
reduction at $0$, $M^{(0)}: I_{n} \rightarrow A$ as well 
as a local reduction at $\infty$, 
$M^{(\infty)}: I_{n} \rightarrow A$. One can therefore
associate to $A$ the triple $(I_{n},M,I_{n})$ in
$\mathcal{C}_{\Sigma}$, with $\Sigma = \mathcal{S}(A)$,
$M = M = \left(M^{(\infty)}\right)^{-1} M^{(0)}$.
Moreover, in this case, $M$ is ``the'' connection matrix $P$ 
and it is elliptic. \\

\textbf{3.1.3.2 Corollary: the Galois group of a regular equation. -}
\emph{The Galois group at any point $\omega^{(0)}_{z_{0}}$
is the Zariski closure of the subgroup generated by the values 
$\left(P(a)\right)^{-1}P(b)$
for $a,b \not\in q^{\mathbf{Z}} \Sigma$.} \\

Indeed, from the equalities $A^{(0)} = A^{(\infty)} = I_{n}$,
we draw that the local Galois groupoids at $0$ and at $\infty$
of the equation are trivial (there are only identity arrows
between any two points). The conclusion now follows from
theorem 3.1.2.3. 
\hfill $\Box$ \\

This is the case tackled by Etingof in \cite{Etingof}, and this 
proposition is his main result.


\subsection{Birkhoff's method}

We shall follow here Birkhoff's method more closely, using
the connection matrix $P$ itself (together with local
linear data) to encode fuchsian equations, then, trying
to interpret its values as monodromy data. However,
the bad multiplicative properties of any canonical choice 
of solutions, hence of the matrix $P$ itself, lead us to twist 
first the tensor product in the category of connection data, 
second the connection matrix itself into a matrix $\breve{P}$ in
order to get galoisian properties. The relation with the point of view
of 3.1 is explained in 3.2.3. Proofs and details can be found
in \cite{JSthese} and \cite{JSGAL}.


\subsubsection*{3.2.1~~~Equivalences of tannakian categories}

We encode a fuchsian equation $\sigma X = A X$ 
by its local solutions at $0$ and $\infty$:
$X^{(0)} = M^{(0)} e_{q,A(0)}$ and 
$X^{(\infty)} = M^{(\infty)} e_{q,A(\infty)}$
and its connection matrix 
$P = \left(X^{(\infty)}\right)^{-1}X^{(0)}$,
as defined in 1.3.3. We shall consistently use
these notations herebelow, without further
notice. These data are not unique,
so that we use an intermediate category of solutions
to link equations and connection triples. 
By necessity, one does not take the natural tensor
product $X_{1}^{(0)} \otimes X_{2}^{(0)}$ on solutions.
This comes from the fact that 
$e_{q,A} \otimes e_{q,B} \not= e_{q,A \otimes B}$, thereby
destroying our normal forms for solutions. The defect of
equality has been analyzed in 1.2.2.3. We thus give a special 
notation for the twisted tensor product:
$$
X_{1}^{(0)} \underline{\otimes} X_{2}^{(0)} =
\left(M_{1}^{(0)} \otimes M_{2}^{(0)}\right) 
e_{q,A_{1}^{(0)} \otimes A_{2}^{(0)}},
$$
and similarly at $\infty$. This is related to the natural 
tensor product in the following way:
$$
X_{1}^{(0)} \underline{\otimes} X_{2}^{(0)} =
\left(X_{1}^{(0)} \otimes X_{2}^{(0)}\right)
\Phi(A_{1}^{(0)},A_{2}^{(0)}).
$$
By necessity, one is thus led to twist the natural tensor product
of connection matrices:
$$
P_{1} \otimes P_{2} =
\left(X_{1}^{(\infty)} \otimes X_{2}^{(\infty)}\right)^{-1}
\left(X_{1}^{(0)} \otimes X_{2}^{(0)}\right)
$$
in the following way:
$$
P_{1} \Otimes P_{2} =
\Phi(A_{1}^{(\infty)} , A_{2}^{(\infty)})
\; (P_{1} \otimes P_{2}) \;
\left(\Phi(A_{1}^{(0)} , A_{2}^{(0)})\right)^{-1}.
$$
This notation is slightly ambiguous, since the right hand side
does not really depend on $P_{1}$ and $P_{2}$ alone, but also involves
the local linear data 
$A_{1}^{(0)},A_{2}^{(0)},A_{1}^{(\infty)} , A_{2}^{(\infty)}$.
Note that, in the case that one of $A_{1}^{(0)},A_{2}^{(0)}$ is unipotent 
and the same holds at $\infty$, we have 
$P_{1} \Otimes P_{2} = P_{1} \otimes P_{2}$. This is the case,
e.g. for \emph{regular} equations (i.e. such that the matrices at
$0$ and $\infty$ are $I_{n}$). \\

\textbf{3.2.1.1 The category $\mathcal{S}'$ of solutions. -}
Its objects are quadruples
$$
(A^{(0)},M^{(0)},A^{(\infty)},M^{(\infty)}) \in
Gl_{n}(\mathbf{C}) \times Gl_{n}(\mathcal{M}(\mathbf{C})) \times 
Gl_{n}(\mathbf{C}) \times Gl_{n}(\mathcal{M}(\mathbf{C}_{\infty}))
$$
such that $X^{(0)} = M^{(0)} e_{q,A(0)}$ and 
$X^{(\infty)} = M^{(\infty)} e_{q,A(\infty)}$ are connected
in a sense similar to 3.1.1.2. One can likewise adapt the
definitions in such a way as to get a neutral tannakian category. \\

\textbf{3.2.1.2 The category $\mathcal{C}'$ of connection data. -}
Birkhoff's classification theorem
(\emph{see} \cite{Birkhoff1}, \cite{JSAIF}) amounts to
saying that the data $(A^{(0)},P,A^{(\infty)})$ are enough
to recover $A$ up to rational equivalence. We shall now
give of it a categorical formulation.
\begin{itemize}

\item{\textbf{Objects of $\mathcal{C}'$. -}
They are the triples:
$$
(A^{(0)},P,A^{(\infty)}) \in 
Gl_{n}(\mathbf{C}) \times Gl_{n}(\mathcal{M}(\mathbf{E}_{q})) 
\times Gl_{n}(\mathbf{C}).
$$
}

\item{\textbf{Morphisms in $\mathcal{C}'$. -}
The morphisms from the object $(A^{(0)},P,A^{(\infty)})$ 
of order $n$ to the object $(B^{(0)},Q,B^{(\infty)})$ of order $p$  
are the pairs
$$
(R^{(0)},R^{(\infty)}) \in M_{p,n}(\mathbf{C}) \times M_{p,n}(\mathbf{C})
$$
such that 
$$
\begin{cases}
R^{(0)} \overline{A^{(0)}} = \overline{B^{(0)}} R^{(0)} \\
R^{(\infty)} P = Q R^{(0)} \\
R^{(\infty)} \overline{A^{(\infty)}} = 
\overline{B^{(\infty)}} R^{(\infty)}
\end{cases}
$$
This can be justified by the properties proved in 1.2.4
(for more details, \emph{see} \cite{JSthese}).
Note the following consequences of the definition:
$$
\begin{cases}
R^{(0)} \overline{A_{s}^{(0)}} = \overline{B_{s}^{(0)}} R^{(0)} \\
R^{(\infty)} \overline{A_{s}^{(\infty)}} = 
\overline{B_{s}^{(\infty)}} R^{(\infty)}
\end{cases}
\quad
\begin{cases}
R^{(0)} A_{u}^{(0)} = B_{u}^{(0)} R^{(0)} \\
R^{(\infty)} A_{u}^{(\infty)} = B_{u}^{(\infty)} R^{(\infty)}
\end{cases}
$$
}

\item{\textbf{Tensor structure on $\mathcal{C}'$. -}
The tensor product of two objects 
$(A_{1}^{(0)},P_{1},A_{1}^{(\infty)})$ and 
$(A_{2}^{(0)},P_{2},A_{2}^{(\infty)})$
is defined to be
$$
(A_{1}^{(0)},P_{1},A_{1}^{(\infty)}) \otimes 
(A_{2}^{(0)},P_{2},A_{2}^{(\infty)}) =
(A_{1}^{(0)} \otimes A_{2}^{(0)},
P_{1} \Otimes P_{2},
A_{1}^{(\infty)} \otimes A_{2}^{(\infty)}).
$$
The tensor product of two morphisms
$(R_{1}^{(0)},R_{1}^{(\infty)})$ and 
$(R_{2}^{(0)},R_{2}^{(\infty)})$
is defined componentwise, from the usual tensor product.
That the tensor product of two objects is one is obvious ; that
the tensor product of two morphisms is one is not tautological,
but follows from the properties of morphisms in $\mathcal{S}'$
and $\mathcal{C}'$.}

\end{itemize}
Defining as before two fibre functors 
$\mathcal{C}' \rightarrow Vect_{\mathbf{C}}^{f}$
by sending $(A^{(0)},P,A^{(\infty)})$
to the $\mathbf{C}^{n}$ underlying $A^{(0)}$ (resp. $A^{(\infty)}$)
and $(S^{(0)},S^{(\infty)})$ to $S^{(0)}$ (resp. $S^{(\infty)}$),
one obtains again a neutral tannakian category
over $\mathbf{C}$. \\

\textbf{3.2.1.3 The equivalence of $\mathcal{E}_{f}$, $\mathcal{S}'$ 
and  $\mathcal{C}'$. -} 
As in 3.1, two functors can be define; first, one from solutions 
to equations:
$$
\begin{cases}
(A^{(0)},M^{(0)},A^{(\infty)},M^{(\infty)}) \leadsto A \\
(F,S^{(0)},S^{(\infty)}) \leadsto F
\end{cases}
$$
Next, one from solutions to connection triples:
$$
\begin{cases}
(A^{(0)},M^{(0)},A^{(\infty)},M^{(\infty)}) \leadsto
(A^{(0)},P,A^{(\infty)}) \\
(F,S^{(0)},S^{(\infty)}) \leadsto (S^{(0)},S^{(\infty)})
\end{cases}
$$
Both are exact $\otimes$-preserving $\mathbf{C}$-linear 
equivalence of categories.
Note that this equivalence is compatible with
the fibre functors previously introduced. \\

\textbf{3.2.1.4 Singularities and exponents. -}
In order to compute the connection component of the Galois 
groupoid we shall need to evaluate the connection matrix at 
various points of $\mathbf{C}^{*}$, avoiding its singularities 
(since $P$ is elliptic, these singularities may actually be seen 
as points in $\mathbf{E}_{q}$). For an object $(A^{(0)},P,A^{(\infty)})$
of $\mathcal{C}'$ coming from an object $A$ of $\mathcal{E}_{f}$, 
the singularities of $P$ respectively come from:
\begin{enumerate}

\item{$\mathcal{S}(A)$ for the meromorphic part 
$M = (M^{(\infty)})^{-1} M^{(0)}$. Precisely, 
$\mathcal{S}(M) \subset q^{\mathbf{Z}} \mathcal{S}(A)$.}

\item{$Sp(A(0)),Sp(A(\infty))$ for the semi-simple components 
of the $e_{q,A(-)}$ parts. Since reduction to a constant matrix 
is not unique, these exponents of $A$ are actually defined
up to a factor in $q^{\mathbf{Z}}$.}

\item{$l_{q}$, the $q$-logarithm, in case $A(0)$ or $A(\infty)$
is not semi-simple.}

\end{enumerate}
On the other hand, we want to define $\otimes$-stable
categories only, so that we have to allow for multiplication 
of the exponents (and the inverse, to have stability under
dual taking). Thus, for $\Sigma$ a finite subset of the open set 
$\mathbf{C}^{*}$ and $C$ be a finitely generated subgroup of 
the group $\mathbf{C}^{*}$, we consider the full subcategory 
$\mathcal{E}_{f,\Sigma,C}$ of $\mathcal{E}_{f}$
whose objects have all singularities in $q^{\mathbf{Z}} \Sigma$ 
and all exponents in $q^{\mathbf{Z}} C$. From the precised
version of Birkhoff's lemma we draw that an equation $A$ 
is in $\mathcal{E}_{f,\Sigma,C}$ if and only if its exponents
are in $C$ and $\mathcal{S}(M) \subset q^{\mathbf{Z}} \Sigma$. \\

Write $\overline{\Sigma}$, resp. $\overline{C}$ for the image 
in $\mathbf{E}_{q}$ of a finite subset $\Sigma$ of $\mathbf{C}^{*}$, 
resp. a finitely generated subgroup $C$ of $\mathbf{C}^{*}$.
We then consider the full subcategory $\mathcal{C}'_{\Sigma,C}$
of $\mathcal{C}'$ whose objects are the triples 
$(A^{(0)},P,A^{(\infty)})$ such that $Sp(A(0)),Sp(A(\infty))$
are subsets of $q^{\mathbf{Z}} C$ and 
$\mathcal{S}(M) \subset q^{\mathbf{Z}} \Sigma$. 
For such objects,
$\mathcal{S}(P) \subset \overline{\Sigma} \cup \overline{C}$.
Moreover, $\mathcal{C}'_{\Sigma,C}$ is a strictly full tannakian
subcategory of $\mathcal{C}'$ and it is equivalent to 
$\mathcal{E}_{f,\Sigma,C}$. \\

For an object $\mathcal{X} = (A^{(0)},P,A^{(\infty)})$ 
of $\mathcal{C}'$ , denote by $<\mathcal{X}>$ the tannakian 
subcategory generated by $\mathcal{X}$. If $\mathcal{X}$ actually 
belongs to the subcategory $\mathcal{C}'_{\Sigma,C}$, 
this entails $<\mathcal{X}> \subset \mathcal{C}'_{\Sigma,C}$. 
The object being given, he minimal choice is to take for $C$ 
the subgroup of $\mathbf{C}^{*}$ generated by $Sp(A(0))$ and 
$Sp(A(\infty))$ and for $\Sigma$ the singular locus of $M$.


\subsubsection*{3.2.2~~~The Galois groupoid}

Sticking to the previous definitions, we consider a groupoid $G$
with base points $0$ and $\infty$ and with corresponding 
arrow sets $Aut^{\otimes}(\omega^{(0)})$,
$Aut^{\otimes}(\omega^{(\infty)})$,
$Iso^{\otimes}(\omega^{(0)},\omega^{(\infty)})$,
$Iso^{\otimes}(\omega^{(\infty)},\omega^{(0)})$.

\xymatrix{
& & & & & &
\omega^{(0)} \ar@<1ex>[rr]^-{Iso^{\otimes}(\omega^{(0)},\omega^{(\infty)})} 
\ar@(dl,ul)[]^-{\txt{$Aut^{\otimes}(\omega^{(0)})$}}
& &
\omega^{(\infty)} \ar@<1ex>[ll]^-{Iso^{\otimes}(\omega^{(\infty)},\omega^{(0)})}
\ar@(ur,dr)[]^-{\txt{$Aut^{\otimes}(\omega^{(\infty)})$}}
}

In this section, we shall build elements of the Galois groupoid,
that is, $\otimes$-automorphisms of $\omega^{(0)}$ and of 
$\omega^{(\infty)}$ and $\otimes$-isomorphisms from $\omega^{(0)}$ 
to $\omega^{(\infty)}$. We shall consistently denote
$\mathcal{X} = (A^{(0)},P,A^{(\infty)})$ 
a generic object of $\mathcal{C}'$. We then write
$A^{(0)} = A_{s}^{(0)} A_{u}^{(0)}$ and
$A^{(\infty)} = A_{s}^{(\infty)} A_{u}^{(\infty)}$
the Dunford decompositions. \\

\textbf{3.2.2.1 Local automorphisms of the fibre functor. -}
>From the general facts recalled at the end of the introduction,
we easily deduce the following:
\begin{enumerate}

\item{Let $f$ be a map: $\mathbf{C}^{*} \rightarrow \mathbf{C}^{*}$. 
Then $\mathcal{X} \leadsto f(\overline{A_{s}^{(0)}})$
(resp. $\mathcal{X} \leadsto f(\overline{A_{s}^{(\infty)}})$)
is an automorphism of $\omega^{(0)}$
(resp. of $\omega^{(\infty)}$). If we take
$f \in Hom_{grp}(\mathbf{C^{*}},\mathbf{C^{*}})$
and (by necessity) such that $f(q) = 1$, we get a $\otimes$-compatible 
automorphism.}

\item{Let $\lambda \in \mathbf{C}$. Then 
$\mathcal{X} \leadsto (A_{u}^{(0)})^{\lambda}$
(resp. $\mathcal{X} \leadsto (A_{u}^{(\infty)})^{\lambda}$)
is a $\otimes$-automorphism of $\omega^{(0)}$
(resp. of $\omega^{(\infty)}$).}

\end{enumerate}
We thus obtain subgroups
$G^{(0)} \subset Aut^{\otimes}(\omega^{(0)})$ and
$G^{(\infty)} \subset Aut^{\otimes}(\omega^{(\infty)})$.
We recognize the local Galois groups found in chapter 2.
They are isomorphic to each other and are commutative 
proalgebraic groups with unipotent component $\mathbf{C}$ 
and semi-simple component:
$$
G_{s}^{(0)} \simeq G_{s}^{(\infty)} \simeq
\{f \in Hom_{grp}(\mathbf{C^{*}},\mathbf{C^{*}}) \;/\; f(q) = 1\} \simeq
Hom_{grp}(\mathbf{C^{*}}/q^{\mathbf{Z}},\mathbf{C^{*}}).
$$
This is just the dual $\check{\mathbf{E}_{q}}$
of the abstract group $\mathbf{E}_{q}$.
In this description, our local Galois groups are identified to
a subgroup of $\mathbf{Z}^{alg}$ but they are there embedded
transversally to the natural monodromy group $\mathbf{Z}$:
their intersection with the latter is the trivial subgroup. \\

\textbf{3.2.2.2 Building elements of the connection component. -}
We restrict here the fibre functors $\omega^{(0)},\omega^{(\infty)}$
to some subcategory $\mathcal{C}'_{\Sigma,C}$ (see 3.2.1.4). 
We put $\Sigma' = q^{\mathbf{Z}}(\Sigma \cup C)$ and fix
$a \in \mathbf{C}^{*} - \Sigma'$.
It stems tautologically from our definition of morphisms
in the category $\mathcal{C}'$ that, for any such $a$,
$\mathcal{X} \leadsto P(a)$ is a functorial isomorphism
$\omega^{(0)} \rightarrow \omega^{(\infty)}$.
However, it is not, in general, a \emph{tensor} isomorphism,
because $P_{1}(a) \otimes P_{2}(a) \not= P_{1}(a) \Otimes P_{2}(a)$.
There is of course an exception if $C$ is trivial,
e.g. for regular equations. The right and left excess factors
are $\Phi(A_{1,s}^{(0)},A_{2,s}^{(0)})$ 
and $\Phi(A_{1,s}^{(\infty)},A_{2,s}^{(\infty)})$. 
They can be exactly compensated by taking
$e_{q,A_{s}^{(\infty)}} P(a) (e_{q,A_{s}^{(0)}})^{-1}$
instead of $P(a)$. However, this depends on $A_{s}^{(0)}$
and $A_{s}^{(\infty)}$ and not only on
$\overline{A_{s}^{(0)}}$ and $\overline{A_{s}^{(\infty)}}$
(see the first half of the first fact in 3.2.2.1),
so that it is no longer a functorial isomorphism. \\

In order to twist the connection matrix, one chooses,
for each $a \in \mathbf{C}^{*}$, a group homomorphism
$g_{a} \in Hom_{grp}(\mathbf{C^{*}},\mathbf{C^{*}})$
such that $g_{a}(q) = 1$. We have exhibited 
such a family $(g_{a})_{a \in \mathbf{C}^{*}}$ in 2.2.3
and we shall make this choice more precise further below.
One then puts, for $c \in \mathbf{C}^{*}$,
$\psi_{a}(c) = \frac{e_{c}(a)}{g_{a}(c)}$, so that:
$$
\begin{cases}
\psi_{a}(c_{1}) \psi_{a}(c_{2}) = \psi_{a}(c_{1}c_{2}) \phi(c_{1},c_{2})(a) \\
\psi_{a}(c) \text{ depends only on } \overline{c}
\end{cases}
$$
Now extend $\psi_{a}$ to matrices, so that:
$$
\begin{cases}
\psi_{a}(A_{1}) \otimes \psi_{a}(A_{2}) = 
\left(\psi_{a}(A_{1}) \Otimes \psi_{a}(A_{2})\right) \Phi(A_{1},A_{2})(a) \\
\psi_{a}(A) \text{ depends only on } \overline{A}
\text{ (actually, on $\overline{A_{s}}$)} 
\end{cases}
$$
We have built our twisting factor. 
It is made up of two ingredients: 
one is due to the twisting of the tensor product, 
itself due to the non canonical choice of solutions.
The other comes from the artificial concentration
of the local groupoid at $0$, in a unique base point
(and a unique local group), as shown in the figure at the end of 3.2
and in proposition 3.2.3.1 and remark 3.2.3.2
(all the points $\omega^{(0)}_{a}$ are artificially concentrated
at the unique point $\omega^{(0)}$).
Let $a \in \mathbf{C}^{*} - \Sigma'$. Then
$$
\mathcal{X} \leadsto \breve{P}(a) =
\left(\psi_{a}(A_{s}^{(\infty)})\right)^{-1} P(a) 
\left(\psi_{a}(A_{s}^{(0)})\right)
$$
is a $\otimes$-isomorphism:
$\omega^{(0)} \rightarrow \omega^{(\infty)}$.
Note that choosing another group homomorphism $g_{a}$ changes it
by a factor in $\check{\mathbf{E}_{q}}$ and therefore changes
our twisted connection matrix $\breve{P}$
by a left and a right factor in the semi-simple component of the local
Galois groups. Similarly, it does not matter that we have taken
the same $g_{a}$ to twist on the $0$ and on the $\infty$ side.
We therefore take a slightly different choice for this family.
We first write $\mathbf{C}^{*} = \mathbf{U} \times q^{\mathbf{R}}$,
thereby meaning that we have chosen a logarithm of $q$:
$q = e^{- 2 \imath \pi \tau} \;,\; Im(\tau) \gt 0$. 
We thus write every $z \in \mathbf{C}^{*}$: 
$z = u q^{y} = u e^{- 2 \imath \pi \tau y}$ with $|u| = 1$
and $y \in \mathbf{R}$. This allows us to define, for any
$\alpha \in \mathbf{C}$, a group homomorphism:
$$
\delta_{\alpha}:
\begin{cases}
\mathbf{C}^{*} \rightarrow \mathbf{C}^{*} \\
u q^{y} \mapsto q^{\alpha y} = e^{- 2 \imath \pi \tau \alpha y}
\end{cases}
$$
To define $g_{a}$, we now choose a logarithm of $a$. We first
define the function $\log_{q}$ on the whole of $\mathbf{C}^{*}$
by the following conditions: it is to be holomorphic on 
$\mathbf{C}^{*} - q^{\mathbf{R}}$, one has $\log_{q}(q^{y}) = y$ 
and the discontinuity is just \emph{before} the cut when turning 
counterclockwise around $0$. Lastly, we put $g_{a} = \delta_{\alpha}$
where $\alpha = \log_{q}(a)$. This definition is consistent with 
that in 2.2.3.6, we just deal here with continuity and cuts. \\

\xymatrix{
& & & & & &
\omega^{(0)} \ar@<1ex>[rr]^-{\breve{P}(a)} 
\ar@(dl,ul)[]^-{\txt{$f\left(\overline{A_{s}^{(0)}}\right)$\\
$\left(A_{u}^{(0)}\right)^{\lambda}$}}
& &
\omega^{(\infty)} \ar@<1ex>[ll]^-{\breve{P}(b)^{-1}}
\ar@(ur,dr)[]^-{\txt{$f\left(\overline{A_{s}^{(\infty)}}\right)$\\
$\left(A_{u}^{(\infty)}\right)^{\lambda}$}} \\
}

It is an important fact that $\breve{P}$ is not an elliptic
function. Here are the effects of the two fundamental loops
of $\pi_{1}(\mathbf{E}_{q})$:
\begin{itemize}

\item{\textbf{Automorphy due to the monodromy of the logarithm. -}
Let $a$ turn counterclockwise once around $0$. Then 
$g_{a}(c) = \delta_{\alpha}(u q^{y}) = e^{- 2 \imath \pi \tau \alpha y}$
is multiplied by $e^{-2 \imath \pi y}$ and, since $e_{q,c}$
is uniform on $\mathbf{C}^{*}$, $\psi_{a}(c)$ is multiplied by 
$\gamma_{2}(c)$. This $\gamma_{2}$ sends $q$ to $1$, 
so that it defines an element $\gamma_{2}^{(0)}$ of 
the local Galois group at $0$ (and the like at $\infty$):
clearly, it represents the plain classical loop around $0$ 
in $\mathbf{C}^{*}$.}

\item{\textbf{Automorphy due to the defect 
of ellipticity of $\breve{P}$. -}
>From the equality $\log_{q}(qa) = \log_{q}(a) + 1$, we draw:
$\frac{g_{qa}}{g_{a}} = 
\frac{\delta_{\alpha+1}}{\delta{\alpha}} = 
\delta_{1}$,
whence 
$\frac{\psi_{qa}(c)}{\psi_{a}(c)} = 
\frac{e_{q,c}(qa)/g_{qa}(c)}{e_{q,c}(a)/g_{a}(c)} =
\frac{c}{\delta_{1}(c)} = \gamma_{1}(c)$. 
This $\gamma_{1}$ defines again an element $\gamma_{1}^{(0)}$ of 
the local Galois group at $0$, similarly at $\infty$.
We thus have:
$$
\breve{P}(qa) = 
\left(\gamma_{1}^{(\infty)}(A^{(\infty)})\right) P(a) 
\left(\gamma_{1}^{(0)}(A^{(0)})  \right)^{-1}.
$$
While $\breve{P}$ is not elliptic (except in the regular case), 
its left and right automorphy factors under the action of 
$q^{\mathbf{Z}}$ are elements of the local Galois groups.}

\end{itemize}

\textbf{3.2.2.3 A density lemma. -}
We again restrict ourselves to a category
$\mathcal{C}'_{\Sigma,C}$, and, occasionnally, 
to $\mathcal{C}'_{\Sigma,reg}$.
The arguments apply as well to the tannakian subcategory 
$<\mathcal{X}>$ generated by an object 
$\mathcal{X} = (A^{(0)},P,A^{(\infty)})$. 
We consider a $\{0,\infty\}$-subset $E$
of the Galois groupoid $\breve{G}$, with base $\{0,\infty\}$
($\simeq Spec(\mathbf{C} \times \mathbf{C})$),
with arrow sets $Aut^{\otimes}(\omega^{(0)})$,
$Aut^{\otimes}(\omega^{(\infty)})$, 
$Iso^{\otimes}(\omega^{(0)},\omega^{(\infty)})$ and
$Iso^{\otimes}(\omega^{(\infty)},\omega^{(0)})$ 
with the following constraints:
\begin{itemize}

\item{The component $E(0)$ above $0$ contains the unipotent
loop $\mathcal{X} \leadsto A_{u}^{(0)}$, and a family of
semi-simple loops 
$\mathcal{X} \leadsto \gamma_{i}(\overline{A_{s}^{(0)}})$
where the $\gamma_{i} \in Hom_{grp}(\mathbf{C^{*}},\mathbf{C^{*}})$
are such that $\underset{i}{\bigcap} \ker \gamma_{i} = q^{\mathbf{Z}}$.
Alternatively, if viewed as elements of 
$Hom_{grp}(\mathbf{C^{*}}/q^{\mathbf{Z}},\mathbf{C^{*}})$,
the $\gamma_{i}$ are such that $\underset{i}{\bigcap} \ker \gamma_{i} = \{1\}$.The component $E(\infty)$ above $\infty$ contains the
corresponding elements at $\infty$.}

\item{The component $E(0,\infty)$ above $0,\infty$ contains the
paths $\mathcal{X} \leadsto \breve{P}(a)$ for all 
$a \in \mathbf{C}^{*} - \Sigma'$.}

\end{itemize}
Then, the $\{0,\infty\}$-set $E$ generates 
a Zariski dense subgroupoid of $\breve{G}$. 


\subsubsection*{3.2.3~~~Relation to the results in 3.1}

We now relate the fibre functors and Galois groupoid studied here
with those described previously.
Let $F: A \rightarrow B$ be a morphism of flat objects, meromorphic
on $\mathcal{C}'$. Then, $F e_{q,A} = e_{q,B} S$, where $S$ has
constant coefficients and $S \overline{A} = \overline{B} S$. 
Thus, $F(z_{0}) e_{q,A}(z_{0}) = e_{q,B}(z_{0}) S$, that is,
$\mathcal{X} \leadsto e_{q,A}(z_{0})$ is a natural transformation
from $\omega^{(0)}$ (the fibre functor in chapter 3) to 
$\omega^{(0)}_{z_{0}}$. It is however not $\otimes$-compatible.
On the other side, the relation $S \overline{A} = \overline{B} S$
implies that, for any map 
$\psi: \mathbf{C}^{*} \rightarrow \mathbf{C}^{*}$
such that $\psi(c)$ depends only on $\overline{c}$, 
$\mathcal{X} \leadsto \psi(\overline{A^{(0)}})$ provides
a natural isomorphism (not a $\otimes$-isomorphism) 
from $\omega^{(0)}$ to itself. Hence:
$$
\mathcal{X} \leadsto 
\left(\psi(\overline{A^{(0)}})\right)^{-1} e_{q,A}(z_{0})
$$
is again a natural transformation from $\omega^{(0)}$ 
to $\omega^{(0)}_{z_{0}}$. For it to be $\otimes$-preserving,
it is necessary and sufficient that the map
$g: c \mapsto \frac{e_{q,c}(z_{0})}{\psi(c)}$
be a group homomorphism $\mathbf{C}^{*} \rightarrow \mathbf{C}^{*}$.
Otherwise said, $\psi(c) = \frac{e_{q,c}(z_{0})}{g(c)}$,
where $g$ is a group homomorphism 
$\mathbf{C}^{*} \rightarrow \mathbf{C}^{*}$
and the condition $\psi(q c) = \psi(c)$ says that $g(q) = z_{0}$.
We have found again the $g_{a}$ and $\psi_{a}$ of 3.2.
We conclude that any such group homomorphism $g_{z_{0}}$
provides a $\otimes$-isomorphism:
$$
\mathcal{X} \leadsto 
g_{z_{0}}(\overline{A^{(0)}}) = g_{z_{0}}(\overline{A^{(0)}_{s}})
$$
from $\omega^{(0)}$ to $\omega^{(0)}_{z_{0}}$. Of course,
composition of such isomorphisms exactly gives the morphisms 
$\omega^{(0)}_{z_{0}} \rightarrow \omega^{(\infty)}_{z_{0}}$
already found. \\

Now, as regards the connection component, we have already seen that
the connection matrix corresponding to $(A^{(0)},M,A^{(\infty)})$ is
$P = \left(e_{q,A^{(\infty)}}\right)^{-1} M e_{q,A^{(0)}}$. 
Therefore, the morphism
$\mathcal{X} \leadsto \breve{P}(z_{0})$ from $\omega^{(0)}$ to 
$\omega^{(\infty)}$ provided by the \emph{twisted}
connection matrix $\breve{P}$ is but the composition:
$$
\begin{CD}
\omega^{(0)}
@>{g_{z_{0}}(\overline{A^{(0)}})}>>
\omega^{(0)}_{z_{0}}
@>{M(z_{0})}>>
\omega^{(\infty)}_{z_{0}}
@<{g_{z_{\infty}}(\overline{A^{(\infty)}})}<<
\omega^{(\infty)}.
\end{CD}
$$
To summarize the relationship between the description of the Galois
groupoid given here, in 3.2, and the previous descriptions, given
in 2.2 and 3.1, we must introduce some more notations (which won't
be used elsewhere). We shall call $G$ the former groupoid (its base set
is the disjoint union $\mathbf{C}^{*} \amalg \mathbf{C}^{*}$) 
and $\breve{G}$ the latter one
(its base set is $\{0,\infty\}$). Recall that 
$\Gamma_{a} \in Iso^{\otimes}(\omega^{(0)}_{a},\omega^{(\infty)}_{a})$
was defined in proposition 3.1.2.2. The corresponding ``twisted''
element of $Iso^{\otimes}(\omega^{(0)},\omega^{(\infty)})$, defined
in 3.2.2.2 as $\mathcal{X} \leadsto \breve{P}(a)$, we denote by
$\breve{\Gamma}_{a}$. Similarly, for the local groupoids, we define,
for $a \in \mathbf{C}^{*}$:
$$
\Gamma_{a}^{(0)} \in 
Iso^{\otimes}(\omega^{(0)},\omega^{(0)}_{a})
\text{~by~}
\mathcal{X} \leadsto g_{a}(\overline{A^{(0)}_{s}})
\text{~and~}
\Gamma_{a}^{(\infty)} \in 
Iso^{\otimes}(\omega^{(\infty)}_{a},\omega^{(\infty)})
\text{~by~}
\mathcal{X} \leadsto g_{a}(\overline{A^{(\infty)}_{s}}).
$$ \\

\textbf{3.2.3.1 Proposition. -}
\emph{For all $a$,
$\breve{\Gamma}_{a} = \Gamma_{a}^{(\infty)} \circ 
\Gamma_{a} \circ \Gamma_{a}^{(0)}$.}
\hfill $\Box$ \\

It is moreover clear that the correspondence thus obtained between
the elements of both Galois groupoids preserve the continuity of the
underlying morphisms of groups, and that the fundamental groupoids
are mapped to each other through this correspondance. \\

The following commutative diagram exhibits the relations
linking elements of $G$ to elements of $\breve{G}$.
The leftmost and rightmost fibre functors $\omega^{(0)}$
and $\omega^{(\infty)}$ (here evaluated on an object 
$\mathcal{X} = (A^{(0)},M,A^{(\infty)})$) should be thought of
as base points of $\breve{G}$, while the inner vertical triangles
respectively belong to the subgroupoids $G^{(0)}$ and $G^{(\infty)}$ 
of $G$.

{\tiny
\xymatrix{
& &
& {\omega^{(0)}_{a}}(\mathcal{X}) \ar[rrr]^{M(a)} 
\ar'[d][ddd]^{g_{\frac{b}{a}}(\overline{A^{(0)}_{s}})}
& & & 
{\omega^{(\infty)}_{a}}(\mathcal{X}) 
\ar'[d][ddd]^{g_{\frac{b}{a}}(\overline{A^{(\infty)}_{s}})} 
\ar[drrr]^{g_{a}(\overline{A^{(\infty)}_{s}})} & & & \\
\omega^{(0)}(\mathcal{X}) 
\ar[urrr]^{g_{a}(\overline{A^{(0)}_{s}})} 
\ar[rr]^{g_{c}(\overline{A^{(0)}_{s}})} 
\ar[ddrrr]^{g_{b}(\overline{A^{(0)}_{s}})} & &
{\omega^{(0)}_{c}}(\mathcal{X}) \ar[rrr]^{M(c)} 
\ar[ur]_{g_{\frac{a}{c}}(\overline{A^{(0)}_{s}})} & & & 
{\omega^{(\infty)}_{c}}(\mathcal{X}) 
\ar[ur]^{g_{\frac{a}{c}}(\overline{A^{(\infty)}_{s}})} 
\ar[rrrr]^{g_{c}(\overline{A^{(\infty)}_{s}})} &  
& & & \omega^{(\infty)}(\mathcal{X}) \\
& &
& & & & 
& & & \\
& &
& {\omega^{(0)}_{b}}(\mathcal{X}) \ar[rrr]^{M(b)} 
\ar[uul]_{g_{\frac{c}{b}}(\overline{A^{(0)}_{s}})} 
& & & 
{\omega^{(\infty)}_{b}}(\mathcal{X}) 
\ar[uul]^{g_{\frac{c}{b}}(\overline{A^{(\infty)}_{s}})} 
\ar[uurrr]_{g_{b}(\overline{A^{(\infty)}_{s}})}
& & & \\
}
}

\bigskip

\textbf{3.2.3.2 Remark. -}
The above diagram can also be understood as explicitly identifying $G$ 
with the groupoid induced by $\breve{G}$ through the canonical
projection $\mathbf{C}^{*} \amalg \mathbf{C}^{*} \rightarrow \{0,\infty\}$,
in the sense of \cite{DF}, 1.6.


\subsection{The connection component for abelian regular equations}

In 3.2, we have given a more concrete description
of the ``connection component''. Note however that
neither description can be considered as a result
of topological nature, as it requires uncountably 
many generating paths. A partial solution will be 
proposed here.\\

For a regular equation $A$, the Galois group at $0$
is Zariski-generated by the $P(a)^{-1} P(b)$, where
$a,b$ run over $\mathbf{C}^{*} - q^{\mathbf{Z}} \mathcal{S}(A)$.
Generally speaking, it is generated by the local component at $0$,
$G^{(0)}$, \emph{one} conjugate 
$\breve{P(a)}^{-1} G^{(\infty)} \breve{P}(a)$
of the local component at $\infty$, and the connection component,
that is, the group generated by the $\breve{P}(a)^{-1} \breve{P}(b)$
where $a,b$ run over $\mathbf{C}^{*} - \Sigma'$. In the case of 
torsionless equations (see 3.3.3.1), the twisting factors of 
$\breve{P}$ belong to the local Galois groups and one can replace
the true connection component by the fake one, generated by the 
$P(a)^{-1} P(b)$. \\

We shall describe here the connection component (hence the
Galois group) for \emph{regular abelian equations}; by this,  
we mean those regular systems such that all the values of 
the connection matrix commute with each other, that is, 
such that the connection component is a commutative group.


\subsubsection*{3.3.1~~~Summary of some results in \cite{SerreGacc}}

We here apply to our context results of chapters 3 and 4 of Serre's book
\cite{SerreGacc}, on which this section heavily relies. We identify 
the complex torus $\mathbf{E}_{q} = \mathbf{C}^{*}/q^{\mathbf{Z}}$
with the corresponding elliptic curve and the latter
with the set of its complex points. We also identify the rational 
function field $k(\mathbf{E}_{q})$ of the algebraic curve 
$\mathbf{E}_{q}$ with the field $\mathcal{M}(\mathbf{E}_{q})$ of elliptic 
functions. The connection matrix $P$ defines a meromorphic function on 
$\mathbf{E}_{q}$ and we call $S$ its singular locus (made up of its poles
along with those of $P^{-1}$, see 1.2.2), a finite subset of 
$\mathbf{E}_{q}$. \\

We fix once for all a base point $a_{0} \in \mathbf{E}_{q} - S$. 
The meromorphic mapping:
$$
\begin{cases}
\mathbf{E}_{q} \rightarrow Gl_{n}(\mathbf{C}) \\
a \mapsto (P(a_{0}))^{-1} P(a)
\end{cases}
$$
can be seen as a rational map 
$f: \mathbf{E}_{q} \rightarrow Gl_{n}(\mathbf{C})$
and the holomorphy of $P$ on $\mathbf{E}_{q} - S$ implies that
$f$ is regular on the curve $\mathbf{E}_{q} - S$. Hence the Galois group 
$G \subset Gl_{n}(\mathbf{C})$ of $A$ is Zariski-generated
by the image of the regular map 
$f: \mathbf{E}_{q} - S \rightarrow Gl_{n}(\mathbf{C})$. \\

Being parameterized by a Zariski-dense subset of an irreductible
projective curve, $G$ is a connected algebraic group and we have
assumed it to be commutative. It is therefore the product of an 
algebraic torus
\footnote{To avoid any mishap, we shall systematically call
\emph{complex torus} an elliptic curve over $\mathbf{C}$ and 
\emph{algebraic torus} a torus in the sense of the theory of
linear algebraic groups.}
and an affine space:
$$
G \simeq \mathbf{G_{m}}^{k} \times \mathbf{G_{a}}^{l}
$$
(\emph{see} \cite{Borel}, 3.8 and 4.8 and \cite{Springer}, 3.4). \\

According to \cite{SerreGacc}, theorem 1 of chapter 3, a \emph{module}:
$$
\frak{M} = \sum_{p \in S} n_{p} [p] \;,\; (\text{all } n_{p} \gt 0)
$$
is associated to $f$, that is, an effective divisor on $\mathbf{E}_{q}$ 
with support exactly $S$. To this module is associated a
\emph{generalized jacobian}:
$$
\Phi_{\frak{M}}: \mathbf{E}_{q} \rightarrow J_{\frak{M}},
$$
where $J_{\frak{M}}$ is a commutative algebraic group and
$\Phi_{\frak{M}}$ is a rational map defined up to a 
translation in $J_{\frak{M}}$ ; for instance, fixing a base
point $a_{0} \in \mathbf{E}_{q} - S$ and requiring that it be mapped to
$0 \in J_{\frak{M}}$ uniquely determines $\Phi_{\frak{M}}$
(\emph{see} chapter 5 of \cite{SerreGacc}). 
We shall henceforth do so. \\

The generalized jacobian $J_{\frak{M}}$ has the following
universal property: if $\frak{M}$ is associated to $f$,
then $f$ has a unique factorization $f = F \circ J_{\frak{M}}$
with $F: J_{\frak{M}} \rightarrow G$ a (regular) morphism 
of algebraic groups. From the universality and the general
properties of algebraic groups, it follows that $Im (F) = G$. \\

As a matter of fact, the module $\frak{M}$ is not uniquely
determined by $f$: any $\frak{M}' \geq \frak{M}$ will do.
There is then a corresponding map 
$J_{\frak{M}'} \rightarrow J_{\frak{M}}$. Thus, we can
factor $f$ through the projective limit of all $J_{\frak{M}}$
with modules supported by $S$. This can be done with a fixed
base point outside $S$. Call $J_{S}$ this projective limit.
We shall not confuse it with the $J_{\frak{M}}$ corresponding
to the module $\sum_{p \in S} [p]$. \\

According to \cite{SerreGacc}, p. 99, $J_{\frak{M}}$ is an extension of
the jacobian $J$ of $\mathbf{E}_{q}$ by a linear group $L_{\frak{M}}$,
the structure of which will be made explicit below:
$$
0 \rightarrow L_{\frak{M}} \rightarrow 
J_{\frak{M}} \rightarrow J \rightarrow 0.
$$
Taking the projective limits, there is a corresponding extension:
$$
0 \rightarrow L_{S} \rightarrow J_{S} \rightarrow J \rightarrow 0.
$$

In our case, $J = \mathbf{E}_{q}$. Moreover, the map 
$F: J_{\frak{M}} \rightarrow G$
is totally determined by its restriction to $L_{\frak{M}}$: indeed,
for two maps $F,F': J_{\frak{M}} \rightarrow G$ coinciding on
$L_{\frak{M}}$, one would get $F^{-1} F': J_{\frak{M}} \rightarrow G$, 
trivial on $L_{\frak{M}}$, thus factoring through a regular map from 
the projective curve $J$ to the affine group $G$ , hence trivial. \\

To summarize, to every regular abelian object of rank $n$ with 
singularities in $S$, we associate a regular map from $L_{S}$
to $Gl_{n}(\mathbf{C})$ the image of which is its Galois group.
this correspondence is one to one and we shall hereafter make
it more explicit.


\subsubsection*{3.3.2~~~The abelianized of the regular fundamental group}

Still following \cite{SerreGacc}, we introduce, for $p \in S$ 
and $n \in \mathbf{N}^{*}$, the following groups:
$$
\begin{cases}
U_{p} = \{g \in k(\mathbf{E}_{q})^{*} \;/\; v_{p}(g) \geq 0 \} \\
U_{p}^{(n)} = \{g \in k(\mathbf{E}_{q})^{*} \;/\; v_{p}(1 - g) \geq n \} \\
V_{p}^{(n)} = U_{p}^{(1)} / U_{p}^{(n)}
\end{cases}
$$
The latter is a $(n-1)$-dimensional affine space. 
In characteristic zero, it can be parametrized using
the exponential of truncated power series, so that we
can (and shall) see it as the group:
$$
V_{p}^{(n)} \simeq 
\{\exp(- a_{1} t - \cdots - a_{n-1} \frac{t^{n-1}}{n-1}) \},
$$
with $t$ a local parameter at $p$. This parametrization
will make easier the description of morphisms to $\mathbf{G_{a}}$
in 3.3.2.2. Also note that, writing $g = (g/g(0)) \times g(0)$, 
one has:
$$
\begin{cases}
U_{p} \simeq U_{p}^{(1)} \times \mathbf{G_{m}} \\
\frac{U_{p}}{U_{p}^{(n)}} \simeq V_{p}^{(n)} \times \mathbf{G_{m}}
\end{cases}
$$
We now define, for 
$\frak{M} = \sum_{p \in S} n_{p} [p] \;,\; (\text{all } n_{p} \gt 0)$:
$$
\begin{cases}
R_{\frak{M}} = \prod_{p \in S} \frac{U_{p}}{U_{p}^{(n_{p})}}
             \simeq \mathbf{G_{m}}^{S} \times \prod_{p \in S} \mathbf{G_{a}}^{n_{p}-1} \\
\Delta = \{(x,\ldots,x) \in \mathbf{G_{m}}^{S}\} \text{ (the diagonal)} \\
L_{\frak{M}} = \frac{R_{\frak{M}}}{\Delta}
       \simeq \frac{\mathbf{G_{m}}^{S}}{\Delta} \times \prod_{p \in S} \mathbf{G_{a}}^{n_{p}-1}
\end{cases}
$$
Then, going to the projective limit, we get: 
\begin{equation}
\begin{cases}
L_{S} = L_{S,s} \times L_{S,u} \;,\; \text{where} \\
L_{S,s} = 
\frac{\mathbf{G_{m}}^{S}}{\Delta} \simeq \mathbf{G_{m}}^{|S|-1} 
\text{~and}\\
L_{S,u} \simeq \prod_{p \in S} \left(1 + t_{p} \mathbf{C}[[t_{p}]] \right)
\end{cases}
\end{equation}
The groups $L_{S,s}$ and $L_{S,u}$ are respectively the semisimple
and the unipotent factor of the Jordan decomposition of the
commutative algebraic group $L_{S}$ (\emph{see} \cite{Borel},I.4.5).
Here, we have, for each $p \in S$, selected a local parameter $t_{p}$ 
at $p$ and identified the projective limit of the $\mathbf{G_{a}}^{n_{p}-1}$
to 
$$
1 + t_{p} \mathbf{C}[[t_{p}]] =
\exp \{- a_{1} t - \cdots - a_{n-1} \frac{t^{n-1}}{n-1} - \cdots \}.
$$

To all our regular abelian objects of order $n$ with singular 
locus on $S$, we have associated injectively a regular morphism 
of algebraic groups from $L_{S}$ to $Gl_{n}(\mathbf{C})$. 
To find precisely our candidate for the \emph{abelianized regular
fundamental group with singularities in $S$}:
$$
\pi_{ab,S,reg}^{1} = \frac{\pi_{S,reg}^{1}}{[\pi_{S,reg}^{1},\pi_{S,reg}^{1}]},
$$
we have to check which morphisms: $L_{S} \rightarrow Gl_{n}(\mathbf{C})$
actually arise from abelian objects in $\mathcal{C}_{\Sigma,reg}$.
Dealing with commutative groups, we just have to find all maps to $\mathbf{G_{m}}$ and 
to $\mathbf{G_{a}}$, that is, $1$-dimensional and unipotent $2$-dimensional
objects. The following (again) comes from \cite{SerreGacc} (paragraph 18
and the description of local symbols in chapter 3). \\

\textbf{3.3.2.1 The semi-simple component. -}
Let $f: \mathbf{E}_{q} \rightarrow \mathbf{G_{m}}$, which we identify 
with an elliptic function with poles and zeroes in $S$.
The corresponding map on $L_{S,u}$ is trivial. 
On the $p$-component ($p \in S$) of $L_{S,s}$, 
it is given by $x \mapsto x^{v_{p}(f)}$, the 
triviality on the diagonal $\Delta$ being forced
by the residue formula: 
$\underset{p \in S}{\sum} v_{p}(f) = 0$. \\

Such an elliptic function is characterized, up to a factor
in $\mathbf{C}^{*}$, by its divisor $\sum v_{p}(f) [p]$. 
The latter is bound by the following conditions:
$$
\begin{cases}
\underset{p \in S}{\sum} v_{p}(f) = 0 \\
\underset{p \in S}{\sum} v_{p}(f) p = 0_{\mathbf{E}_{q}}
\end{cases}
$$
Conversely, these conditions characterize the divisors of
elliptic functions.
We want to get rid of all elements of $L_{S,s}$
that are killed by such divisors. Therefore, we put:
$$
\begin{cases}
Rel_{\mathbf{E}_{q}}(S) =
\{(n_{p})_{p \in S} \in \mathbf{Z}^{S} \;/\;
\underset{p \in S}{\sum} n_{p} p = 0_{\mathbf{E}_{q}} \}
\\
L_{S,s}' = \text{ image in } L_{S,s} \text{ of }
\{(x_{p})_{p \in S} \;/\; 
\forall (n_{p})_{p \in S} \in Rel_{\mathbf{E}_{q}}(S) \;,\;
\underset{p \in S}{\prod} x_{p}^{n_{p}} = 1 \}
\end{cases}
$$
And we can now put:
$$
\pi_{ab,S,reg,s}^{1} = \frac{L_{S,s}}{L_{S,s}'}.
$$

\textbf{3.3.2.2 The unipotent component. -}
We make it explicit by considering unipotent rank $2$ objects:
$\begin{pmatrix} 1 & f \\ 0 & 1 \end{pmatrix}$, where
$f: \mathbf{E}_{q} \rightarrow \mathbf{G_{a}}$ is rational with all poles
on $S$ (zeroes don't matter here). \\

The corresponding effect on $L_{S,s}$ is trivial. 
The effect on $L_{S,u}$ is trivial only at $p$-components
such that $p$ is a pole: $v_{p}(f) = - k , k \gt 0$.
Then, using our previous ``logarithmic'' parametrization,
it is given by: $(a_{n})_{n \geq 1} \mapsto a_{k}$.
Since we can prescribe arbitrarily the orders of the poles of $f$
just by putting zeroes elsewhere, we get the whole dual of
$L_{s,u}$ and may conclude:
$$
\pi_{ab,S,reg,u}^{1} = L_{S,u}.
$$

\textbf{3.3.2.3 Theorem. -}
\emph{The abelian regular objects with singularities in $S$ are 
classified by the representations of the following algebraic group
(see equation (3)):
$$
\pi_{ab,S,reg}^{1} = \frac{L_{S,s}}{L_{S,s}'} \times L_{S,u}.
$$
}

This group can be seen as the \emph{abelianized regular fundamental 
group with singular locus carried by $S$}. \\

\textbf{3.3.2.4 Example: dimension 1. -}
Here is an explicit computation in dimension $1$. 
We consider the equation $\sigma_{q} y = a y$, where:
$$
a(z) = 
a_{0} \prod_{i = 1}^{r} \frac{1 - u_{i}^{-1}z}{1 - v_{i}^{-1}z} =
a_{\infty} \prod_{i = 1}^{r} \frac{1 - u_{i}w}{1 - v_{i}w}
$$
One has used $w = \frac{1}{z}$ ; the above requires that 
$a_{\infty} \prod u_{i} = a_{0}\prod v_{i}$. 
Then the \emph{connection number}:
$$
p(z) = \frac{e_{q,a_{0}}(z)}{e_{q,a_{\infty}^{-1}}(w)}
\prod_{i = 1}^{r} \frac{u_{i} \Theta_{q}(z/u_{i})}{v_{i}\Theta_{q}(z/v_{i})}
$$
is elliptic. In the regular case, one has 
$a_{0} = a_{\infty} = 1$, $\prod u_{i} = \prod v_{i}$
and the connection number is:
$$
p(z) = 
\prod_{i = 1}^{r} \frac{u_{i} \Theta_{q}(z/u_{i})}{v_{i}\Theta_{q}(z/v_{i})}.
$$
The connection component is the subgroup of 
$\mathbf{C}^{*}$ generated by the values 
$\frac{p(b)}{p(a)}$, where $a,b$ run through
$\mathbf{C}^{*} - \{u_{1},\ldots,u_{r},v_{1},\ldots,v_{r}\}$.
One can of course fix $a$. This group is clearly connected,
so it has to be $\mathbf{C}^{*}$ (the general case) or
trivial. The latter occurs if $p(z)$ is constant, that
is, if the given equation is (equivalent to) the trivial
equation $\sigma_{q} f = f$. \\

\textbf{3.3.2.5 Example: dimension 2, unipotent connection component. -}
One considers the system:
$$
A(z) = \begin{pmatrix} 1 & a(z) \\ 0 & 1 \end{pmatrix},
$$
where $a(z) \in \mathbf{C}(z)$ is such that $a(0) = a(\infty) = 0$.
Then the connection matrix is:
$$
P(z) = \begin{pmatrix} 1 & p(z) \\ 0 & 1 \end{pmatrix},
$$
where $p(z) = \underset{n \in \mathbf{Z}}{\sum} a(q^{n} z)$.
The connection component is:
$$
\left\{\begin{pmatrix} 1 & \alpha \\ 0 & 1 \end{pmatrix} \;/\; 
\alpha \in G\right\},
$$
where $G$ is the subgroup of $\mathbf{C}$ generated by 
the $p(a) - p(b)$. It is generally equal to $\mathbf{C}$
and so is the connection component. The only exceptional
case is when $p$ is constant, which means that $A$ is
rationally equivalent to the trivial equation 
$\sigma_{q} X = X$. \\

\textbf{3.3.2.6 Remark: torsionless equations. -}
We consider the case of equations with torsion-free 
local components. Various conditions on the exponents 
deserve that name ; we shall require that the set of
exponents at $0$ (resp. at $\infty$) modulo $q^{\mathbf{Z}}$
be a free subset of $\mathbf{E}_{q}$. Then, for each map
$g: \mathbf{C}^{*}/q^{\mathbf{Z}} \rightarrow \mathbf{C}^{*}$,
one can find a \emph{group homomorphism}
$f: \mathbf{C}^{*}/q^{\mathbf{Z}} \rightarrow \mathbf{C}^{*}$
such that $g(A_{s}^{(0)}) = f(\overline{A_{s}^{(0)}})$, so that
the former belongs to $G^{(0)}$ ; similarly at $\infty$.
This implies that the twisting factors in each 
particular $\breve{P}(a)$ belong to the local Galois groups,
so that the values $P(a)$ of the \emph{untwisted} connection
matrix belong to the Galois groupoid, and moreover generate it
along with the local groups. So we can replace the true
connection component by the fake connection component
generated by the $P(a)^{-1}P(b)$. To the latter, the content of 
3.3.2 applies word for word. For more details, \emph{see} 
\cite{JSthese}, second part, 3.1.2.2. 



\section{Additional results}


\subsection{Confluence of galoisian automorphisms}

This paragraph extends to the Galois group (more precisely,
to the fundamental groupoid as we have defined it) the
confluence results obtained in \cite{JSAIF}. These results
are closely related to semicontinuity results obtained by
Yves Andr\'e, \emph{see} \cite{YA2}, in an algebro-geometric setting.
Our results are less general but more explicit, since we 
follow specified elements along $q$-paths in $\mathbf{C}^{*}$.
However, these results are not very complete since we do
not know the whole story about \emph{relations}. \\

According to the general assumptions in \cite{JSAIF}, chapters
3 and 4, we shall consider the matrix $A$ of a fuchsian $q$-difference
equation, depending on $q$ in such a way that:
$$
\frac{A - I_{n}}{q - 1} \to \tilde{B},
$$
where the differential equation is fuchsian and non
resonant at $0$ and $\infty$, which entails that the 
$q$-difference equation also is for $q$ close enough to $1$. 
Calling $\tilde{z}_{1},...,\tilde{z}_{r}$, we assume
the convergence to be uniform on any compact subset of 
$\mathbf{C}^{*} - \underset{1 \leq j \leq r}{\bigcup} \tilde{z}_{j}
q_{0}^{\mathbf{R}}$. Last, we assume the Jordan structures
at $0$ and $\infty$ to vary ``flatly'' (see \emph{loc. cit.}
for a precise formulation). \\

We shall then attach to $A$ the \emph{canonical} triple 
$(A(0),P,A(\infty))$ as defined in \emph{loc. cit.} 
in the non resonant case. Of course, we shall have
this triple vary along with $\epsilon$ and $q$.


\subsubsection*{4.1.1~~~General conventions}

First, we shall slightly modify the choice of $q$-characters
and $q$-logarithm, so as to get simpler determinations of
their limits as $q \to 1$. We take:
$$
\begin{cases}
e_{q,c}(z) = 
z^{\epsilon(c)} \frac{\Theta_{q}(-z)}{\Theta_{q,\overline{c}}(-z)} \\
l_{q}(z) = -z \frac{\Theta_{q}'(-z)}{\Theta_{q}(-z)}
\end{cases}
$$

Following the conventions of \cite{JSAIF}, we have $q$ tend 
to $1$ \emph{along a fixed logarithmic spiral}. We fix
$q_{0} = e^{- 2 \imath \pi \tau_{0}}$, with $Im(\tau_{0}) \gt 0$
and take $q = q_{0}^{\epsilon} = e^{- 2 \imath \pi \tau}$, 
where $\tau = \tau_{0} \epsilon$, $\epsilon \gt 0$. We shall
have $\epsilon$ tend to $0$ along $\mathbf{R}_{+}^{*}$. 
The following assertions are proven in \emph{loc. cit.}:
\begin{enumerate}

\item{For $\epsilon \to 0^{+}$, let $c_{\epsilon} \in \mathbf{C}^{*}$
be such that 
$\frac{c_{\epsilon} - 1}{q - 1} \to \gamma \in \mathbf{C}$.
Then $e_{q,c_{\epsilon}}(z) \to z^{\gamma}$.}

\item{In the same circumstances, $(q-1)l_{q}(z) \to \log z$.}

\end{enumerate}
Here, we take $log z = 2 \imath \pi x$ and 
$z^{\gamma} = e^{2 \imath \pi \gamma x}$, where we have written
$z = e^{2 \imath \pi x} , x = u + v \tau_{0} , 
u \in ]-\frac{1}{2},\frac{1}{2}[$. Said otherwise, we have
taken a cut along $- q_{0}^{\mathbf{R}}$. \\

Then, we must choose loops in the local components. We start 
with the ``fundamental level'' value $q = q_{0} , \epsilon = 1$ ;
afterwise, we shall need a calibration (or ``renormalisation'')
to handle the ``level $\epsilon$''. We split:
$$
\mathbf{C} = \mathbf{R} \oplus \mathbf{R} \tau_{0}
$$
and, writing $x = u + v \tau_{0} , u , v \in \mathbf{R}$,
we define :
$$
\begin{cases}
x \overset{p_{1}}{\mapsto} u \\
x \overset{\lambda p_{2}}{\mapsto} \lambda v
\end{cases}
$$
These are group homomorphisms $\mathbf{C} \rightarrow \mathbf{R}$
sending $\mathbf{Z}$ to $\mathbf{Z}$, thereby defining:
$$
\begin{cases}
z \overset{\gamma_{1}}{\mapsto} e^{2 \imath \pi u} \\
z \overset{\gamma_{2}^{\lambda}}{\mapsto} e^{2 \imath \pi \lambda v}
\end{cases}
$$
(we have written $z = e^{2 \imath \pi x}$). These are 
group homomorphisms: $\mathbf{C}^{*} \rightarrow \mathbf{C}^{*}$,
with images in the unit circle $\mathbf{U}$. \\

We want to relate these loops
\footnote{We shall concentrate on the component at $0$,
the case of $\infty$ being obviously the same.} 
to our previous fundamental loops
in the semi-simple local components at level $q = q_{0}^{\epsilon}$:
$$
G_{q}^{(0)} = Hom_{grp}(\mathbf{C}^{*}/q^{\mathbf{Z}},\mathbf{C}^{*}).
$$
We see, writing $x = u + v \tau_{0} = u + \frac{v}{\epsilon} \tau$
that the latter are precisely $\gamma_{1}$ and $\gamma_{2}^{\epsilon}$. \\

In the same way, we have to ``renormalize'' the twisting factors in
$\breve{P}$, mainly the $g_{a}$. This is done writing, as before,
$a = e^{2 \imath \pi \alpha}$ and taking as $g_{a}$ the group
homomorphism from $\mathbf{C}^{*}$ to $\mathbf{C}^{*}$ induced by
$$
\begin{cases}
\mathbf{C} \rightarrow \mathbf{C} \\
u + v \tau_{0} \mapsto - \frac{v}{\epsilon} \alpha
\end{cases}
$$
(the latter clearly sends $\mathbf{Z}$ to $\mathbf{Z}$). \\

We must now choose galoisian automorphisms at the level
$\epsilon = 0$,  that is, for the limit differential equation.
The latter has a pole at $0$, so that its local Galois group
is an image in $GL_{n}(\mathbf{C})$ of 
$\pi_{1}((\mathbf{C}^{*},0),.)^{alg} = \mathbf{Z}^{alg}$.
We only take care here of the semi-simple component
$\mathbf{Z}_{s}^{alg} = Hom_{grp}(\mathbf{C}^{*},\mathbf{C}^{*})$.
Unhappily, we shall not arrive at the usual fundamental loop
$1 \in \mathbf{Z}$, here identified with 
$Id_{\mathbf{C}^{*}} \in Hom_{grp}(\mathbf{C}^{*},\mathbf{C}^{*})$.
To define specific elements, we split:
$$
\mathbf{C} = \frac{1}{\tau_{0}} \mathbf{R} \oplus \mathbf{R}
$$
and, writing $x' = \frac{u'}{\tau_{0}} + v' , u' , v' \in \mathbf{R}$,
we define:
$$
\begin{cases}
x' \overset{w \tilde{p}_{1}}{\mapsto} w u' \\
x' \overset{\tilde{p}_{2}}{\mapsto} v'
\end{cases}
$$
These group homomorphisms $\mathbf{C} \rightarrow \mathbf{C}$
send $\mathbf{Z}$ to $\mathbf{Z}$ and, writing 
$z' = e^{2 \imath \pi x'}$, we can define:
$$
\begin{cases}
z' \overset{\tilde{\gamma}_{1}^{w}}{\mapsto} e^{2 \imath \pi w u'} \\
z' \overset{\tilde{\gamma}_{2}}{\mapsto} e^{2 \imath \pi v'} \\
\end{cases}
$$
These loops at $0$ define elements of the semi-simple component
of the local Galois group of the equation
$$
\delta \tilde{X} \underset{def}{=} z \frac{d}{dz} \tilde{X} = 
\tilde{B} \tilde{X}
$$
through the matrices $\gamma_{i}(e^{2 \imath \pi \tilde{B}(0)})$.
Here, the differential equation is assumed to be fuchsian at $0$, 
so that $\tilde{B}(0) \in M_{n}(\mathbf{C})$.
These matrices generate a Zariski-dense subgroup of the local
Galois group, though not the monodromy group: the latter is
generated by $e^{2 \imath \pi \tilde{B}(0)}$, which comes from
$\tilde{\gamma}_{1}^{\frac{1}{\tau_{0}}} \tilde{\gamma}_{2}$.


\subsubsection*{4.1.2~~~Confluence of the connection component}

>From \cite{JSAIF}, we know that $P$ tends to $\tilde{P}$, 
a matrix that is locally constant on the non connected open 
subset 
$$
\tilde{\Omega} = 
\mathbf{C}^{*} - \underset{0 \leq j \leq r}{\bigcup} \tilde{z}_{j}
q_{0}^{\mathbf{R}}
$$
of $\mathbf{S}$, where we have put, for simplicity, 
$\tilde{z}_{0} = 1$. Of course, 
$q_{0}^{\mathbf{R}} = e^{- 2 \imath \pi \tau_{0} \mathbf{R}}$.



\bigskip \hrule \bigskip

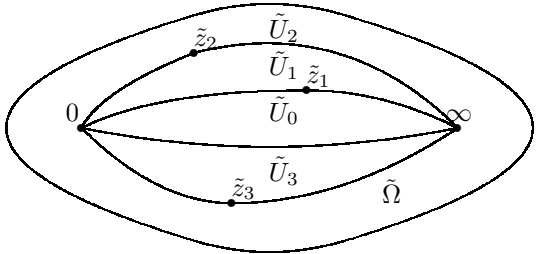
\unitlength = 1cm

\bigskip

\begin{picture} (15,3)


\qbezier (4,0)(9,1.5)(4,3)     

\qbezier (2,0)(-3,1.5)(2,3)      

\qbezier (2,3)(3,3.3)(4,3)   

\qbezier (2,0)(3,-0.3)(4,0)  


\put(0.5,1.5){\circle*{0.1}}
\put(0.3,1.6){$0$}

\put(5.5,1.5){\circle*{0.1}}
\put(5.35,1.6){$\infty$}

\put(2,2.5){\circle*{0.1}}
\put(2,2.6){$\tilde{z}_{2}$}

\put(3.5,2){\circle*{0.1}}
\put(3.5,2.1){$\tilde{z}_{1}$}

\put(2.5,0.5){\circle*{0.1}}
\put(2.5,0.6){$\tilde{z}_{3}$}


\qbezier (2,2.5) (4,3) (5.5,1.5)

\qbezier (3.5,2) (4.5,2) (5.5,1.5)

\qbezier (2.5,0.5) (4,0.5) (5.5,1.5)

\qbezier (0.5,1.5) (3,1) (5.5,1.5)

\qbezier (0.5,1.5) (0.75,2) (2,2.5)

\qbezier (0.5,1.5) (1.5,2) (3.5,2)

\qbezier (0.5,1.5) (1.5,0.5) (2.5,0.5)


\put(3,1.6){$\tilde{U}_{0}$}

\put(3,2.2){$\tilde{U}_{1}$}

\put(3,2.7){$\tilde{U}_{2}$}

\put(3,0.8){$\tilde{U}_{3}$}




\put(4.5,0.5){$\tilde{\Omega}$}


\put(-1,-0.8) {The boundary of $\tilde{\Omega}$ is made up of the $q$-spirals
generated by $1$ and the singularities $\tilde{z}_{i}$ of $\tilde{B}$}

\end{picture}

\bigskip \bigskip \bigskip \hrule \bigskip

The matrix $\tilde{P}$ takes a finite number of values
$\tilde{P}(a_{i}) , 0 \leq i \leq r$ and the 
$\left(\tilde{P}(a_{i})\right)^{-1} \tilde{P}(a_{i-1}) , 1 \leq i \leq r$
are the monodromy operators at singularities other than $0,\infty$.
But we have built our Galois isomorphisms with $\breve{P}$
instead of $P$, so we have to study the fate of 
$\psi_{a}(A_{s}^{(0)})$ and $\psi_{a}(A_{s}^{(\infty)})$ as $q \to 1$. \\

Under the confluence assumptions, the exponents $c$ of $A$ 
at $0$, resp. at $\infty$ are such that
$\frac{c - 1}{q - 1} \to \tilde{c}$, 
the exponents of $\tilde{B}$ at $0$, resp. at $\infty$.
We shall do the computation with $c = q^{\tilde{c}}$, which,
according to the lemma of \cite{JSAIF}, 3.1, does not matter. 
Following 4.1.1, we find that $e_{q,c}(a) \to e^{2 \imath \pi \tilde{c} a}$
then $g_{a}(c) \to e^{2 \imath \pi v'}$, where we have written
$\tilde{c} = \frac{u'}{\tau_{0}} + v' , u' , v' \in \mathbf{R}$.
Last, we obtain
$\psi_{a}(c) \to e^{2 \imath \pi \frac{\alpha u'}{\tau_{0}}}$:
the twisting factor tends to 
$\tilde{\gamma}_{1}^{\frac{\alpha}{\tau_{0}}}$. To be precise:
$$
\breve{P}(a) \to 
\tilde{\gamma}_{1}^{\frac{\alpha}{\tau_{0}}}
(e^{2 \imath \pi \tilde{B}(\infty)})
\tilde{P}(a)
\left(\tilde{\gamma}_{1}^{\frac{\alpha}{\tau_{0}}}
(e^{2 \imath \pi \tilde{B}(0)})\right)^{-1}.
$$
Therefore, up to factors from the local Galois groups
at $0$ and $\infty$, we get the whole system of monodromy
factors at other singularities.


\subsubsection*{4.1.3~~~Confluence of the local components}

\textbf{4.1.3.1 Unipotent part. -}
The unipotent loop at $0$ defined at the beginning of 2.2.3
gives rise to a continuous
family of Galois automorphisms $A_{u}^{(0)}$. We renormalize\
and follow instead the $\left(A_{u}^{(0)}\right)^{\frac{-1}{\tau}}$:
the limit is plainly the unipotent Galois automorphism at level $0$,
$e^{2 \imath \pi \tilde{B}_{n}(0)}$, obtained from the nilpotent
component in the additive Dunford decomposition of $\tilde{B}(0)$. \\

\textbf{4.1.3.2 Semi-simple part: generators. -}
As noticed before, the exponents $c$ of $A$ at $0$ are such that 
$\frac{c - 1}{q - 1} \to \tilde{c}$, the exponents of $\tilde{B}$ 
at $0$. Again, we compute with the innocuous assumption that
$c = q^{\tilde{c}}$ and we write
$\tilde{c} = \frac{u'}{\tau_{0}} + v' , u' , v' \in \mathbf{R}$.
Then:
$$
\begin{cases}
\gamma_{1}(c) = e^{2 \imath \pi u' \epsilon} \to 1 \\
\gamma_{2}(c) = e^{2 \imath \pi v'} \to \tilde{\gamma}_{2}(\tilde{c})
\end{cases}
$$
We obtain eventually:
$$
\begin{cases}
\gamma_{1}(A_{s}^{(0)}) \to I_{n} \\
\gamma_{2}(A_{s}^{(0)}) = \tilde{\gamma}_{2}(e^{2 \imath \pi \tilde{B}(0)})
\end{cases}
$$
The loop $\gamma_{2}$ (which we interpreted as the plain loop 
around $0$ in $\mathbf{C}^{*}$) turns infinitely fast, thus compensating 
the trivialization of the exponents (which tend to $1$).
The loop $\gamma_{1}$ (which we interpreted as the start of the move 
to infinity) turns at constant speed, so that we must accelerate it 
to compensate for the trivialization. We therefore consider 
$\gamma_{1}^{E(\frac{1}{\epsilon})}$ and find that this will do:
$$
\gamma_{1}(A_{s}^{(0)})^{\left[\frac{1}{\epsilon}\right]} 
\to \tilde{\gamma}_{1}(e^{2 \imath \pi \tilde{B}(0)}).
$$
The trip is not so smooth, involving jumps at $1/m , m \in \mathbf{N}^{*}$. 
In the end, we have reached the whole subgroup generated
by $\tilde{\gamma}_{1}$ and $\tilde{\gamma}_{2}$. \\

\textbf{4.1.3.3 Semi-simple part: relations. -}
Again, we consider exponents
that vary along $q$-spirals: $q^{\gamma_{1}},\ldots,q^{\gamma_{n}}$.
One must compare the multiplicative relations of the $q^{\gamma_{i}}$ 
with the additive relations of the $\gamma_{i}$ modulo $\mathbf{Z}$. 
We therefore introduce the module of relations:
$$
L = \{(m_{1},\ldots,m_{n}) \in \mathbf{Z}^{n} \;/\;
m_{1} \gamma_{1} + \cdots + m_{n} \gamma_{n} \in \mathbf{Z}\}.
$$
Writing $\gamma_{i} = \frac{a_{i}}{\tau_{0}} + b_{i}$,
with $a_{i},b_{i} \in \mathbf{R}$, the above condition 
is equivalent to:
$$
\begin{cases}
m_{1} a_{1} + \cdots + m_{n} a_{n} = 0 \\
m_{1} b_{1} + \cdots + m_{n} b_{n} \in \mathbf{Z}
\end{cases}
$$
Now define the ``exceptional set'': 
$$
E = \{\epsilon \gt 0 \;/\; 
\frac{1}{\epsilon} \in \mathbf{Q} a_{1} + \cdots \mathbf{Q} a_{n}\}.
$$
This is an enumerable set and it is clear that, 
for $\epsilon \not \in E$, the local $q$-difference 
Galois group at level $\epsilon$ has no more relations than 
the local differential Galois group at level $0$.


\subsubsection*{4.1.4~~~Description of the monodromy action
with a fixed base point}

Consider a fixed base point $a_{0} \in \tilde{U}_{0}$ 
such that $|a_{0}| \lt |\tilde{z}_{i}|$ for $0 \leq i \leq r$.
In each slice $\tilde{U}_{i}$, choose $a_{i}$ such that
$|a_{i}| = |a_{0}|$. Then, for $i = 1, \ldots,r$,
we can define a loop with base point $a_{0}$ in the following way:
it goes from $a_{0}$ to $a_{i-1}$ along a simple circle arc
with center $0$, counterclockwise ; it turns once counterclockwise 
around $\tilde{z}_{i}$, crossing $q_{0}^{\mathbf{R}} \tilde{z}_{i}$
exactly twice ; it comes back from $a_{i-1}$ to $a_{0}$ through
the same circle arc. Thus, we get well defined elements:
$$
\Gamma_{i} \in 
\pi_{1}(\mathbf{C}^{*} - \{\tilde{z}_{1},\ldots,\tilde{z}_{r}\};a_{0}) \;,\;
i = 1,\ldots,r .
$$
Together with $\Gamma_{0}$, the class of the simple positive
circle around $0$, they form a family of free generators of
$\pi_{1}(\mathbf{C}^{*} - \{\tilde{z}_{1},\ldots,\tilde{z}_{r}\};a_{0})$. \\

The monodromy action of $\Gamma_{0}$ on the space of solutions
of the differential equation $\delta \tilde{X} = \tilde{B} \tilde{X}$
as well as the action of the simple loop around $\infty$, 
$\left(\Gamma_{0} \Gamma_{1} \cdots \Gamma_{r}\right)^{-1}$,
are obtained by confluence of the local $q$-difference Galois groups,
as seen in 4.1.3. To be precise, only differential galoisian
automorphisms were reached this way, but with the same Zariski
closures as these fundamental loops. \\

The monodromy action of $\Gamma_{i}$ for $1 \leq i \leq r$
has matrix: 
$$
\left(\tilde{P}(a_{i})\right)^{-1} \tilde{P}(a_{i-1}) =
\underset{q \to 1}{\lim} \left(P(a_{i})\right)^{-1} P(a_{i-1}).
$$
To compare it with the galoisian automorphisms
$\left(\breve{P}(a_{i})\right)^{-1} \breve{P}(a_{i-1})$,
one just has to insert the twisting factors shown in 4.1.2.


\subsection{Extension to the $p$-adic case}

We indicate here briefly how most of the previous results
can be extended to the context of $p$-adic $q$-difference
equations. The possibility to do this rests on Tate's theory
of the uniformization of rigid elliptic curves;
it was suggested by Yves Andr\'e
\footnote{Marius van der Put told us that the results in \cite{SVdP}
could be similarly extended.}.
More details are to be found in \cite{JSthese}, \cite{JSGAL},
along with detailed references to the litterature. 


\subsubsection*{4.2.1~~~Classification}

We take as a base field the completion $\mathbf{C}_{p}$
of the algebraic closure $\overline{\mathbf{Q}_{p}}$ of
the field $\mathbf{Q}_{p}$ of $p$-adic numbers. It is
an algebraically closed complete non archimedian valued 
field. One can define, for $q \in \mathbf{C}_{p}^{*}$
such that $|q| \lt 1$, an analytic curve
$\mathbf{E}_{q} = \mathbf{C}_{p}^{*}/q^{\mathbf{Z}}$ 
whose meromorphic function field $\mathcal{M}(\mathbf{E}_{q})$
is an elliptic field (i.e. algebraic function field of
genus $1$, \emph{see} \cite{ChevalleyAlg}), so that $\mathbf{E}_{q}$
can be identified to an elliptic curve over $\mathbf{C}_{p}$.
One gets in this way exactly those elliptic curves whose
modular invariant $j(\mathbf{E}_{q})$ is not an integer.
The uniformization of such an elliptic curve is obtained
with the help of the $p$-adic theta function:
$$
\Theta(z) = \prod_{n \geq 0} (1 - q^{n}z) \prod_{n \geq 1} (1 - q^{n}z^{-1}).
$$
This has all the properties we used to define our fundamental
solutions of constant coefficient systems. Hence our abelian 
$\mathbf{C}_{p}$-linear rigid tensor categories $\mathcal{E}_{f}$,
$\mathcal{S}$ and $\mathcal{C}$ can be defined, as well as the
exact $\mathbf{C}_{p}$-linear $\otimes$-functors $SE$ and $SC$,
our equivalence theorems  remain valid here, as well as the choice 
for the fibre functors.


\subsubsection*{4.2.2~~~The connection component}

To build as in 3.2 the matrix $\breve{P}$, we needed a morphism 
$g_{a}: \mathbf{C}_{p}^{*} \rightarrow \mathbf{C}_{p}^{*}$ sending
$q$ to $a$. Here, for lack of an exponential, we shall resort to
a more Zornian construction. \\

\textbf{4.2.2.1 Lemma. -}
\emph{Let $K$ be an algebraically closed field of characteristic $0$
and let $x \in K^{*} - \mu_{\infty}(K^{*})$, where $\mu_{\infty}(K^{*})$
is the torsion subgroup (roots of unity) of $K^{*}$. Then $K^{*}$
and $K^{*}/x^{\mathbf{Z}}$ are respectively isomorphic to
$(\mathbf{Q}/\mathbf{Z}) \times \mathbf{Q} \times V$ and to 
$(\mathbf{Q}/\mathbf{Z}) \times (\mathbf{Q}/\mathbf{Z}) \times V$, 
where $V$ is a  $\mathbf{Q}$-vector space and where $x \in K^{*}$
corresponds to the element $(\overline{0},1,0)$
of $(\mathbf{Q}/\mathbf{Z}) \times \mathbf{Q} \times V$.} 
\hfill $\Box$ \\

This lemma guarantees the existence of the group homomorphisms
$g_{a}$. Thus, our construction of $\breve{P}$ and our density
lemma in 3.2 remain valid. 


\subsubsection*{4.2.3~~~The local components}

We now make more precise our choice of the $g_{a}$.
let $a \in \mathbf{C}_{p}^{*}$ correspond to
$(\overline{\alpha},\beta,\xi) \in 
(\mathbf{Q}/\mathbf{Z}) \times \mathbf{Q} \times V$
and choose a lifting (a logarithm !) $\alpha$ of 
$\overline{\alpha}$ in $\mathbf{Q}$. Then the morphism:
$$
\begin{cases}
(\mathbf{Q}/\mathbf{Z}) \times \mathbf{Q} \times V \rightarrow
(\mathbf{Q}/\mathbf{Z}) \times \mathbf{Q} \times V \\
(\overline{\alpha'},\beta',\xi') \mapsto 
(\overline{\beta' \alpha},\beta' \beta,\beta' \xi) 
\end{cases}
$$
is well defined and corresponds to a group homomorphism
$\mathbf{C}_{p}^{*} \rightarrow \mathbf{C}_{p}^{*}$
sending $q$ to $a$. We now obtain naturally our fundamental
semi-simple loops $\gamma_{1}$ and $\gamma_{2}$:

\begin{enumerate}

\item{If we change the ``logarithm'' $\alpha$ to $\alpha + 1$,
$g_{a}$ is changed to $g'_{a}$ in such a way that 
$\frac{g'_{a}}{g_{a}}$ corresponds to:
$$
(\overline{\alpha'},\beta',\xi') \mapsto 
(\overline{\beta'},0,0).
$$
This morphism we take as
$\gamma_{1}: \mathbf{C}_{p}^{*} \rightarrow \mathbf{C}_{p}^{*}$.
It does send $q$ to $1$.}

\item{If we compute $\frac{\psi_{q a}}{\psi_{a}}$, we get 
$c \mapsto \frac{c}{g_{qa}(c)/g_{a}(c)}$, corresponding to:
$$
(\overline{\alpha'},\beta',\xi') \mapsto 
(\overline{\alpha'},0,\xi').
$$
This morphism we take as
$\gamma_{2}: \mathbf{C}_{p}^{*} \rightarrow \mathbf{C}_{p}^{*}$.
It also sends $q$ to $1$.}

\end{enumerate}

It is obvious from the description with $(\overline{\alpha'},\beta',\xi')$
that $\Ker \gamma_{1} \cap \Ker \gamma_{2} = q^{\mathbf{Z}}$. 
Thus, our description of the local monodromy groups in 3.2
is still valid here \emph{almost naturally}, that is,
up to the choice of a logarithm.


\newpage

\section*{Appendices}

\appendix


\section{The proalgebraic hull of $\mathbf{Z}$}

For general definitions and results on tannakian categories 
and Tannaka's duality, we used \cite{DM} and \cite{DF},
mainly the former. A much more detailed version of
what follows can be found in \cite{JSthese}.

\subsection*{A.1~~~Proof of the assertions in 2.1.3}

We consider the category 
$\mathcal{R} = Rep_{\mathbf{C}}(\mathbf{Z})$ of
finite dimensional complex representations of $\mathbf{Z}$.
Its objects can be seen as  pairs $(X,\xi)$ of a finite dimensional
complex vector space $X$ and a linear automorphism $\xi$
of $X$ (the image of $1$ by the representation).
$\mathcal{R}$ is a $\mathbf{C}$-linear abelian category, 
and it is endowed with an obvious tensor structure, 
making it a tannakian category. Duals and internal Homs
also are the natural ones. \\

We note that, as an abelian category, $\mathcal{R}$
is exactly the same as the category $Mod^{f}(R)$ of
finite length modules over the principal ideal domain
$R = \mathbf{C}[T,T^{-1}]$. This, of course, does not
apply to the tensor structure. \\

The forgetful functor $\omega$ is a fibre functor
over $\mathbf{C}$, thus making $\mathcal{R}$ a neutral
tannakian category. It is therefore equivalent to the
category of representations of the proalgebraic group
$Gal(\mathcal{R}) = \underline{Aut}^{\otimes}(\omega)$.
This is the proalgebraic hull of $\mathbf{Z}$ and we
call it $\mathbf{Z}^{alg}$. In this appendix, we shall
compute $\mathbf{Z}^{alg}$. More precisely, we shall
find the group $Aut^{\otimes}(\omega)$ of its complex 
points and then indicate briefly its scheme structure. 

\subsubsection*{A.1.1~~~Natural transformations of $\omega$}

We consider a functorial morphism: $\omega \rightarrow \omega$,
that is an element of $End(\omega)$. For the moment,
we do not assume it to be $\otimes$-compatible. 
It is defined  by giving for each representation $(X,\xi)$ 
a linear automorphism $\overline{\xi}$ of the underlying 
vector space $X$ in such a way that, for any representation
morphism $u : (X,\xi) \rightarrow (Y,\eta)$, that is, for
any linear map $u : X \rightarrow Y$ such that $u \xi = \eta u$,
one has the naturality condition: 
$u \overline{\xi} = \overline{\eta} u$.
Now, applying the naturality condition to the automorphism
$\xi$ of $(X,\xi)$, one finds that $\overline{\xi}$ commutes
with $\xi$: in the $R$-module interpretation, this means
that $\overline{\xi}$ is a morphism of $R$-modules. \\

Applying the naturality condition to the projections
of a product to its components, one sees that 
$\xi \mapsto \overline{\xi}$ is compatible with
the decomposition into blocks (direct sums).
Using Jordan decomposition, each $(X,\xi)$ 
can be decomposed into blocks having a single
eigenvalue. Moreover, there is no non trivial
morphism between such blocks admitting distinct
eigenvalues. We can therefore restrict the study 
to the subcategory of objects $(X,\xi)$ with 
the single eigenvalue $\lambda \in \mathbf{C}^{*}$.
These correspond exactly to the $R$-modules with 
$\pi_{\lambda}$-torsion, where $\pi_{\lambda}$ is 
the prime ideal generated by $T - \lambda$. 
Jordan blocks correspond to cyclic modules (torsion modules 
with one generator). Therefore, we can consider only 
morphisms between cyclic modules. \\

Each eigenvector $x$ of $\xi$ may be considered as
a morphism: $(\mathbf{C},\lambda) \rightarrow (X,\xi)$.
Applying the naturality condition shows that
$\overline{\xi}(x) = \overline{\lambda} x$. This entails
$Sp(\overline{\xi}) = \{\overline{\lambda}\}$. \\

We then consider $(X,\xi)$ and $(Y,\eta)$, cyclic modules
of ranks $n$ and $p$, therefore, respectively isomorphic
to $\frac{R}{\pi_{\lambda}^{n}}$ and $\frac{R}{\pi_{\lambda}^{p}}$. 
Thus, any morphism $u$ between them is a homothety $a \mapsto ay$
and the naturality condition is equivalent to 
$\overline{\eta}(y) = \overline{\xi}(1) y$ (the multiplication
makes sense because of the torsion condition on $y$). 
>From this, we draw first that the morphism $\overline{\xi}$
attached to the representation $(X,\xi)$ corresponding
to the $R$-module $\frac{R}{\pi_{\lambda}^{n}}$ can be
identified with an \emph{element} 
$\overline{\xi_{\lambda,n}} \in \frac{R}{\pi_{\lambda}^{n}}$ ;
second, that the image of $\overline{\xi_{\lambda,n}}$ in 
$\frac{R}{\pi_{\lambda}^{n-1}}$ is $\overline{\xi_{\lambda,n-1}}$.
It is not hard, making explicit all morphisms between cyclic modules,
to see that there are no more conditions to be checked. \\

We conclude that the natural transformation 
$(X,\xi) \mapsto \overline{\xi}$ of $\omega$ 
is completely described by a compatible family 
$(\overline{\xi_{\lambda}})_{\lambda \in \mathbf{C}^{*}}$ 
where each 
$\overline{\xi_{\lambda}} \in 
\underset{\leftarrow}{\lim} \frac{R}{\pi_{\lambda}^{n}}$.
This means that $End(\omega)$ is canonically isomorphic
to $\hat{R} = \underset{\lambda \in \mathbf{C}^{*}}{\prod}
\mathbf{C}[[T - \lambda]]$.
Concretely, a representation $(X,\xi)$ gives rise to a finite
length $R$-module, hence to a $\hat{R}$-module (since $\hat{R}$
is the projective limit of all quotient rings $\frac{R}{I}$).
Then $\overline{\xi}$ is described by a homothety on each
of its torsion components, by an element of the corresponding
$\mathbf{C}[[T - \lambda]]$.

\subsubsection*{A.1.2~~~Transformations preserving the tensor structure}

We now consider a $\otimes$-compatible transformation
$(X,\xi) \mapsto \overline{\xi}$ of $\omega$. It is
automatically an isomorphism, thus, an element of
$Aut^{\otimes}(\omega)$. Clearly, it is enough to
exploit this multiplicative condition for Jordan blocks. \\

The multiplicative condition applied in dimension $1$
shows that $\lambda \mapsto \overline{\lambda}$ is a
group homomorphism. Then, tensoring $(\mathbf{C},\lambda)$
with the unipotent Jordan block of rank $n$ gives
the the Jordan block of rank $n$ of eigenvalue $\lambda$.
It therefore suffices to exploit the multiplicative condition
on unipotent Jordan blocks. \\

So call $(\mathbf{C}^{n},\xi_{1,n})$ the standard model
for the unipotent Jordan block of rank $n$. ``A little plethysm''
gives:
$$
(\mathbf{C}^{n},\xi_{1,n}) \otimes (\mathbf{C}^{p},\xi_{1,p}) \simeq
\bigoplus (\mathbf{C}^{q},\xi_{1,q})
$$
where $q$ runs through $\{n+p-1,n+p-3,...,|n-p|+3,|n-p|+1\}$.
>From this, one concludes that the element $P$ of $\mathbf{C}[[T - 1]]$
that corresponds to our transformation satisfies:
$P(U + V + UV) = P(U) P(V)$. It is an easy exercice to show that
$P(U) = (1 + U)^{a}$ for some complex number $a$, so that, for
unipotent representations $(X,\xi)$, one has 
$\overline{\xi} = \xi^{a}$ (which makes sense for a unipotent
automorphism). This achieves the proof of the description
of $\mathbf{Z}^{alg}$ given in section 2.1.3. We propose to call
the resulting group the ''super-proalgebraic hull''.

\subsubsection*{A.1.3~~~Group scheme structure}

The group scheme $Gal(\mathcal{R}) = \underline{Aut}^{\otimes}(\omega)$
represents the functor in $\mathbf{C}$-algebras:
$S \leadsto Aut^{\otimes}(\omega \circ \Phi_{S})$ where $\Phi_{S}$
is the functor $M \leadsto M \otimes_{\mathbf{C}} S$. Call $G$
this group scheme. An element of $G(S)$ is the datum, for each
representation $(X,\xi)$, of a $S$-automorphism $\overline{\xi}$
of $X \otimes S$ with naturality and multiplicative conditions.
To make it explicit, we must again consider Jordan blocks, and
separate the one dimensional effect from the effect on unipotent
blocks. For unipotent objects, one finds that $\overline{\xi} = \xi^{a}$
for some $a \in S$. This gives the unipotent part of $G$: it
is defined by $G_{u}(S) = S$, hence it is the affine group $G_{a}$. \\

In dimension $1$, one finds a group homomorphism 
$\lambda \mapsto \overline{\lambda}$ from $\mathbf{C}^{*}$
to $S^{*}$. This gives the semi-simple part:
$G_{s}(S) = Hom_{gr}(\mathbf{C}^{*},S^{*})$. It is the
projective limit of the functors in groups 
$S \leadsto  Hom_{gr}(H,S^{*})$, where $H$ runs through
the inductive system of finitely generated subgroups
of $\mathbf{C}^{*}$. For such a subgroup, it is easy to
see that the functor is represented by the affine scheme
with algebra $\mathbf{C}[H]$, the group algebra of $H$
over $\mathbf{C}$, which has a natural Hopf algebra structure.
We conclude that 
$G_{s} = Spec \left(\mathbf{C}[\mathbf{C}^{*}]\right)$.

\subsection*{A.2~~~Additional informations}

\subsubsection*{A.2.1~~~Chevalley duality}

In \cite{ChevalleyLie}, Chevalley constructs the complex
algebraic group attached to a compact real Lie group from
the category of its representations. This can be considered 
as the true birth of Tannaka duality. The construction
is simplified a little by the fact that the representations
of a compact group are completely reducible, so that the
naturality conditions we used above boil down to compatibility
with direct sums and conjugations. \\

Therefore, one may wonder what that method would give 
for $\mathbf{Z}$: a bigger group that $\mathbf{Z}^{alg}$,
for sure, but how much bigger ? It is proved in \cite{JSthese}
(along with an explicit computation) that the group thus obtained 
has $\mathbf{Z}^{alg}$ as a subgroup of order $2$. This is due
to the ``plethysm'' formula above and its consequence that
the knowledge of $\overline{\xi_{1,2}}$ only completely
determines an element ; and the unique eigenvalue of
$\overline{\xi_{1,2}}$ has to be $\pm 1$ (again because of
plethysm).

\subsubsection*{A.2.2~~~The proalgebraic hull of $\mathbf{Z}$ over other base fields}

Our computation of $\mathbf{Z}^{alg}$ is easily seen to work the 
same if one replaces $\mathbf{C}$ by any algebraically closed
field of null characteristic. \\

We now consider a field $K$ of such that $char K = p \gt 0$. 
One can use again Jordan reduction and obtain the semi-simple
component in a similar way: it is $Hom_{gr}(K^{*},K^{*})$,
and the group scheme structure can be determined as above. \\

As for the unipotent component, it is determined by 
the ``formal group'' we called $P$, that is, a formal 
power series such that $P(U + V + UV) = P(U) P(V)$ 
and $P(0) = 1$. One can show that $P(U) = (1+U)^{a}$
where $a$ is a $p$-adic number: 
$a = \sum k_{i} p^{i} \in \hat{\mathbf{Z}}_{p}$.
Then the product 
$(1+U)^{\sum_{i \geq 0} k_{i} p^{i}} = \prod_{i \geq 0} (1+U)^{k_{i}p^{i}}$
is $U$-adically convergent. Thus, our unipotent component
is now the profinite discrete group $\hat{\mathbf{Z}}_{p}$. \\

It is easy to see that the unipotent components of the proalgebraic
hulls of $\mathbf{Z}$ over algebraically closed fields of all
characteristics actually are the fibres over the corresponding
geometric points of a single scheme over $\mathbf{Z}$. One
can naturally ask what the latter classifies. Such a modular 
interpretation has been given by Bertrand Toen in \cite{Toen}.

\subsubsection*{A.2.3~~~Bibliographic and folkloric informations}

While the above basic ``exercice'' does not appear in the 
litterature, one can find some related information in the 
following references:

\begin{enumerate}

\item{
The book in preparation of Jean-Pierre Ramis \cite{RamisGalDiff}
on the inverse problem in differential Galois theory contains
a computation of $\mathbf{Z}^{alg}$ by a more conceptual method
(but with more prerequisites about algebraic groups).}

\item{In section 5.2. of \cite{SVdP} one finds the detailed  
description of the algebraic quotients of $\mathbf{Z}^{alg}$
but neither proof nor description of the action.}

\item{According to Y. Andr\'e:\\ 
\emph{En fait  $Z^{alg}$
fait partie du folklore, Serre l'avait decrit dans son cours au College sur
les motifs, et il attribuait cette description a Chevalley.}}

\item{The paper ``Differential equations in characteristic $p$'' 
(M. van der Put,Compositio Math. 97, 227-251, 1995) section 2.3,
contains a result similar in spirit but nevertheless different          
what is computed there is a kind of additive version of $\mathbf{Z}^{alg}$,
with ``Lie-like'' elements instead of ``group-like'' elements.}

\end{enumerate}

Besides, the latter authors have proposed swifter ways to get
the description of $\mathbf{Z}^{alg}$:

\begin{enumerate}

\item{According to Y. Andr\'e:\\
\emph{Le point
essentiel est le suivant: soit  $u \in GL_n(K)$, $K$ de caracterisitique
nulle. Soit $A_u$ le plus petit sous-groupe algebrique de $GL_n$ contenant
$u$. Alors si $u$ est semi-simple,  $A_u$ est un groupe de type
multiplicatif, de torsion cyclique (son groupe des caracteres est le
sous-groupe de $\bar K ^*$ engendre par les valeurs propres de $u$).
Si $u$ n'est pas semi-simple, on considere sa "decomposition de Dunford"
dans $GL_n$, et la partie unipotente donne lieu a un facteur additif $G_a$
dans $A_u$.   $Z^{alg}$ est fabrique par limite projective filtrante, et se
decompose en $(G_a) \times (\mu_{\infty}) \times$ un protore.
}
}

\item{According to M. van der Put:\\
\emph{
For an algebraically closed field of characteristic 0, one considers $H:=k^*$
as abstract group and forms its group ring ${\bf Z}[H]$. Then 
$G:=Spec({\bf Z}[H])\times {\bf G}_a=Spec({\bf Z}[H][t])$ is an affine group 
scheme over $\bf Z$. One has to show that there are ``canonical'' isomorphisms
$G(R)\rightarrow {\rm Aut}^{\otimes}(\omega )(R)$ for any $k$-algebra $R$
(i.e., an isomorphism of the functors $G$ and ${\rm Aut}^{\otimes}(\omega )$
on the category of all $k$-algebras).
In the case that $k$ is an algebraically closed field of characteristic $p>0$,
the group scheme $G$ is a little different, namely
$G=Spec({\bf Z}[H])\times \hat{\bf Z}_p$, where $H$ denotes $k^*$ as abstract
group and $\hat{\bf Z}_p$ denotes the constant group scheme over $\bf Z$,
corresponding to the abstract group $\hat{\bf Z}_p$. Again one has to show
that the two functors $G$ and ${\rm Aut}^{\otimes}(\omega )$ are isomorphic
as functors on the category of all $k$-algebras.}
}

\end{enumerate}



\section{Automorphisms of the "field of solutions"}

We compute here a group that bears similarities at the same time
with our Galois group and with the group computed by van der Put
and Singer in \cite{SVdP}, chap. 12, pp. 150-151 (we detail in
\cite{JSthese} the important differences of these two groups). 
This is the closest we get from Picard-Vessiot theory ! 
Moreover, some aspects of the following computation 
may shed some light on the factor $g_{a}$ that helped building 
the ``mistigri'' $\psi_{a}$ in 3.2. 

\subsection*{B.1~~~Fields of solutions}

We improperly call, in the title of this appendix, ``field of solutions''
the extension:
$$
\mathbf{K} = \mathbf{k}(l_{q},(e_{q,c})_{c \in \mathbf{C}^{*}})
$$ 
of $\mathbf{k} = \mathbf{C}(z)$ obtained by adding 
the basic construction tools needed to build solutions. 
Note however that we rather seeked the said solutions in the field
$\mathbf{K}_{0} = \mathbf{K} \otimes_{\mathbf{k}} \mathbf{k}_{0}$
where $\mathbf{k}_{0} = \mathcal{M}(\mathbf{C})$. 
Since the characters are enough to generate all elliptic functions,
it is plain that actually:
$$
\mathbf{K} = \mathcal{M}(\mathbf{E})(l_{q},(e_{q,c})_{c \in \mathbf{C}^{*}}).
$$ 
It is a subfield of $\mathcal{M}(\mathbf{C}^{*})$, invariant 
under the automorphism $\sigma_{q}$. We shall compute the group:
$$
G = \{\phi \in Aut(\mathbf{K}/\mathbf{k}) \;/\; 
\phi \circ \sigma_{q} = \sigma_{q} \circ \phi \}.
$$

Here are three basic facts:

\begin{enumerate}

\item{An element $\phi \in G$ leaves invariant
$\mathbf{K}^{\sigma_{q}} = \mathcal{M}(\mathbf{E})$.}

\item{It sends a non zero $c$-character $f$ to a non zero $c$-character,
so that then, $\frac{\phi(f)}{f} \in \mathbf{K}^{\sigma_{q}}$.}

\item{Similarly, it sends a $q$-logarithm $g$ to a $q$-logarithm,
so that then $\phi(g) - g \in \mathbf{K}^{\sigma_{q}}$.}

\end{enumerate}

\subsection*{B.2~~~Automorphisms}

\subsubsection*{B.2.1~~~$G$ as a fibered product}

>From basic fact 1 above, one draws 
that there is a restriction morphism
$G \rightarrow G_{0} = Aut(\mathcal{M}(\mathbf{E})/\mathbf{C})$ ;
from basic facts 2 and 3, that the fields
$\mathcal{M}(\mathbf{E}) ((e_{q,c})_{c \in \mathbf{C}^{*}})$ and
$\mathcal{M}(\mathbf{E}) (l_{q})$ are invariant under $G$. Define:
$$
\begin{cases}
G_{1} = \{\phi \in Aut(\mathbf{K_{1}}/\mathbf{k_{1}}) \;/\; 
\phi \circ \sigma_{q} = \sigma_{q} \circ \phi \} \\
G_{2} = \{\phi \in Aut(\mathbf{K_{2}}/\mathbf{k_{2}}) \;/\; 
\phi \circ \sigma_{q} = \sigma_{q} \circ \phi \}
\end{cases}
\quad
\text{where}
\quad
\begin{cases}
\mathbf{K_{1}} = \mathcal{M}(\mathbf{E}) ((e_{q,c})_{c \in \mathbf{C}^{*}})\\
\mathbf{K_{2}} = \mathcal{M}(\mathbf{E}) (l_{q})
\end{cases}
$$
and where $\mathbf{k_{i}} = \mathbf{K_{i}} \cap \mathbf{k}$ ($i = 1,2$).
This gives rise to two morphisms: $G \rightarrow G_{1}$ and 
$G \rightarrow G_{2}$. From the results in \cite{JSAIF}, appendix A, one
draws:
$$
G = G_{1} \times_{G_{0}} G_{2}.
$$
We have to study the two projection morphisms. \\

>From now on, we shall write $\overline{\phi}$ the restriction 
of an automorphism $\phi \in G$ to $\mathcal{M}(\mathbf{E})$.

\subsubsection*{B.2.2~~~Action on the $q$-logarithm}

Writing, for $\phi \in G_{2}$, $a_{\phi} = \phi(l_{q}) - l_{q}$,
one notes that $\phi$ is determined by the pair
$(\overline{\phi},a_{\phi}) \in G_{0} \times \mathcal{M}(\mathbf{E})$.
Since $\psi(\phi(l_{q})) = l_{q} + a_{\psi} + \overline{\psi}(a_{\phi})$,
one draws that $G_{2}$ is isomorphic to the semi-direct product
of $G_{0}$ by $\mathcal{M}(\mathbf{E})$. It is not a direct product. 
We remark that the component $\mathcal{M}(\mathbf{E})$ was invisible
in our Galois group.

\subsubsection*{B.2.3~~~Action on $q$-characters}

In a similar way, $\phi \in G_{1}$ is determined 
by $\overline{\phi}$ and the family of the 
$f_{c} = \frac{\phi(e_{q,c})}{e_{q,c}} \in \mathcal{M}(\mathbf{E})^{*}$.
One must have $f_{qc} = f_{c}$ since
$\frac{e_{q,qc}}{e_{q,c}} \in \mathbf{k}$ is fixed by $\phi$.
We shall compute the image and the kernel of $G_{1} \rightarrow G_{0}$. \\

\textbf{B.2.3.1 The kernel of $G_{1} \rightarrow G_{0}$.-}
An element of the kernel has $f_{c} f_{d} = f_{cd}$, thereby defining
a group morphism $\mathbf{C}^{*} \rightarrow \mathcal{M}(\mathbf{E})^{*}$;
composing, we get 
$\mathbf{C}^{*} \rightarrow \mathcal{M}(\mathbf{E})^{*}/\mathbf{C}^{*}$,
a morphism from a divisible group to a free one, hence trivial. 
Since $f_{qc} = f_{c}$, we have obtained a morphism 
$\mathbf{C}^{*} \rightarrow \mathbf{C}^{*}$ such that $q \mapsto 1$.
This obviously works the other way round, and we conclude:
$$
Ker(G_{1} \rightarrow G_{0}) = 
Hom_{gr}(\mathbf{C}^{*}/q^{\mathbf{Z}},\mathbf{C}^{*}).
$$

\textbf{B.2.3.2 The structure of $G_{0}$.-}
>From the classical theory of elliptic curves 
(\emph{see} \cite{Silverman}), one draws:
$$
G_{0} = Aut(\mathbf{E}).
$$
This group consists in \emph{changes of variables}: $z \mapsto az$,
$z \mapsto \frac{1}{z}$ and maybe , for some special values
of $q$, $z \mapsto z^{\imath}$ and $z \mapsto z^{\rho}$. \\

\textbf{B.2.3.3 The image of $G_{1} \rightarrow G_{0}$.-}
Studying divisors, one can prove that none of the last 
three types of change of variables can come from an element
of $G_{1}$. On the other hand, $z \mapsto az$ is in the image
and we now determine its antecedents. We specify such an
antecedent by writing $\phi(e_{q,c}) = \frac{e_{q,c}(az)}{g(c)(z)}$. 
Then we find that $g$ is a morphism 
$\mathbf{C}^{*} \rightarrow \mathbf{C}^{*}$ such that $q \mapsto a$.
We conclude: \\

\emph{The image $G'_{0}$ of $G_{1} \rightarrow G_{0}$ is made up of the
translations $\overline{\phi_{a}}$, hence isomorphic to $\mathbf{E}$.
The antecedents of $\overline{\phi_{a}}$ correspond bijectively to 
the morphisms $\mathbf{C}^{*} \rightarrow \mathbf{C}^{*}$ such that 
$q \mapsto a$.} \\

\textbf{B.2.3.4 The structure of $G_{1}$. -}
To summarize, any $\phi \in G_{1}$ can be written:
$$
\phi_{g}: e_{q,c} \mapsto \frac{e_{q,c}(az)}{g(c)}
$$
where $g$ is a morphism $\mathbf{C}^{*} \rightarrow \mathbf{C}^{*}$
such that $q \mapsto a$. Moreover, one finds that:
$\phi_{g} = \phi_{h} \Leftrightarrow
\exists r \in \mathbf{Z} \;:\:
\forall c \in \mathbf{C}^{*} \;,\; g(c) = c^{r} h(c)$.
The group $G_{1}$ can therefore be identified to the quotient of
$Hom_{gr}(\mathbf{C}^{*},\mathbf{C}^{*})$ by the subgroup of the
$z \mapsto z^{r}$ ($r \in \mathbf{Z}$). Identifying
$G'_{0}$ to $\mathbf{E}$, the restriction morphism
$G_{1} \rightarrow G'_{0} \subset G_{0}$ is given by
$g \mapsto g(q) \pmod{q^{\mathbf{Z}}}$. \\

The kernel can be identified to
$Hom_{gr}(\mathbf{C}^{*}/q^{\mathbf{Z}},\mathbf{C}^{*}/q^{\mathbf{Z}})$
quotiented by the subgroup of the $z \mapsto z^{r}$ ($r \in \mathbf{Z}$).
This is isomorphic to
$Hom_{gr}(\mathbf{C}^{*}/q^{\mathbf{Z}},\mathbf{C}^{*}) =
\check{\mathbf{E}}$.



\section {Tamely irregular equations}

We study here irregular systems that admit convergent
fundamental solutions. We shall moreover assume that
no ramification is needed. We intend to return later
to the ramified case, which involves some algebraic
complications. Here, we shall expound the bare minimum
needed to handle some useful applications.


\subsection*{C.1~~~Convergent solutions}

Say that a system has \emph{level $\nu \in \mathbf{N}^{*}$} 
if all the slopes of its Newton polygons at $0$ and $\infty$ 
belong to $\frac{1}{\nu} \mathbf{Z}$. 


\subsubsection*{C.1.1~~~The categories $\mathcal{E}_{1}$ and $\mathcal{E}_{t,1}$}

>From the general properties of Newton polygons 
alluded in 1.2.1, it follows that the full
subcategory of $\mathcal{E}$ whose objects have 
level $1$ is a tannakian subcategory ; we call
it $\mathcal{E}_{1}$. \\

We know each object of $\mathcal{E}_{1}$ to be
formally equivalent at $0$ to a direct sum of
pure objects, and the same at $\infty$. We
call such an object \emph{tamely irregular}
(of level $1$) if these equivalences can be
realized through meromorphic gauge transformations.
This again defines a tannakian subcategory 
$\mathcal{E}_{t,1}$ of $\mathcal{E}_{1}$. 
We shall see it to be neutral.


\subsubsection*{C.1.2~~~Local fundamental solutions and connection matrix}

An equation of $\mathcal{E}_{t,1}$ has local solutions
$X^{(0)}$ and $X^{(\infty)}$ that are direct sums of
``pure'' solutions, of the form $M \Theta^{\mu} e_{q,A}$,
where $\mu$ runs through the slopes, $M$ is meromorphic
and $A$ is constant invertible. Then the connection matrix
$P = \left(X^{(\infty)}\right]^{-1}X^{(0)}$ is elliptic. \\

The same construction for a level $\nu$ equation would
also have given rise to an elliptic matrix, but with
respect to another elliptic curve, isogenous to 
$\mathbf{E}$.


\subsection*{C.2~~~Tame connection data}

To the tamely irregular equation $A$ of level $1$,
we associate a \emph{connection triple}:
$$
\left(\sum_{\mu \in \mathbf{Z}} z^{\mu} A^{(0)}_{\mu},P,
\sum_{\mu \in \mathbf{Z}} z^{\mu} A^{(\infty)}_{\mu}\right).
$$


\subsubsection*{C.2.1~~~The category $\mathcal{C}_{t,1}$}

As before, we thus form an abelian tensor category, that
we call $\mathcal{C}_{t,1}$. We do that just by adding the rule 
$z^{\mu} A \otimes z^{\mu'} A' = z^{\mu + \mu'} A \otimes A'$.
The same arguments as before show $\mathcal{C}_{t,1}$ to be
equivalent to $\mathcal{E}_{t,1}$: indeed, we may again
use properties of fundamental solutions with coefficients
in the field $\mathcal{M}(\mathbf{C}^{*})$.


\subsubsection*{C.2.2~~~The groups $G^{(0)}_{t,1}$ and $G^{(\infty)}_{t,1}$}

With respect to the fuchsian case, the main difference is 
that the local components are now $\mathbf{Z}$-graded. 
This means that the local Galois groups inherit a factor
$\mathbf{C}^{*}$:
$$
\begin{cases}
G_{t,1}^{(0)} = G^{(0)} \times \mathbf{C}^{*} \\
G_{t,1}^{(\infty)} = G^{(\infty)} \times \mathbf{C}^{*}
\end{cases}
$$
For a specific equation, an element $\lambda$ of the factor
$\mathbf{C}^{*}$ acts by $\lambda^{\mu}$ on a block of slope
$\mu$ ; in particular, its action is trivial for the slope $0$.


\subsubsection*{C.2.3~~~Devissage and Stokes operators}

We emphasize the following interesting consequence of
Adam's lemma. Since the first slope (say, at $0$) admits
a full complement of convergent solutions, any irregular
system admits a devissage by tamely irregular equations. \\

Now, suppose one has found a local fibre functor that
restricts to $\omega^{(0)}$ on tamely irregular equations.
Then there is a restriction morphism of the local Galois
groups: $G_{i,1}^{(0)} \rightarrow G_{t,1}^{(0)}$.
It is natural to consider as the Stokes operators caused
by the irregularity at $0$ the elements killed by this
morphism. But it follows from the above devissage that the
kernel of this morphism is unipotent. \\

More generally, in the local category, every irregular 
object comes with a canonical filtration by the slopes
such that the quotients are pure. We can therefore attach
to it its "graded module", a tamely irregular object.
The fact that a slope is a quotient of a degree by a rank
suggests a strong analogy with the filtration of a vector
bundle by semi-stable bundles, although, strangely, the
slopes are here in the reverse order.


\subsection*{C.3~~~Examples}

We just want to show that tamely irregular equations 
spontaneously appear in nature, thereby motivating
a special study. We shall not treat the examples in
great detail.

\subsubsection*{C.3.1~~~An example arising from the bispectral problem}

This example was communicated to us by Emil Horozov. 
It comes from the theory of bispectral operators, 
\emph{see} \cite{BHY}, appendix A.
It arises when one looks for fixpoints of the so-called
``$q$-deformed generalized Bessel operators''. In our
notations, the equation is:
$$
\left[{\left(\sigma_{q} - q^{\beta_{1}}\right) \cdots 
\left(\sigma_{q} - q^{\beta_{N}}\right) - z^{N}}\right]f = 0.
$$
Note that it contains as a particular case the equation
satisfied by the $q$-Bessel function $J_{\nu}^{(1)}(x;q)$
studied by Changgui Zhang in \cite{ZhangBessel} 
(see equation 0.3 there). \\

The equation is clearly fuchsian at $0$, with exponents
the $q^{\beta_{i}}$, so that its local Galois group at $0$,
depends only on the multiplicative relations of these
exponents modulo $q$. \\

At $\infty$, it is pure of slope $1$, hence tamely irregular.
Writing 
$$
\left(T - q^{\beta_{1}}\right) \cdots \left(T -q^{\beta_{N}}\right) =
\sum_{i = 0}^{N} \alpha_{i} T^{N-i}
$$
the characteristic equation, putting $w = \frac{1}{z}$,
$f(z) = g(w) = \Theta_{q}(w) h(w)$, one has the equation:
$$
h(q^{N}w) + \sum_{k = 0}^{N-1} 
\frac{(-1)^{k-N} q^{\frac{k(k+1) - N(N+1)}{2}} \alpha_{k} w^{k}}
     {\alpha_{N} w^{N} - q^{-N^{2}}} h(q^{k} w) = 0.
$$
This is fuchsian, with exponents the 
$(-q^{\frac{N+1}{2}} e^{\frac{2 \imath \pi k}{N}}) , 0 \leq k \leq N-1$.
Thus, the local Galois group at $\infty$ is made up of 
a $\mathbf{C}^{*}$ factor another factor, namely 
the fuchsian local Galois group of the latter equation.

\subsubsection*{C.3.2~~~An example from $q$-integrability theory}

This example was communicated to us by Alfred Ramani. 
It comes from the work of Grammatikos and Ramani on
$q$-integrability, in theoretical physics. \\

One linearises the homographic equation 
$\sigma_{q} f = \frac{f - az}{c(f-c)}, ac \not= 0$ 
by putting $f = \frac{P}{Q}$, where:
$$
\sigma_{q} \begin{pmatrix} P \\ Q \end{pmatrix} = 
A \begin{pmatrix} P \\ Q \end{pmatrix}, \quad
\text{where:} \quad
A(z) = h(z) \begin{pmatrix} 1 & -az \\ c & -c^{2} \end{pmatrix}.
$$
The choice of $h$ is arbitrary. It can change the fuchsianity
of the system at $0$ or $\infty$, but it cannot change the
tameness: indeed, the presence of $h$ amounts to multiplying
solutions by a function $H$ such that $\sigma_{q} H = h H$,
and this does not produce divergent series. \\

We take $h = 1$, whence a system fuchsian at $0$.
The equation at $\infty$ is given by the matrix:
$$
A\left(\frac{1}{qw}\right)^{-1} =
\frac{q}{ac - q c^{2}w}
\begin{pmatrix} - c^{2} w & \frac{a}{q} \\ -c w & w \end{pmatrix}.
$$
Again, the functional factor has no effect on tameness
and we get rid of it. \\

Now, for any system of rank $2$:
$\begin{pmatrix} \alpha & \beta \\ \gamma & \delta \end{pmatrix}$,
an equivalent equation of order $2$ is:
$$
\sigma_{q}^{2}F - 
\left(\sigma_{q}(\alpha) + \delta \frac{\sigma_{q}(\beta)}{\beta}\right)
\sigma_{q} F + 
\frac{\sigma_{q}(\beta)}{\beta} (\alpha \delta - \beta \gamma) F = 0.
$$
We therefore find the equation:
$$
F(q^{2}w) + (q c^{2} - 1) w F(qw) + 
c w \left(- c w + \frac{a}{q}\right) F(w) = 0.
$$
This is pure of slope $\frac{1}{2}$, 
hence tamely irregular of level 2.


\emph{See also:}
\verb+http://www.ups-tlse.fr/~sauloy/+



\begin{thebibliography}{37}

\bibitem{YA2} \textbf{Andr\'e Y., 2001.}
Diff\'erentielles non-commutatives et th\'eorie de Galois
diff\'erentielle ou aux diff\'erences,
\emph{Ann. Scient. Ec. Norm. Sup}, no 34, 685-739.

\bibitem{Arnold} \textbf{Arnold V.I., 1980.}
\emph{Ordinary Differential Equations}, in \emph{Dynamical Systems},
\emph{Encyclopaedia of Mathematical Sciences}, Vol. 1, Springer Verlag.

\bibitem{BG} \textbf{Baranovsky V. and Ginzburg V., 1996.}
Conjugacy Classes in Loop Groups and $G$-Bundles
on Elliptic Curves,
International Mathematics Research Notes, no 15.

\bibitem{Bertrand} \textbf{Bertrand D., 1986.}
\emph{Groupes alg\'ebriques lin\'eaires et 
th\'eorie de Galois diff\'erentielle},
\emph{Cours de troisi\`eme cycle}, Universit\'e Paris VI.

\bibitem{Birkhoff1} \textbf{Birkhoff G.D., 1913.} 
The generalized Riemann problem for linear differential equations and the
allied problems for linear difference and $q$-difference equations,
\emph{Proc. Amer. Acad.}, 49, pp. 521-568.

\bibitem{Borel} \textbf{Borel A., 1991.}
\emph{Linear Algebraic Groups}, $2^{\text{nd}}$ edition, Springer Verlag.

\bibitem{CanoRamis} \textbf{Cano J. et Ramis J.-P., 1999.}
\emph{Th\'eorie de Galois diff\'erentielle},              
\emph{Book in preparation}.             

\bibitem{ChevalleyAlg} \textbf{Chevalley C., 1963.} 
\emph{Introduction to the theory of algebraic functions of one variable}, 
Mathematical Surveys, No. VI. American Mathematical Society, Providence, R.I.

\bibitem{Cohen} \textbf{Cohen R., 1965.}
\emph{Difference Algebra}, Interscience Press.

\bibitem{Deligne} \textbf{Deligne P., 1970.}                     
\emph{Equations diff\'erentielles \`a points singuliers r\'eguliers},
\emph{Lecture Notes in Mathematics}, 163, Springer Verlag.

\bibitem{DF} \textbf{Deligne P., 1990.}
Cat\'egories Tannakiennes, 
in \emph{Grothendieck Festschrift} (Cartier \& al. eds),
Vol. II, Birkh\"{a}user.

\bibitem{DM} \textbf{Deligne P. and Milne J., 1989.}
\emph{Tannakian Categories},
in \emph{Hodge Cycles, Motives and Shimura Varieties} (Deligne \& al. eds),
\emph{Lecture Notes in Mathematics}, 900, Springer Verlag.

\bibitem{LDVpreprint} \textbf{Di Vizio L., 2000.}
Arithmetic theory of $q$-difference equations. 
The $q$-analogue of Grothendieck-Katz conjecture
on $p$-curvatures. 
Pr\'epublication de l'Institut de Math\'ematiques
de Jussieu, no 286. \emph{Also} to appear in 
\emph{Invent. Math.}

\bibitem{Etingof} \textbf{Etingof P.I., 1995.}
Galois Groups and Connection Matrices of $q$-difference Equations,
\emph{Electronic Research Announcements of the A.M.S.}, Vol. 1, Issue 1.

\bibitem{GR} \textbf{Gasper G. and Rahman M., 1990.}
\emph{Basic hypergeometric series},
\emph{Encyclopedia of Mathematics}, Vol. 35, Cambridge University Press.

\bibitem{Hendriks} \textbf{Hendriks P.A., 1996.}
Algebraic aspects of linear differential and difference equations, 
\emph{Thesis}, University of Groningen.

\bibitem{Ince} \textbf{Ince E.L., 1956.}
\emph{Ordinary Differential Equations}, Dover Publications.

\bibitem{Katz} \textbf{Katz N.M., 1987.}
On the calculation of some differential galois groups,
Invent. math. 87, 13-61.

\bibitem{Praagman} \textbf{Praagman C., 1983.}
The formal classification of linear difference equations,
\emph{Proc. Kon. Ned. Ac. Wet. ser. a}, 86.

\bibitem{SVdP} \textbf{van der Put M. and Singer M.F., 1997.}
\emph{Galois theory of difference equations},
\emph{Lecture Notes in Mathematics}, 1666, Springer Verlag.

\bibitem{Ramanujan} \textbf{Ramanujan S., 1927.}
\emph{Collected Works}, Chelsea.

\bibitem{RamisGrowth} \textbf{Ramis J.-P., 1992.}
About the growth of entire functions solutions to linear algebraic 
$q$-difference equations,
\emph{Annales de Fac. des Sciences de Toulouse}, S\'erie 6, Vol. I, no 1,
pp. 53-94.

\bibitem{RamisJPRTraum} \textbf{Ramis J.-P., 1990.}
Fonctions $\theta$ et \'equations aux $q$-diff\'erences, 
\emph{Unpublished notes}, Strasbourg.

\bibitem{RamisGalDiff} \textbf{Ramis J.-P., 1999.}
\emph{About the Inverse Problem in Differential Galois Theory.
The Differential Abhyankar Conjecture},
\emph{Book in preparation}. A summary appears in:
\textbf{Braaksma B.L.J., Immink G.K. and van der Put M., 1996.}
\emph{The Stokes Phenomenon and Hilbert's 16th Problem},
World Scientific.

\bibitem{RSZ} \textbf{Ramis J.-P., Sauloy J. and Zhang C., 2001.}
Local analytic classification of irregular $q$-difference equations,
\emph{Article in preparation}.

\bibitem{JSthese} \textbf{Sauloy J., 1999.}
Th\'eorie de Galois des \'equations aux $q$-diff\'erences fuchsiennes,
\emph{Th\`ese}, Universit\'e Paul Sabatier, Toulouse.

\bibitem{JSAIF} \textbf{Sauloy J., 2000.}
Syst\`emes aux $q$-diff\'erences singuliers r\'eguliers :
classification, matrice de connexion et monodromie,
\emph{Annales de l'Institut Fourier}, 
Tome 50, fasc. 4, pp. 1021-1071.

\bibitem{JSCRAS3} \textbf{Sauloy J., 2002.}
La filtration canonique par les pentes d'un module
aux $q$-diff\'erences, 
C.R. Acad. Sci. Paris, janvier 2002.

\bibitem{JSAIF2} \textbf{Sauloy J., 2002.}
La filtration canonique par les pentes des modules aux $q$-diff\'erences
et le gradu\'e associ\'e,
\emph{Article in preparation}.

\bibitem{JSIRR} \textbf{Sauloy J., 2002.}
Local Galois theory of irregular $q$-difference equations,
\emph{Article in preparation}.

\bibitem{JSGAL} \textbf{Sauloy J., 2002.}
Galois Theory of fuchsian $q$-difference equations,
\emph{Long version of the present paper},
\verb^url: picard.ups-tlse.fr/~sauloy^.

\bibitem{JSGTQDIF} \textbf{Sauloy J., 2001.}
La filtration canonique par les pentes des modules aux $q$-diff\'erences
et le gradu\'e associ\'e,
R\'edaction d'expos\'es au Groupe de Travail 
``Equations aux $q$-diff\'erences'', 
\verb^url: picard.ups-tlse.fr/~sauloy^.

\bibitem{SerreGacc} \textbf{Serre J.-P., 1959.}
\emph{Groupes alg\'ebriques et corps de classes},
Hermann.

\bibitem{Seshadri} \textbf{Seshadri C.S., 1982.}
\emph{Fibr\'es vectoriels sur les courbes alg\'ebriques},
Ast\'erisque No. 96, Soci\'et\'e Math\'ematique de France.

\bibitem{Springer} \textbf{Springer T.A., 1998.}
\emph{Linear Algebraic Groups, Second Edition},
Birkh\"{a}user.

\bibitem{Wasow} \textbf{Wasow W., 1965.}
\emph{Asymptotic Expansions for Ordinary Differential Equations},
Dover Publications.

\bibitem{Weil} \textbf{Weil A., 1938.}
G\'en\'eralisation des fonctions ab\'eliennes,
J. Math. Pures et Appl., no 17, pp. 47-87.

\bibitem{ZhangGroningen} \textbf{Zhang C., 2001.}
Une sommation discr\`ete pour des \'equations aux $q$-diff\'erences
lin\'eaires et \`a coefficients analytiques: th\'eorie g\'en\'erale
et exemples,
to appear in the Proceedings of the Workshop
``Differential Equations and Stokes Phenomenon'', Groningen, may 2001,
also as a preprint of the Universit\'e Paul Sabatier, Toulouse.

\end{thebibliography}

\begin{thebibliography}{54}

\bibitem{BHY} \textbf{Bakalov B., Horozov E. and Yakimov M., 1996.}
General methods for constructing bispectral operators,
Physics Letters A 222, pp. 59-66.

\bibitem{ChevalleyLie} \textbf{Chevalley C., 1957.} 
\emph{Theory of Lie groups, I}, Princeton.

\bibitem{DF} \textbf{Deligne P., 1990.}
Cat\'egories Tannakiennes, 
in \emph{Grothendieck Festschrift} (Cartier \& al. eds),
Vol. II, Birkh\"{a}user.

\bibitem{DM} \textbf{Deligne P. and Milne J., 1989.}
\emph{Tannnakian Categories},
in \emph{Hodge Cycles, Motives and Shimura Varieties} (Deligne \& al. eds),
\emph{Lecture Notes in Mathematics}, 900, Springer Verlag.

\bibitem{SVdP} \textbf{van der Put M. and Singer M.F., 1997.}
\emph{Galois theory of difference equations},
\emph{Lecture Notes in Mathematics}, 1666, Springer Verlag.

\bibitem{RamisGalDiff} \textbf{Ramis J.-P., 1999.}
\emph{About the Inverse Problem in Differential Galois Theory.
The Differential Abhyankar Conjecture},
\emph{A para\^itre}.

\bibitem{JSthese} \textbf{Sauloy J., 1999.}
Th\'eorie de Galois des \'equations aux $q$-diff\'erences fuchsiennes,
\emph{Th\`ese}, Universit\'e Paul Sabatier, Toulouse.

\bibitem{JSAIF} \textbf{Sauloy J., 2000.}
Syst\`emes aux $q$-diff\'erences singuliers r\'eguliers:
classification, matrice de connexion et monodromie,
\emph{Annales de l'Institut Fourier}, 
Tome 50, fasc. 4, pp. 1021-1071.

\bibitem{Silverman} \textbf{Silverman J.H., 1985.}
\emph{The Arithmetic of Elliptic Curves}, Springer Verlag.

\bibitem{Toen} \textbf{Toen B., 2000.}
\emph{Dualit\'e de Tannaka sup\'erieure I:
Structures monoidales}, 
Preprint of the Max Planck Institute,
MPI-2000-57, \emph{see also:}
\verb+http://www.mpim-bonn.mpg.de+.

\bibitem{ZhangBessel} \textbf{Zhang C., 2000.}
Sur les fonctions $q$-Bessel de Jackson,
Preprint, Universit\'e de La Rochelle.

\end{thebibliography}
\end{document}